\documentclass[12pt]{memo-l}
\usepackage{graphicx,amssymb,latexsym,eepic,bm,hyphenat,mathtools}
\usepackage[dvips,includehead,width=14.85cm,height=24.5cm,twoside,centering]{geometry}
\usepackage{color}
\usepackage[numbers,sort&compress]{natbib}

\usepackage[bookmarks=true, bookmarksopen=false,%
    colorlinks=true,%
    linkcolor=darkblue,%
    citecolor=darkblue,%
    filecolor=darkblue,%
    menucolor=darkblue,%
    urlcolor=darkblue
]{hyperref}

\usepackage{tikz,tikz-cd}

\definecolor{darkred}{rgb}{1,0,0}
\definecolor{lightred}{rgb}{1,0.4,0}
\definecolor{darkblue}{cmyk}{1,0.4,0,0.4}  
\definecolor{lightblue}{cmyk}{1,0.4,0,0}  
\definecolor{darkgreen}{cmyk}{1,0.5,1,0}
\definecolor{lightgreen}{cmyk}{1,0,1,0}
\newcommand{\darkred}[1]{\color{darkred}{#1}\color{black}}
\newcommand{\lightred}[1]{\color{lightred}{#1}\color{black}}
\newcommand{\darkblue}[1]{\color{darkblue}{#1}\color{black}}
\newcommand{\lightblue}[1]{\color{lightblue}{#1}\color{black}}
\newcommand{\darkgreen}[1]{\color{darkgreen}{#1}\color{black}}
\newcommand{\lightgreen}[1]{\color{lightgreen}{#1}\color{black}}

\def\mathcenter#1{\vcenter{\hbox{$#1$}}}

\def\mfig#1{\mathcenter{\includegraphics{#1}}}
\def\mfigb#1{\mathcenter{\includegraphics[trim=-1 -1 -1 -1]{#1}}}

\newcommand{\arxiv}[1]{\href{http://arxiv.org/abs/#1}{\nolinkurl{arXiv:#1}}}
\numberwithin{section}{chapter}
\numberwithin{equation}{chapter}

\usepackage{chngcntr}
\counterwithout{figure}{chapter}

\newtheorem{theorem}{Theorem}[chapter]
\newtheorem{proposition}[theorem]{Proposition}

\newtheorem{corollary}[theorem]{Corollary}
\newtheorem{lemma}[theorem]{Lemma}

\theoremstyle{definition}

\newtheorem{remark}[theorem]{Remark}
\newtheorem{question}[theorem]{Question}
\newtheorem{example}[theorem]{Example}
\newtheorem{definition}[theorem]{Definition}

\newcommand{\abs}[1]{\lvert#1\rvert}
\newcommand{\Arcs}{\mathbf{A}}
\newcommand{\Exch}{\mathbf{E}}
\newcommand{\Mark}{\mathbf{M}}

\newcommand{\notch}{\scriptstyle\bowtie}

\newcommand{\Surf}{\mathbf{S}}

\newcommand{\SM}{(\Surf,\Mark)}
\newcommand{\barSM}{(\overline{\Surf,\Mark})}

\newcommand{\barM}{\overline\Mark}
\newcommand{\barL}{{\overline{L}}}
\newcommand{\barLL}{{\overline{\LL}}}

\newcommand{\Teich}{\mathcal{T}}
\newcommand{\dTeich}{\widetilde{\Teich}}
\newcommand{\dPunct}{\tilde{P}}

\newcommand{\AP}{\mathbf{A}^{\hspace{-1pt}\circ}}
\newcommand{\AT}{\mathbf{A}^{\hspace{-1pt}\notch}}
\newcommand{\EP}{\mathbf{E}^\circ}
\newcommand{\ET}{\mathbf{E}^{\notch}}
\newcommand{\DP}{\Delta^{\hspace{-1pt}\circ}}
\newcommand{\DT}{\Delta^{\hspace{-1pt}\notch}}

\newcommand{\DTSM}{\DT(\Surf,\Mark)}
\newcommand{\ESM}{\Exch(\Surf,\Mark)}
\newcommand{\ETSM}{\ET(\Surf,\Mark)}

\newcommand{\DPSM}{\DP\SM}
\newcommand{\EPSM}{\EP\SM}
\newcommand{\APSM}{\AP\SM}
\newcommand{\ATSM}{\AT(\Surf,\Mark)}

\newcommand{\barAPSM}{\AP\barSM}
\newcommand{\barATSM}{\AT\barSM}
\newcommand{\BSM}{\mathbf{B}\SM}
\newcommand{\TTSM}{\overline{\Teich}\SM}
\newcommand{\TTSML}{\overline{\Teich}(\Surf,\Mark,\LL)}
\newcommand{\baralpha}{\overline{\alpha}}
\newcommand{\barbeta}{\overline{\beta}}
\newcommand{\bargamma}{\overline{\gamma}}
\newcommand{\bardelta}{\overline{\delta}}
\newcommand{\bareta}{\overline{\eta}}
\newcommand{\bartheta}{\overline{\theta}}

\newcommand{\Acal}{\mathcal{A}}

\newcommand{\Fcal}{\mathcal{F}}
\newcommand{\Rcal}{\mathcal{R}}

\newcommand{\ExPat}{\mathcal{M}}
\newcommand{\ZP}{\ZZ\PP}

\newcommand{\CC}{\mathbb{C}}
\newcommand{\PP}{\mathbb{P}}
\newcommand{\ZZ}{\mathbb{Z}}
\newcommand{\RR}{\mathbb{R}}
\newcommand{\HH}{\mathbb{H}}
\newcommand{\LL}{\mathbf{L}}
\newcommand{\Trop}{{\rm Trop}}

\newcommand{\Star}{{\operatorname{star}}}
\newcommand{\oplusk}{{\,\oplus_k\,}}

\newcommand{\atag}{\scriptstyle\bowtie}

\DeclareMathOperator{\Diff}{Diff}
\DeclareMathOperator{\Gr}{Gr} 

\def\xx{\mathbf{x}}
\def\yy{\mathbf{y}}
\def\zz{\mathbf{z}}
\def\kk{\mathbf{k}}
\def\pp{\mathbf{p}}

\begin{document}

\frontmatter

\title{Cluster~algebras~and~triangulated~surfaces
Part~II:~Lambda~lengths}

\author{Sergey Fomin}
\address{
Department of Mathematics, University of Michigan,
Ann Arbor, MI 48109, USA}
\email{fomin@umich.edu}

\author{Dylan Thurston}
\address{
Department of Mathematics, Indiana University,
  Bloomington, IN 47405, USA}
\email{dpthurst@indiana.edu}
\thanks{Partially supported by NSF grants DMS-0555880 (S.~F.), 
DMS-1101152 (S.~F.), DMS-1361789 (S.~F.), DMS-1008049 (D.~T.), and DMS-1507244 (D.~T.),
and a Sloan Fellowship~(D.~T.).
}

\date{
October 19, 2012. Revised March 15, 2016}



\begin{abstract}
For any cluster algebra
whose underlying combinatorial data can be encoded by a
bordered surface with marked
points, we construct a geometric realization
in terms of
suitable decorated Teichm\"uller space of the surface.
On the geometric side, this requires opening the surface at each
interior marked
point~into~an additional geodesic boundary component.
On the algebraic side, it relies on the notion of a non-normalized
cluster algebra and the machinery of tropical lambda lengths.

Our model allows for an arbitrary choice of coefficients
which translates into~a choice of a family of integral laminations 
on the surface.
It provides an intrinsic \linebreak[3] 
interpretation of cluster variables as
renormalized lambda lengths of arcs on the surface.
Exchange relations 
are written in terms of the shear coordinates of the laminations, and
are interpreted as generalized Ptolemy relations for lambda lengths.

This approach gives alternative proofs for the main structural results
from our previous paper, removing unnecessary assumptions on the
surface.
\end{abstract}

\maketitle

\setcounter{tocdepth}{0}

\tableofcontents

\mainmatter

\chapter{Introduction}

\section*{Overview of the main results}

This paper continues the study
of cluster algebras associated with marked Riemann surfaces with holes
and punctures.
The emphasis of our first paper~\cite{cats1}, written 
in collaboration
with Michael Shapiro, was on the combinatorial construction
of the \linebreak[3]
cluster complex, which we showed to be closely related to the
\emph{tagged arc complex}~of 
the surface, a particular extension of
the classical simplicial complex of arcs connect\-ing marked points.
The focus of the current paper is on geometry, specifically on 
providing an explicit geometric interpretation for the cluster
variables in \emph{any} cluster 
algebra (of geometric type)
whose exchange matrix can be associated with a~tri\-angulated surface.
More concretely, we demonstrate that each cluster variable can be
viewed as a properly
normalized \emph{lambda length} of the corresponding (tagged)~arc. 
Introduced by
R.~Penner \cite{penner-decorated, penner-bordered, penner-lambda},
lambda lengths serve as coordinates on appropriate 
\emph{decorated Teichm\"uller spaces};
a point in such a space is a hyperbolic metric on a
surface, together with some additional decoration. 
In fact, we use an extension of lambda lengths to \emph{opened
  surfaces}, surfaces with some extra geodesic boundary.

The main underlying idea of this line of inquiry---already present in
the pioneering works of V.~Fock and
A.~Goncharov~\cite{fock-goncharov1, fock-goncharov2},
and of M.~Gekhtman, M.~Shapiro, and A.~Vainshtein~\cite{gsv2}---is to
interpret a decorated Teichm\"uller space as the real positive part of an
algebraic variety. The coordinate ring $\mathcal{A}$ of this variety is
generated by the variables corresponding to the  lambda lengths;
these variables satisfy certain algebraic relations among lambda
lengths, known as generalized \emph{Ptolemy relations}.
The rings $\mathcal{A}$ arising from various versions of this
construction possess, in a very natural way, a cluster algebra
structure: the lambda lengths become cluster variables,
and clusters correspond to tagged triangulations.

\smallskip

The main results obtained in this paper can be succinctly summarized as
follows.
\begin{enumerate}
\item
We investigate a broad class of cluster algebras that includes any
cluster algebra of geometric type whose (skew-symmetric) exchange
matrix is a signed adjacency matrix of a triangulated bordered
surface. Each such cluster algebra is naturally associated with a
collection of integral \emph{laminations} on the surface.
Every exchange relation is then readily written in terms of
\emph{shear coordinates} of these laminations with respect to a given
tagged triangulation.

\item
We describe the cluster variables in these cluster algebras as
generalized lambda lengths. To this end, we (a) extend the lambda
length construction to tagged arcs, (b) define Teichm\"uller spaces of
opened surfaces and the associated lambda lengths, (c) define
\emph{laminated Teichm\"uller spaces} and \emph{tropical lambda
  lengths}, and (d) combine all these elements together to realize
cluster variables as certain rescalings of lambda lengths of tagged~arcs.
\pagebreak[3]

\item
The underlying combinatorics of a cluster algebra is governed by its
\emph{cluster complex}.
In our first paper~\cite{cats1}, we described this complex in terms of
tagged arcs, and investigated its basic structural properties.
For technical reasons, some results 
in~\cite{cats1} required an exception that excluded
closed surfaces with exactly two punctures.
In this paper, we remove this restriction, extending the main results
of~\cite{cats1} to
arbitrary bordered surfaces with marked points.
\end{enumerate}

\section*{Structure of the paper}

Chapters~\ref{sec:non-normalized}--\ref{sec:geom-type} are devoted to
preliminaries on cluster algebras and exchange patterns.
It turns out that the proper algebraic framework for our main
construction is provided by the axiomatic setting of
\emph{non-normalized} cluster algebras. 
This setting, although the original one 
for the cluster algebra theory~\cite{ca1}, 
was all but abandoned in the intervening years, as the main 
developments and applications dealt almost
exclusively with normalized cluster algebras.
Chapter~\ref{sec:non-normalized} contains a review of non-normalized 
cluster algebras and the mutation rules used to define them.
The only (minor) novelty here is Proposition~\ref{pr:yhat}. 
Chapter~\ref{sec:rescaling} discusses rescaling of cluster variables
and the technicalities involved in constructing a normalized pattern
by rescaling a non-normalized one.
Chapter~\ref{sec:geom-type} recalls the notion of cluster algebras of
geometric type, and discusses realizations of these algebras in which
both cluster and coefficient variables are represented by positive
functions on a topological space, in anticipation of
Teichm\"uller-theoretic applications.


In Chapter~\ref{sec:cats1-recap}, we review the main constructions and
some of the main results of the prequel~\cite{cats1} to this paper:
ordinary and tagged arcs on a bordered surface with marked points;
triangulations, flips, and arc complexes (both ordinary and tagged);
associated exchange graphs; and signed adjacency matrices and their
mutations.


In Chapter~\ref{sec:structural}, we formulate the first batch of our
results, which describe the structure of cluster algebras whose
exchange matrices come from triangulated surfaces.
(The proofs come much later.)
In brief, we extend all main structural theorems of~\cite{cats1} to arbitrary
surfaces, including closed surfaces with two marked points,
not covered in~\cite{cats1} because of the exceptional nature of the
fundamental groups of the corresponding exchange graphs.


In Chapters~\ref{sec:hyperbolic-cluster}--\ref{sec:non-normalized-patterns-from-surfaces},
we develop the hyperbolic geometry tools required for our main
construction.
Chapter~\ref{sec:hyperbolic-cluster} presents Penner's concept of
lambda length, and its basic algebraic properties.
In Chapter~\ref{sec:lambda-tagged}, we adapt this concept to the
tagged setting, and write down the appropriate versions of Ptolemy
relations. This enables us to describe, in
Theorem~\ref{th:cluster-lambda}, the main structural features of
a particular class of exchange patterns whose coefficient variables
come from boundary segments.
(This class of cluster algebras already appeared in
\cite{fock-goncharov1, fock-goncharov2, gsv2}.)
To handle general coefficient systems, we need another geometric
idea, introduced in Chapter~\ref{sec:opened-surfaces}: the concept of an
\emph{opened surface} obtained by replacing interior marked points by
geodesic circular boundary components.
Chapter~\ref{sec:lambda-opened} presents the corresponding version of
lambda lengths (for the \emph{lifted arcs} on the opened surface),
and an appropriate variation of the decorated Teichm\"uller space.
In Chapter~\ref{sec:non-normalized-patterns-from-surfaces}, we show
that suitably rescaled lambda lengths of lifted arcs
form a non-normalized exchange pattern, providing a crucial building
block for our main construction (to be completed in
Chapter~\ref{sec:laminated-teich}).


Chapters~\ref{sec:shear}--\ref{sec:tropical-lambda} are devoted to
the combinatorics of general coefficient systems of geometric type
(equivalently, extended exchange matrices).
As noted by Fock and Goncharov~\cite{fg-dual-teich},
W.~Thurston's shear coordinates for simple closed curves transform
under ordinary flips
according to the rules of matrix mutations.
In Chapter~\ref{sec:shear}, we review this beautiful theory;
in Chapter~\ref{sec:shear-tagged},
we extend its main results, namely Thurston's coordinatization
theorem for integral laminations and the matrix mutation rule,
to the tagged setting.
Chapter~\ref{sec:tropical-lambda} discusses the notion of
\emph{tropical lambda length}, a discrete analogue of Penner's
concept, in which hyperbolic lengths are replaced by a combination of
transverse measures with respect to a family of laminations. 
Tropical lambda lengths satisfy the tropical version of exchange
relations, which makes them ideally suited for the role of
rescaling factors in our main construction.


This construction is presented in
Chapter~\ref{sec:laminated-teich}. Here we introduce the notion of a
\emph{laminated Teichm\"uller space} whose defining data include, in
addition to the surface, a fixed multi-lamination on~it.
This space can be coordinatized by yet another version of lambda
lengths, obtained by dividing the ordinary lambda lengths of lifted
arcs by the tropical ones. (The ratio does not depend
on the choice of a lift of the arc. Although it does depend on
the choice of lifted laminations, this choice does not affect the
resulting cluster algebra structure, up to a canonical isomorphism.)
These ``laminated lambda lengths'' 
form an exchange pattern (thus generate a cluster algebra) of the
required kind. In other words, they satisfy the exchange relations
of Ptolemy type whose coefficients are encoded by the shear
coordinates of the chosen multi-lamination with respect to the current
tagged triangulation. This geometric realization allows us to prove
all our claims, made in earlier chapters, pertaining to the structural
properties of the cluster algebras under consideration.


The next two chapters are devoted to applications and examples.
In Chapter~\ref{sec:topol-coord-ring}, we provide topological
models for cluster-algebraic structures in coordinate rings of various
algebraic varieties, such as certain Grassmannians and affine
base spaces. In each case, the relevant cluster structure is encoded
by a particular choice of a bordered surface with a
collection of integral laminations on it. These examples provide 
fascinating links---awaiting further exploration---between representation
theory
and combinatorial topology.
Chapter~\ref{sec:prin-univ-coeff} treats in concrete detail two
important classes of coefficient systems introduced in~\cite{ca4}: the
\emph{principal coefficients} and the \emph{universal coefficients}
(the latter in finite types $A_n$ and~$D_n$ only).


Appendix~\ref{app:relative-lambda-lengths} contains an informal
discussion of how our main geometric construction of renormalized
lambda lengths on opened surfaces can be obtained by means of
\emph{tropical degeneration} from a somewhat more conventional
setting in which, instead of fixing a multi-lamination,
one picks in advance a decorated hyperbolic structure on an opened
surface.


Appendix~\ref{app:vers-teichmuller} is designed to
help the reader navigate between the various versions of Teichm\"uller
spaces and their respective coordinatizations. 

\smallskip

Several figures in this paper are best viewed in \darkred{color}. 

\smallskip

Our notation and terminology agree with the previous paper in the
series~\cite{cats1}.

\newpage


\section*{Historical notes}

A preliminary version of this paper (64 pages long) was circulated in
May 2008, and posted on our respective web sites. 
While that version contained all the main results and proofs, we were
not satisfied with the exposition and felt the need to double-check
the details of our setup.  Indeed, careful inspection
revealed a number of
flaws, each of them fixable. In particular, the factor $\nu(p)$ in
Definition~\ref{def:lambda-opened-tagged} used to be different, 
which meant
that the older version of Lemma~\ref{lem:exch-opened-digon} did not
work uniformly for plain and notched arcs.  Similarly,
Definition~\ref{def:transverse-measure} for notched arcs was missing
the term $\abs{l_\barL(p)}$, which meant that
Lemma~\ref{lem:exch-digon-trop} was not uniform, causing problems
for the cluster algebra structure.

%
Completing the revision 
took us much longer than originally anticipated: 
the revised preprint~\cite{cats2-arxiv-v1} was only posted
in October 2012. 
We thank the friends and colleagues who kept up the
pressure, urging us to finish the~job. 
This version of the paper is essentially identical
to~\cite{cats2-arxiv-v1}, 
save for minor editorial changes,
  additional historical comments, and expanded bibliography. 
We did not make any attempts to change the presentation to reflect 
related developments in cluster 
algebra theory that took place after the preliminary version 
was circulated in~2008. In~particular, we do not cite the papers
referencing this work, except in the brief historical notes included
below. 
One advantage of this decision is that the reader need not worry that we might 
rely on some work that was in turn dependent on our results. 

We kept the citations in the text to R.~Penner's lecture
notes~\cite{penner-lambda}, 
even though the latter were subsumed by his recent book~\cite{penner-book}. 
We point out that \cite[Lemma~4.4]{penner-lambda} (our
Lemma~\ref{lem:horocyclic-segment}) matches 
\cite[Chapter~1, Lemma~4.9]{penner-book},
whereas \cite[Theorem~5.10]{penner-lambda} (cf.\ our Remark~\ref{rem:penner-5.10})
has become \cite[Chapter~2, Theorem~2.25]{penner-book}. 

\textbf{Subsequent developments.}
At referee's request, we include a brief description of some related
research on cluster algebras associated with surfaces,
focusing on results obtained before the release of the \texttt{arXiv} 
version of this paper~\cite{cats2-arxiv-v1}. 

From the standpoint of  general structure theory of cluster algebras, 
the importance of the surface case studied in this paper and its
prequel~\cite{cats1} was validated by the classification of cluster
algebras of \emph{finite mutation type} obtained
in~\cite{felikson-shapiro-tumarkin, felikson-shapiro-tumarkin-2}. 
In the skew-symmetric case, 
the classification~\cite{felikson-shapiro-tumarkin} states that  
a con\-nected quiver~$Q$ without frozen vertices has finite mutation type
if and only if $Q$ arises from a triangulation of
a marked bordered surface---unless $Q$ has at most two
vertices or is isomorphic to one of 
11~exceptional quivers listed in~\cite{derksen-owen}.
Another closely related result is a complete classification of cluster algebras of
polynomial growth~\cite{felikson-shapiro-tumarkin-growth}. 

Manifestly positive combinatorial formulas for 
Laurent expansions of cluster variables in cluster algebras associated
with surfaces 
were given in \cite{musiker-schiffler, schiffler}
(for surfaces without punctures) 
and~\cite{musiker-schiffler-williams-laurent-positivity}
(in full generality). 
Even more general formulas of this kind,
for elements of certain explicit combinatorial bases in these cluster algebras, 
were given in~\cite{musiker-schiffler-williams-bases} 
(making use of~\cite{musiker-williams}). 
Another family of bases, 
closely related to dual canonical bases for quantum groups, 
was investigated in~\cite{d.thurston-pnas}. 

A construction of quivers with potentials associated to
triangulations was proposed in
\cite{qp-triangulated-1}, and further studied
in~\cite{qp-triangulated-3}. 

A geometric construction that can be regarded as the orbifold version of 
lambda lengths was proposed in~\cite{shapiro-chekhov}. 
This required introducing a 
certain generalization of the concept of a cluster algebra. 

An interested reader may also wish to consult  
the surveys~\cite{icm, williams-survey}, 
the book~\cite{marsh-book}, and the \emph{Cluster Algebras
  Portal}~\cite{portal} 
for additional references and resources.

\pagebreak[3]

\section*{Acknowledgments}

First and foremost, we would like to thank Michael Shapiro for the
collaboration~\cite{cats1} that inspired our work on this paper, for
valuable insights, and for comments on the preliminary version. 
We are also grateful to 
Leonid Chekhov, 
Yakov Eliashberg, 
Vladimir Fock, 
Alexander Goncharov, 
David Kazhdan, 
Daniel Labardini Fragoso, 
Bernard Leclerc,
Gregg Musiker, 
Robert Penner, 
Nathan Reading,
Hugh Thomas, 
Pavel Tumarkin, 
Alek Vainshtein,
Lauren Williams, 
and Andrei Zelevinsky for helpful advice, stimulating discussions, 
and/or comments on the preliminary version of the paper.
We thank the anonymous referees for their encouraging comments and helpful suggestions. 

Most of the work was done while D.~T.\ was at Barnard College,
Columbia University.  
We acknowledge the hospitality of IRMA (Strasbourg), MSRI and UC Berkeley,
and QGM (\AA rhus).

\chapter{Non-normalized cluster algebras}
\label{sec:non-normalized}

The original definition and basic properties of
(non-normalized) cluster algebras were given in~\cite{ca1};
see~\cite{ca3} for further developments.
In this chapter, we recall the basic notions of this theory following
the aforementioned sources (cf.\ especially \cite[Sections~1.1--1.2]{ca3}).
Note that our previous paper~\cite{cats1} used a more restrictive
normalized setup.

The construction of a (non-normalized, skew-symmetrizable)
cluster algebra begins with a \emph{coefficient
  group}~$\PP$,
an abelian group without torsion, written multiplicatively.
Take the integer group ring~$\ZP$, and let $\Fcal$ be (isomorphic to)
a field of rational functions in $n$ independent
variables with coefficients in~$\ZP$.
Later, we are going to call the positive integer $n$ the \emph{rank} of
our yet-to-be-defined cluster algebra, and $\Fcal$ its \emph{ambient field}.
The following definitions are central for the cluster algebra theory.

\begin{definition}[\emph{Seeds}]
\label{def:seed}
A \emph{seed} in $\Fcal$ is a triple $\Sigma = (\xx, \pp, B)$
consisting of:
\begin{itemize}
\item
a \emph{cluster} $\xx\subset \Fcal$,
a set of $n$ algebraically independent elements
(called \emph{cluster variables}) which
generate $\Fcal$ over the field of fractions of~$\ZP$;
\item
a \emph{coefficient tuple} $\pp = (p_x^\pm)_{x\in \xx}$,
a $2n$-tuple of elements of~$\PP$;
\item
an \emph{exchange matrix} $B \!=\! (b_{xy})_{x,y\in\xx}$,
a \emph{skew-symmetrizable} $n\!\times\! n$ integer matrix.
\end{itemize}
That is, $B$ can be made skew-symmetric by rescaling its columns by
appropriately chosen positive integer scalars.
(In all applications in this paper, $B$~will in fact be skew-symmetric.)
\end{definition}

\begin{definition}[\emph{Seed mutations}]
\label{def:seed-mutation}
Let $\Sigma = (\xx, \pp, B)$
be a seed in $\Fcal$, as above.
Pick a cluster variable $z\in\xx$.
We~say that another seed $\overline \Sigma = (\overline \xx,\overline \pp,
\overline B)$ is related to $\Sigma$ by a
\emph{seed mutation} in direction~$z$ if
\begin{itemize}
\item
the cluster $\overline \xx$ is given by
$\overline \xx = \xx - \{z\} \cup \{\overline z\}$,
where the new cluster variable $\overline z \in \Fcal$ is determined
by the \emph{exchange relation}
\begin{equation}
\label{eq:exchange-rel-xx}
z\,\overline z
=p_z^+ \, \prod_{\substack{x\in\xx \\ b_{xz}>0}} x^{b_{xz}}
+p_z^- \, \prod_{\substack{x\in\xx \\ b_{xz}<0}} x^{-b_{xz}} \,;
\end{equation}
\item
the coefficient tuple
$\overline \pp=(\overline p ^\pm_y)_{x\in\overline \xx}\,$
satisfies
\begin{align}
\label{eq:p-mutation1}
\overline p ^\pm_{\overline z}&=p ^\mp_z\,;\\
\label{eq:p-mutation2}
\frac{\overline p ^+_y}{\overline p ^-_y} &=
\begin{cases}
(p_z^+)^{b_{zy}}\,\dfrac{p ^+_y}{p ^-_y}   & \text{if $b_{zy}\geq 0$}\\[.15in]
(p_z^-)^{b_{zy}}\,\dfrac{p ^+_y}{p ^-_y}   & \text{if $b_{zy}\leq 0$}
\end{cases}
\text{\qquad (for $y\neq z$)}
\end{align}

\item
the exchange matrix
$\overline B\!=\!(\overline b_{xy})$ is obtained from $B$ by the \emph{matrix
mutation}~rule
\begin{equation}
\label{eq:B-mutation}
\overline b_{x y} =
\begin{cases}
0 & \text{if $x=y=\overline z$};\\
- b_{xz} & \text{if $y=\overline z\neq x$};\\
- b_{zy} & \text{if $x=\overline z\neq y$};\\
b_{xy} & \text{if $x\neq\overline z$, $y\neq\overline z$, and
$b_{xz}b_{zy}\leq 0$};\\
b_{xy} + |b_{xz}|b_{zy}
& \text{if $x\neq\overline z$, $y\neq\overline z$, and
$b_{xz}b_{zy}\geq 0$.}
\end{cases}
\end{equation}
\end{itemize}
\end{definition}

It is easy to check that the mutation rule is symmetric:
$\Sigma$ is in turn related to $\overline \Sigma$ by a mutation in
direction~$\overline z$.

\begin{remark}
\label{rem:nn-vs-n}
The crucial---and only---difference between
Definition~\ref{def:seed-mutation} and its counterpart 
\cite[Definition~5.1]{cats1} used in our previous paper is that
here we do not require $\PP$ to be endowed with an ``auxiliary addition''
making it into a semifield 
(cf.\ Definition~\ref{def:n-exch-pattern}), and 
consequently eliminate the normalization requirement
(cf.~\eqref{eq:normalization}). 
As noted in~\cite{ca1},
there is a price to be paid for getting rid of the semifield
structure: in the absence of normalization,
the seed $\overline \Sigma$ in Definition~\ref{def:seed-mutation}
is not determined by $\Sigma$ and a choice of~$z$.
Indeed, the $n-1$ monomial formulas~\eqref{eq:p-mutation2}
prescribe the ratios ${\overline p ^+_y}/{\overline p ^-_y}$ only,
leaving us with $n-1$ degrees
of freedom in choosing specific coefficients~${\overline p ^\pm_y}$.
\end{remark}

\begin{remark}
The matrix mutation rule \eqref{eq:B-mutation} has many equivalent
reformulations; see, e.g., \cite[(4.3)]{ca1}, \cite[(2.2)]{ca4}, and
\cite[(2.5)]{ca4}.
In particular, the last two cases in~\eqref{eq:B-mutation} 
(where we assume that $x\neq\overline z$ and $y\neq\overline z$)
can be restated as 
\begin{align}
\label{eq:B-mutation-alt}
\overline b_{xy}
&= b_{xy}+[b_{zy}]_+ b_{xz} + b_{zy} [-b_{xz}]_+ 
\end{align} 
where we use the notation $[b]_+=\max(b,0)$.
Similarly, \eqref{eq:p-mutation2} can be rewritten as
\begin{equation}
\label{eq:p-mutation2-new}
\frac{\overline p ^+_y}{\overline p ^-_y} 
=\frac{p ^+_y}{p ^-_y} 
\Biggl(\frac{p ^+_z}{p ^-_z} 
\Biggr)^{[b_{zy}]_+}
(p ^-_z)^{b_{zy}}.
\end{equation}
\end{remark}

\begin{example}[\emph{Seed mutation in rank~$2$,
    cf.~\cite[Example~2.5]{ca1}}] 
\label{example:seed-mut2}
For $n=2$, the mutation rules simplify considerably. 
Let $\xx=\{z,x\}$, $\pp=(p^+_z,p^-_z,p^+_x,p^-_x)$, and
\[
B=\begin{bmatrix}
0 & b_{zx}\\
b_{xz} & 0
\end{bmatrix}
=\begin{bmatrix}
0 & -b \\
c & 0
\end{bmatrix}
\qquad (b,c>0). 
\]
Performing mutation in direction~$z$, we obtain: 
\begin{itemize}
\item
the cluster $\overline \xx=\{\overline z,x\}$, 
with $\overline z$ determined by the exchange relation 
\begin{equation}
\label{eq:exchange-rel-rk2}
z\,\overline z
=p_z^+ \, x^c + p_z^- ;
\end{equation}
\item
the coefficient tuple
$\overline \pp=(\overline p^+_{\overline z},\overline p^-_{\overline
  z},\overline p^+_x,\overline p^-_x)$ 
satisfying $\overline p ^\pm_{\overline z}=p ^\mp_z$ and 
\begin{equation}
\label{eq:p-mutation-rk2}
\frac{\overline p ^+_x}{\overline p ^-_x} =
(p_z^-)^{-b}\,\dfrac{p ^+_x}{p ^-_x}\,;
\end{equation}
\item
the exchange matrix
\[
\overline B=\begin{bmatrix}
0 & \overline b_{\overline zx}\\
\overline b_{x\overline z} & 0
\end{bmatrix}
=\begin{bmatrix}
0 & b \\
-c & 0
\end{bmatrix}\,. 
\]
\end{itemize}
\end{example}

\begin{definition}[\emph{Exchange pattern}]
\label{def:exch-pattern}
Let $\Exch$ be a connected (unoriented, possibly infinite)
$n$-regular graph.
That is, each vertex $t$ in~$\Exch$ is connected to exactly $n$ other
vertices; the corresponding $n$-element set of edges is denoted by
$\Star(t)$.
An \emph{attachment} of a seed $\Sigma=(\xx,\pp,B)$ at $t$ is a
bijective labeling of the $n$ cluster variables $x\in\xx$ 
(hence the associated coefficient pairs~$p^\pm_x$, 
and the rows/columns of~$B$) by the $n$ edges
in $\Star(t)$; cf.\ \cite[Section~2.2]{ca2}.
With such a seed attachment, we typically use the natural notation
\begin{align*}
\xx=\xx(t)&=(x_e(t))_{e\in\Star(t)}\,, \\
\pp=\pp(t)&=(p^\pm_e(t))_{e\in\Star(t)}\,,&
p^\pm_e(t)&=p^\pm_{x_e(t)}\,,\\
B=B(t)&=(b_{ef}(t))_{e,f\in\Star(t)}\,,&
b_{ef}(t)&=b_{x_e(t),x_f(t)}\,. 
\end{align*}
An \emph{exchange pattern} on~$\Exch$ is, informally
speaking, a collection of seeds attached at the vertices  of~$\Exch$
and related to each other by the corresponding mutations.
Let us now be precise.
An exchange pattern $\ExPat=(\Sigma_t)$ is a collection of
seeds $\Sigma_t=(\xx(t),\pp(t),B(t))$ labeled by the vertices~$t$
in~$\Exch$ and satisfying the following condition.
Consider an edge $e$ in $\Exch$ connecting two vertices
$t$ and~$\overline t$.
Let $\xx(t)=(x_f(t))_{f\in\Star(t)}$ be the cluster at~$t$,
labeled using the seed attachment.
Then the definition requires that the seed
$\Sigma_{\overline t}$ is related to $\Sigma_t$ by a
mutation in the direction of~$x_e(t)$.
In particular, we have 
\[
\xx(\overline t)=\xx(t)-\{x_e(t)\}\cup\{x_e(\overline t)\}, 
\]
while the exchange relation \eqref{eq:exchange-rel-xx}
associated with the edge $e$ takes the form
\begin{equation}
\label{eq:exchange-rel-E}
x_e(t)\,x_e(\overline t)
=p_e^+(t) \, \prod_{\substack{f\in\Star(t) \\ b_{fe}(t)>0}} x_f(t)^{b_{fe}(t)}
+p_e^-(t) \, \prod_{\substack{f\in\Star(t) \\ b_{fe}(t)<0}} x_f(t)^{-b_{fe}(t)}
\,.
\end{equation}
Since the adjacent clusters
$\xx(t)=(x_f(t))_{f\in\Star(t)}$ and
$\xx(\overline t)=(x_{\overline f}(\overline t))_{\overline f\in\Star(\overline t)}$
coincide (as unlabeled sets) except for the replacement of $x_e(t)$ by
$x_e(\overline t)$, we have a bijection
\begin{align*}
\Star(t)\setminus\{e\}&\to\Star(\overline t)\setminus\{e\}\\
f&\mapsto \overline f
\end{align*}
(defined by $x_f(t)=x_{\overline f}(\overline t)$)
between the edges incident to $t$ and~$\overline t$, respectively.
These bijections form a \emph{discrete connection}
(see, e.g., \cite{guillemin-zara}) on~$\Exch$ that keeps
track of the relabelings of each cluster variable.
The corresponding
equations \eqref{eq:p-mutation1}--\eqref{eq:p-mutation2} become
\begin{align}
\label{eq:p-mutation1-E}
\overline p ^\pm_e(\overline t)&=p ^\mp_e(t)\,;\\ 
\label{eq:p-mutation2-E}
\frac{p ^+_{\overline f}(\overline t)}{p ^-_{\overline f}(\overline t)} &=
\begin{cases}
(p_e^+(t))^{b_{ef}(t)}\,\dfrac{p ^+_f(t)}{p ^-_f(t)}
& \text{if $b_{ef}(t)\geq 0$}\\[.2in]
(p_e^-(t))^{b_{ef}(t)}\,\dfrac{p ^+_f(t)}{p ^-_f(t)}
& \text{if $b_{ef}(t)\leq 0$}
\end{cases}
\text{\qquad (for $f\neq e$).}
\end{align}
The matrix mutation rule can be similarly rewritten in this notation. 
\end{definition}

We note that an exchange pattern on a regular graph $\Exch$ can be
canonically lifted to an exchange pattern on its universal cover, an
$n$-regular tree~$\mathbb{T}_n$,
recovering the original definition in~\cite{ca1}.

\begin{figure}[htbp]
\begin{center}
\setlength{\unitlength}{1.8pt}
\begin{picture}(200,8)(0,-2)
\thicklines

\multiput(10,0)(60,0){4}{\circle*{3}}

\put(0,0){\line(1,0){200}}

\put(100,4){\makebox(0,0){$e$}}
\put(40,4){\makebox(0,0){$f$}}
\put(160,4.5){\makebox(0,0){$\overline f$}}

\put(70,6){\makebox(0,0){$t$}}
\put(130,6.5){\makebox(0,0){$\overline t$}}

\end{picture}
\end{center}
\caption{Exchange graph of rank~2}
\label{fig:expat2}
\end{figure}
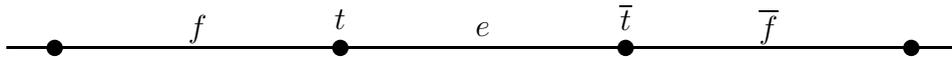

\begin{example}[\emph{Exchange pattern of rank~2}]
\label{example:expat-rk2}
We continue with Example~\ref{example:seed-mut2}.  
Figure~\ref{fig:expat2} shows a fragment of a $2$-regular
graph~$\Exch$.
We attach the seeds $(\xx,\pp,B)$ and
$(\overline\xx,\overline\pp,\overline B)$ at the adjacent
vertices $t$ and~$\overline t$, so that
\[
\begin{array}{l@{\quad}l}
\xx=(z,x)=\xx(t)=(x_e(t),x_f(t)), &
\overline \xx=(\overline z,x)
=\xx(\overline t)=(x_e(\overline t),x_{\overline f}(\overline t)), \\[.1in]
\pp=(p^\pm_z,p^\pm_x)=\pp(t)=(p^\pm_e(t),p^\pm_f(t)), &
\overline \pp=(\overline p^\pm_{\overline z},\overline p^\pm_x)
   =\pp(\overline t)=(p^\pm_e(\overline t),p^\pm_{\overline f}(\overline t)), \\[.1in] 
B\!=\!B(t)\!=\!\begin{bmatrix}
0 & b_{ef}(t)\\
b_{fe}(t) & 0
\end{bmatrix}
\!=\!\begin{bmatrix}
0 & -b \\
c & 0
\end{bmatrix}, &
\overline B\!=\!B(\overline t) \!=\!\begin{bmatrix}
0 & b_{e\overline f}(\overline t)\\
b_{\overline fe}(\overline t) & 0
\end{bmatrix}
\!=\!\begin{bmatrix}
0 & b \\
-c & 0
\end{bmatrix}.
\end{array}
\]
The corresponding exchange relation \eqref{eq:exchange-rel-rk2} becomes
\begin{equation}
\label{eq:exch-rk2}
x_e(t)\,x_e(\overline t)
=p_e^+(t) \, x_f(t)^c + p_e^-(t), 
\end{equation}
while the relation \eqref{eq:p-mutation-rk2}
becomes 
\begin{equation}
\label{eq:p-exch-rk2}
\frac{p ^+_{\overline f}(\overline t)}{p ^-_{\overline f}(\overline t)} =
(p^-_e(t))^{-b}\,\dfrac{p ^+_f(t)}{p ^-_f(t)}\,. 
\end{equation}
\end{example}

\medskip

\pagebreak[3]

\begin{definition}[\emph{Cluster algebra}]
\label{def:nn-ca}
To define a \emph{(non-normalized) cluster algebra}~$\Acal$,
one needs two pieces of data:
\begin{itemize}
\item
an exchange pattern
$\ExPat=(\Sigma_t)$ as in Definition~\ref{def:exch-pattern};
\item
a \emph{ground ring}~$\Rcal$,
a subring with unit in $\ZZ\PP$ that contains all
coefficient tuples~$\pp(t)$, for all seeds $\Sigma_t=(\xx(t), \pp(t),
B(t))$.
\end{itemize}
The cluster algebra $\Acal$ is then defined as
the $\Rcal$-subalgebra of the ambient
field $\Fcal$ generated by the union of all clusters~$\xx(t)$.
\end{definition}

We conclude this chapter by showing that although the rules of
non-normalized seed mutation do not determine the mutated seed
uniquely (see Remark~\ref{rem:nn-vs-n}), certain expressions
derived from each seed in a non-normalized exchange pattern form a discrete
dynamical system, i.e., satisfy a self-contained set of mutation-like
recurrences. 

Proposition~\ref{pr:yhat} below is an extension to the non-normalized case
of an important observation already made in \cite[Lemma~2.11]{fock-goncharov2} (for trivial coefficients),
in \cite[Lemma~1.3]{gsv1} (for tropical coefficients), and
in \cite[Proposition 3.9]{ca4} (for arbitrary normalized
coefficients). 
For a (non-normalized) exchange pattern as in
Definition~\ref{def:exch-pattern}, 
let us denote
\begin{equation}
\label{eq:yhat}
\hat y_e(t)=\frac{p^+_e(t)}{p^-_e(t)}\prod_{f\in\Star(t)}
x_f(t)^{b_{fe}(t)}  
\end{equation}
(cf.\ \cite[(3.7)]{ca4}). 
Note that $\hat y_e(t)$ is nothing but the ratio of the two terms on
the right-hand side of the exchange relation~\eqref{eq:exchange-rel-E}.
Surprisingly, the $n$-tuple  
\[
\hat\yy(t)=(\hat y_e(t))_{e\in\Star(t)}
\]
uniquely determines all adjacent $n$-tuples $\hat\yy(\overline t)$: 

\begin{proposition}
\label{pr:yhat}
Let $e$ be an edge in $\Exch$ connecting $t$ and~$\overline t$. 
Then
\begin{equation}
\label{eq:yhat-mutation}
\hat y_{\overline f}(\overline t) =
\begin{cases}
\hat y_e(t)^{-1}  & \text{if $f=e$};\\[.05in]
\hat y_f(t) \,\hat y_e(t)^{[b_{ef}(t)]_+}
(\hat y_e(t) + 1)^{- b_{ef}(t)} & \text{if $f \neq e$.}
\end{cases}
\end{equation}
\end{proposition}

In the language of \cite{ca4}, the equation~\eqref{eq:yhat-mutation}
means that the quantities $\hat y_e(t)$ form a (normalized)
\emph{$Y\!$-pattern} in~$\Fcal$. 

\begin{proof}
The case $f=e$ is immediate from \eqref{eq:B-mutation} 
and~\eqref{eq:p-mutation1-E}. 
For $f\neq e$, we obtain: 
\begin{align*}
\hat y_{\overline f}(\overline t)
&=\frac{p^+_{\overline f}(\overline t)}
      {p^-_{\overline f}(\overline t)} \, 
  \prod_{\overline g\in\Star(\overline t)} 
        x_{\overline g}(\overline t)^{b_{\overline g\overline
	    f}(\overline t)}
\\
\shortintertext{(by \eqref{eq:yhat})}
&=  \dfrac{p ^+_f(t)}{p ^-_f(t)}
   \left(\dfrac{p_e^+(t)}{p_e^-(t)}\right)^{[b_{ef}(t)]_+}
   \left(\dfrac{x_e (\overline t)}{p_e^-(t)}\right)^{-b_{ef}(t)}
   \prod_{\overline g\neq e} 
        x_{g}(t)^{b_{\overline g\overline
	    f}(\overline t)}
\\
\shortintertext{(by \eqref{eq:p-mutation2-new})}
&=\hat y_f(t)\left(\prod_{g\neq e} x_g(t)^{-b_{gf}(t)+b_{\overline g\overline
	    f}(\overline t)}\right) 
  \left(\dfrac{x_e(t)x_e (\overline t)}{p_e^-(t)}\right)^{-b_{ef}(t)} 
  \left(\dfrac{p_e^+(t)}{p_e^-(t)}\right) ^{[b_{ef}(t)]_+}
\\
\shortintertext{(by \eqref{eq:yhat})}
&=\hat y_f(t)
    \left(\dfrac{p_e^+(t)}{p_e^-(t)}\prod x_g(t)^{b_{ge}(t)}\right)
    ^{[b_{ef}(t)]_+} 
\left(\dfrac{x_e(t)x_e (\overline t)}{p_e^-(t)}\prod x_g(t)^{
  -[-b_{ge}(t)]_+}\right)^{-b_{ef}(t)}
\\
\shortintertext{(by \eqref{eq:B-mutation-alt})}
&=\hat y_f(t) \,\hat y_e(t)^{[b_{ef}(t)]_+}
(\hat y_e(t) + 1)^{- b_{ef}(t)}
\end{align*}
(by \eqref{eq:exchange-rel-E} and \eqref{eq:yhat}).
\end{proof}

\begin{remark}
  \label{rem:mutation-yhat}
  Equation~\eqref{eq:B-mutation-alt} can now be recognized as, in
  some sense, a \emph{tropical} version of Proposition~\ref{pr:yhat}.
  (This statement can be made precise using the construction analogous
  to the one introduced in Definition~\ref{def:tropical}.)
\end{remark}

\begin{example}
Continuing with the rank~2 case of Example~\ref{example:expat-rk2}, we get
\begin{align}
\label{eq:yhat-rk2-1}
\hat y_e(t)&=\frac{p^+_e(t)}{p^-_e(t)} \, x_f(t)^{b_{fe}(t)}  
           =\frac{p^+_e(t)}{p^-_e(t)} \, x_f(t)^c,\\[.05in]
\label{eq:yhat-rk2-2}
\hat y_f(t)&=\frac{p^+_f(t)}{p^-_f(t)} \, x_e(t)^{b_{ef}(t)}  
           =\frac{p^+_f(t)}{p^-_f(t)} \, x_e(t)^{-b}, \\[.05in]
\label{eq:yhat-rk2-3}
\hat y_{\overline f}({\overline t})
           &=\frac{p^+_{\overline f}({\overline t})}
                  {p^-_{\overline f}({\overline t})} \, 
                     x_e(\overline t)^{b_{e\overline f}(\overline t)}  
           =\frac{p^+_{\overline f}({\overline t})}
                  {p^-_{\overline f}({\overline t})} \, 
                     x_e(\overline t)^b. 
\end{align}
Furthermore, \eqref{eq:yhat-mutation} becomes
\[
\hat y_{\overline f}({\overline t})
=\hat y_f(t) 
(\hat y_e(t) + 1)^b, 
\]
which is straightforward to verify using 
\eqref{eq:yhat-rk2-1}--\eqref{eq:yhat-rk2-3}, 
\eqref{eq:exch-rk2}, and \eqref{eq:p-exch-rk2}:  
\begin{align*}
\hat y_{\overline f}({\overline t})
&=\frac{p^+_{\overline f}({\overline t})}
                  {p^-_{\overline f}({\overline t})} \, 
                     x_e(\overline t)^b\\
&=(p^-_e(t))^{-b}\,\dfrac{p ^+_f(t)}{p ^-_f(t)}\,
\left(\dfrac{p_e^+(t) \, x_f(t)^c + p_e^-(t)}{x_e(t)}
\right)^b\\
&=\hat y_f(t)(\hat y_e(t) + 1)^b. 
\end{align*}

\end{example}

\chapter{Rescaling and normalization}
\label{sec:rescaling}

In this chapter, we make a couple of observations related to rescaling
of cluster variables;
these observations will play an important role in the sequel.

The first algebraic observation (see Proposition~\ref{pr:cv-rescaling}
below) is that rescaling an exchange pattern
gives again an exchange pattern. That is, if we replace each cluster
variable by a new one that differs by a constant factor (these constant
factors can be chosen completely arbitrarily for different cluster
variables), and then rewrite each exchange relation in the obvious way in
terms of the new variables, then the coefficients in these new exchange
relations satisfy the monomial relations for an exchange pattern.

To formulate the above statement precisely,
we will need the following natural notion.
We say that a collection $(c_e(t))$
labeled by all pairs $(t,e)$ with $e\in\Star(t))$ is
\emph{compatible with the discrete connection} defined by an exchange
pattern~$\ExPat$ (cf.\ Definition~\ref{def:exch-pattern}) 
if, for any edge $e$ between $t$ and~$\overline t$,
we have $c_f(t)=c_g(\overline t)$ whenever $x_f(t)=x_g(\overline t)$.

\begin{proposition}
\label{pr:cv-rescaling}
Let $\ExPat=(\Sigma_t)$ be an exchange pattern on an
$n$-regular graph~$\Exch$, as in Definition~\ref{def:exch-pattern}.
Let $(c_e(t))$ be a collection of scalars
in~$\PP$ that is compatible with the discrete
connection associated with~$\ExPat$.
Then the following construction yields an exchange pattern
$\ExPat'=(\Sigma'_t)$ on~$\Exch$, with $\Sigma'_t=(\xx'(t),\pp'(t),B(t))$:
\begin{itemize}
\item
the attached cluster $\xx'(t)=(x'_e(t))_{e\in\Star(t)}$
is given by 
\[
x'_e(t)=\frac{x_e(t)}{c_e(t)};
\]
\item
the coefficient tuple $\pp'(t)=(p'^\pm_e(t))$ is defined by
\begin{equation}
\label{eq:rescale-p}
p'^\pm_e(t)
=\frac{p^\pm_e(t)}{c_e(t) c_e(\overline t)} \,
    \prod_{\pm b_{fe}(t)>0} c_f(t)^{\pm b_{fe}(t)},
\end{equation}
where $\overline t$ denotes the endpoint of $e$ different from~$t$;
\item
the exchange matrices $B(t)$ do not change.
\end{itemize}
\end{proposition}

\begin{proof}
It is straightforward to check that
substituting $x_e(t)=x'_e(t)c_e(t)$ into \eqref{eq:exchange-rel-E}
results into the requisite exchange relation in~$\ExPat'$.
Together with the compatibility condition, this ensures that adjacent
(attached) clusters are related by the corresponding mutations.
It remains to demonstrate that the rescaled coefficient tuples $\pp'(t)$
satisfy the requirements in Definition~\ref{def:seed-mutation}.
The condition $p'^\pm_e(\overline t)=p'^\mp_e(t)$
(cf.~\eqref{eq:p-mutation1-E}) is easily verified.
Finally, in order to check the equation~\eqref{eq:p-mutation2-E}
for the coefficients~$p'^\pm_e(t)$, we substitute the
expressions~\eqref{eq:rescale-p} into it, and factor out
\eqref{eq:p-mutation2-E} for the original pattern.
The resulting equation
\begin{align*}
&\prod_{\overline g\in\Star(\overline t)}
c_{\overline g}(\overline t)^{b_{\overline g\overline f}(\overline
  t)}\\
=&
(c_e(t)c_e(\overline t))^{-b_{ef}(t)}
\prod_{g\,:\,b_{ge}(t)\,b_{ef}(t)>0} c_g(t)^{b_{ge}(t)|b_{ef}(t)|}
\prod_{g\in\Star(t)} c_g(t)^{b_{gf}(t)}
\end{align*}
is easily seen to follow from the matrix mutation
rules~\eqref{eq:B-mutation}.
\end{proof}

\begin{remark}
For a more intuitive explanation of why the axioms
of an exchange pattern survive rescaling,
check this property against
the alternative version of the mutation rules
\eqref{eq:p-mutation1}--\eqref{eq:B-mutation} given in
\cite[(2.7)]{ca1} (cf.\ also the ``Caterpillar Lemma''
in~\cite{laurent}).
\end{remark}

We next turn to the issue of \emph{normalization}, that is,
using the rescaling of cluster variables
(as in Proposition~\ref{pr:cv-rescaling})
to obtain a ``normalized'' exchange pattern.
The latter concept requires endowing the coefficient group~$\PP$
with a semifield structure (cf.\ Remark~\ref{rem:nn-vs-n}).

\begin{definition}[\emph{Normalized exchange pattern}]
\label{def:n-exch-pattern}
Suppose that
$(\PP,\oplus,\cdot)$ is a (commutative) \emph{semifield},
i.e., $(\PP,\cdot)$ is an abelian multiplicative group,
$(\PP,\oplus)$ is a commutative semigroup, and the \emph{auxiliary
  addition}~$\oplus$ is distributive with respect to the
multiplication.
(See Definition~\ref{def:tropical} for an example.)
The multiplicative group of any such semifield $\PP$ is torsion-free
\cite[Section~5]{ca1}.
An exchange pattern $\ExPat=(\Sigma_t)$ as in
Definition~\ref{def:exch-pattern}
(or the corresponding cluster algebra)
is called \emph{normalized} if the coefficients $p^\pm_e(t)$ satisfy the
normalization condition
\begin{equation}
\label{eq:normalization}
p^+_e(t)\oplus p^-_e(t)=1\,.
\end{equation}
\end{definition}

Our next algebraic observation is that rescaling of cluster variables in
a non-normalized exchange pattern produces a normalized pattern
if the rescaling factors $c_e(t)$ themselves
satisfy the auxiliary-addition version of the same exchange
relations.

\begin{proposition}
\label{pr:rescale-normalize}
Continuing with the assumptions and constructions of
Pro\-position~\ref{pr:cv-rescaling}, let us furthermore suppose that
the coefficient group $\PP$ is endowed with an additive
operation~$\oplus$ making $(\PP,\oplus,\cdot)$ a semifield.
Then the rescaled pattern $\ExPat'=(\Sigma'_t)$ is normalized
if and only if the scalars $c_e(t)$ satisfy the relations
\begin{equation}
\label{eq:exchange-rel-aux}
c_e(t)\,c_e(\overline t)
=p_e^+(t) \, \prod_{\substack{f\in\Star(t) \\ b_{fe}(t)>0}} c_f(t)^{b_{fe}(t)}
\oplus p_e^-(t) \, \prod_{\substack{f\in\Star(t) \\ b_{fe}(t)<0}}
c_f(t)^{-b_{fe}(t)}
\,.
\end{equation}
\end{proposition}

\begin{proof}
Equation \eqref{eq:exchange-rel-aux} is simply a rewriting of the
normalization condition $p'^+_e(t)\oplus p'^-_e(t)=1$
for the rescaled coefficients $p'^\pm_e(t)$ given by~\eqref{eq:rescale-p}.
\end{proof}

\chapter{Cluster algebras of geometric type and their positive
  realizations}
\label{sec:geom-type}

The most important example of normalized exchange patterns (resp.,
normalized cluster algebras) are the patterns (resp., cluster
algebras) of \emph{geometric type}.

\begin{definition}[{\emph{Tropical semifield, cluster algebra of
	geometric type} \cite[Example~5.6, Definition~5.7]{ca1}}]
\label{def:tropical}
Let $I$ be a finite indexing set, and let
\begin{equation}
\label{eq:PP=Trop}
\PP=\Trop (q_i: i \in I)
\end{equation}
be the multiplicative group of Laurent
monomials in the formal variables $\{q_i : i \in I\}$,
which we call the \emph{coefficient variables}.
Define the auxiliary addition $\oplus$ by
\begin{equation}
\label{eq:tropical addition}
\prod_i q_i^{a_i} \oplus \prod_i q_i^{b_i}  =
\prod_i q_i^{\min (a_i, b_i)}  .
\end{equation}
The semifield $(\PP,\oplus,\cdot)$ is called a \emph{tropical
  semifield}
(cf.\ \cite[Example~2.1.2]{bfz})
A cluster algebra (or the corresponding exchange pattern) is said to
be of \emph{geometric type} if it is
defined by a normalized exchange pattern with coefficients in some
tropical semifield $\PP=\Trop(q_i: i \in I)$,
over the ground ring~$\Rcal=\ZZ[q_i^\pm: i \in I]$
or $\Rcal=\ZZ[q_i: i \in I]$.

This definition differs slightly from the one used in
\cite{ca3,ca1,ca2,cdm}, where only the former choice of~$\Rcal$ was
allowed. 
See \textit{loc.\ cit.}\ for numerous examples. 
\end{definition}

\begin{definition}[\emph{Extended exchange matrix}]
\label{def:B-tilde}
For an exchange pattern of geometric type, the coefficients
$p^\pm_e(t)$ are monomials in the variables~$q_i$.
It is convenient and customary to encode these
coefficients, along
with the exchange matrix~$B(t)$, in a rectangular
\emph{extended exchange matrix} $\tilde B(t)=(b_{ef}(t))$ defined as
follows.
The~columns of $\tilde B(t)$ are, as before, labeled by $\Star(t)$.
The top $n$ rows of $\tilde B(t)$ are also labeled  by $\Star(t)$
while the subsequent rows are labeled by the elements of~$I$.
The top $n\times n$ submatrix of~$\tilde B(t)$ is $B(t)$ (so our
notation for the matrix elements is consistent); the entries of
the bottom $|I|\times n$ submatrix are uniquely determined by the formula
\[
\frac{p_e^+(t)}{p_e^-(t)}
=\prod_{i\in I} q_i^{b_{ie}(t)}.
\]
As observed in~\cite{ca1}, in
the case of geometric type the mutation rules
\eqref{eq:p-mutation1-E}, \eqref{eq:p-mutation2-E}, and
\eqref{eq:normalization} can be restated as saying that the matrices $\tilde B(t)$
undergo a matrix mutation given by the same formulas
\eqref{eq:B-mutation} as before---now with a different set of row
labels.

A pair $(\xx(t)),\tilde B(t))$ consisting of a cluster and the
corresponding extended exchange matrix will be referred to as a \emph{seed}
(of geometric type).
\end{definition}

\begin{example}[\emph{The ring $\CC[\operatorname{SL}_2]$}]
\label{example:mat22}
This was the first example given on the first page of the first paper
about cluster algebras~\cite{ca1}.
The coordinate ring
\[
\Acal=\CC[\operatorname{SL}_2]
=\CC[z_{11},z_{12},z_{21},z_{22}]/\langle z_{11}z_{22}-z_{12}z_{21}-1\rangle
\]
carries a structure of a cluster algebra of geometric type, of rank
$n=1$,
with the coefficient semifield $\Trop(z_{12},z_{21})$,
the cluster variables $z_{11}$ and~$z_{22}$,
and the sole exchange relation
\[
z_{11}z_{22}=z_{12}z_{21}+1.
\]
The clusters are $\xx(t_1)\!=\!\{z_{11}\}$ and $\xx(t_2)\!=\!\{z_{22}\}$. 
The extended exchange matrices~are
\[
\tilde B(t_1)=\begin{bmatrix}
0\\ 1\\ 1
\end{bmatrix},\quad
\tilde B(t_2)=\begin{bmatrix}
0\\ -1\\ -1
\end{bmatrix}.
\]
Cf.\ Example~\ref{example:sl3}.
\end{example}

The following concept is rooted in the original motivations of cluster
algebras, designed in part to study totally positive parts of
algebraic varieties of Lie-theoretic origin, in the sense of
G.~Lusztig~\cite{lusztig-red}
(see also \cite{dbc,cdm,lusztig-intro} and references therein).
In the context of cluster structures arising in Teichm\"uller theory,
similar notions were first considered in
\cite{fock-goncharov1, gsv2}.

\begin{definition}[\emph{Positive realizations}]
\label{def:positive-realization}
The positive realization of a cluster algebra $\Acal$ of geometric type is,
informally speaking, a faithful representation of $\Acal$ in the space of
positive real functions on a topological space $\Teich$ of appropriate
real dimension. An~accurate definition follows.

Let $\Acal$ be a cluster algebra of geometric type
over the semifield $\Trop(q_i: i \in I)$.
As before, let $n$ denote the rank of~$\Acal$, let $\Exch$ be the
underlying $n$-regular graph, and let
$\xx(t)=(x_e(t))_{e\in\Star(t)}$, for $t\in\Exch$, be the
clusters.
A \emph{positive realization} of~$\Acal$ is
a topological space $\Teich$ together with a collection of functions
$x_e(t):\Teich\to\RR_{>0}$ and $q_i:\Teich\to\RR_{>0}$
representing the cluster variables
and the coefficient variables, respectively, so that
\begin{itemize}
\item
these functions satisfy all appropriate exchange relations, and
\item
for each $t\in\Exch$, the map
\begin{equation}
\label{eq:pos-cv-homeo}
\prod_{e\in\Star(t)} x_e(t)\times \prod_{i\in I} q_i
:\Teich\to\RR_{>0}^{n+|I|}
\end{equation}
is a homeomorphism.
\end{itemize}
\end{definition}

The following simple observations will be useful in the sequel.

\pagebreak[3]

\begin{proposition}
\label{prop:positive-realization}
Every cluster algebra of geometric type has a positive realization,
unique up to canonical homeomorphism. 

Conversely, let $(\tilde B(t))$ be 
a collection of extended exchange matrices $\tilde B(t)$ labeled by
the vertices of an $n$-regular graph~$\Exch$ and related to each other
by the corresponding matrix mutations (cf.\
Definition~\ref{def:B-tilde}).
Suppose furthermore that there exists a topological space~$\Teich$  
and positive real functions $(x_e(t))$ and $(q_i)$ 
on $\Teich$ which satisfy the
conditions in Definition~\ref{def:positive-realization}.
Then the matrices $(\tilde B(t))$ define a cluster algebra of
geometric type, and the functions mentioned above provide its positive
realization. 
\end{proposition}

\begin{proof}
To construct a positive realization for a cluster algebra of geometric
type, start with an arbitrary
homeomorphism of the form~\eqref{eq:pos-cv-homeo} for one $t \in\Exch$;
then determine
the rest of the maps $x_e(t):\Teich\to\RR_{>0}$ using exchange
relations.
The key feature of exchange patterns that makes this construction work
is that the re-parametrization maps relating adjacent clusters are
birational and \emph{subtraction-free} (hence positivity preserving).

For the second part, use the fact that two rational functions in $m$
variables are equal if and only if they coincide pointwise as
functions on $\RR_{>0}^m$.
\end{proof}

The fact that positive realizations are unique up to canonical
homeomorphism allows us to speak of \emph{the} positive realization.

\chapter{Bordered surfaces, arc complexes, and tagged arcs}
\label{sec:cats1-recap}

This chapter offers a swift review of the main constructions
in~\cite{cats1}.
For a detailed exposition with lots of examples and pictures,
see \cite[Sections~2--5,~7]{cats1}.

\begin{figure}[htbp]
\begin{center}
\setlength{\unitlength}{3pt}
\begin{picture}(40,40)(0,0)
\thicklines

\put(20,20){\circle{40}}

\qbezier(13,21)(20,13)(27,21)
\qbezier(14.5,20)(20,25)(25.5,20)

\put(20,10){\circle{6}}
\multiput(17,10)(6,0){2}{\circle*{1.2}}
\put(20,13){\circle*{1.2}}

\multiput(14,29)(12,0){2}{\circle*{1.2}}

\thinlines
\put(20,10){\circle{5}}
\put(20,10){\circle{4}}
\put(20,10){\circle{3}}
\put(20,10){\circle{2}}
\put(20,10){\circle{1}}

\end{picture}

\end{center}
\caption{A bordered surface with marked points.
In this example, $\Surf$ is a  torus with a hole;
the boundary $\partial\Surf$ has a single
component with $3$~marked points on it;
and the set $\Mark$ consists of those $3$ points plus $2$~punctures in the
interior of~$\Surf$.}
\label{fig:bordered-surface}
\end{figure}
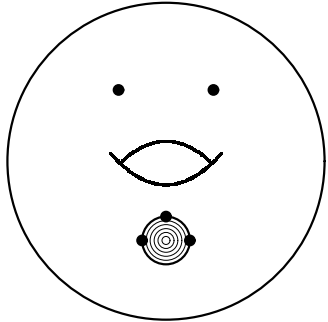

\begin{definition}[\emph{Bordered surface with marked points}]
Let $\Surf$ be a connected oriented $2$\hyp dimensional Riemann
surface with (possibly empty) boundary~$\partial\Surf$.
Fix a non-empty finite set $\Mark$
of \emph{marked points} in~$\Surf$, so that
there is at least one marked point on each connected component
of~$\partial\Surf$.
Marked points in the interior of~$\Surf$ are called \emph{punctures}.
We will want $\Surf$ to have at least one triangulation by a non-empty
set of arcs with endpoints at~$\Mark$.
Consequently, we do not allow $\Surf$ to be a sphere with one or two
punctures; nor an unpunctured or once-punctured monogon;
nor an unpunctured digon or triangle.
We also exclude the case of a sphere with three punctures. 
Such a pair $\SM$ is called a \emph{bordered surface with marked
  points}.
An example is shown in Figure~\ref{fig:bordered-surface}.
\end{definition}

\begin{figure}[htbp]
\begin{center}
\setlength{\unitlength}{3pt}
\begin{picture}(40,40)(0,0)
\thicklines

\darkred{
\qbezier(23,10)(30,13)(30,10)
\qbezier(23,10)(30,7)(30,10)
\qbezier(20,13)(10,18)(17,10)
\put(14,29){\circle{6}}
}

\put(20,20){\circle{40}}

\qbezier(13,21)(20,13)(27,21)
\qbezier(14.5,20)(20,25)(25.5,20)

\put(20,10){\circle{6}}
\multiput(17,10)(6,0){2}{\circle*{1.2}}
\put(20,13){\circle*{1.2}}

\multiput(14,29)(12,0){2}{\circle*{1.2}}

\thinlines
\put(20,10){\circle{5}}
\put(20,10){\circle{4}}
\put(20,10){\circle{3}}
\put(20,10){\circle{2}}
\put(20,10){\circle{1}}

\end{picture}
\qquad\qquad
\begin{picture}(40,40)(0,0)
\thicklines

\darkred{
\put(20,8){\circle{10}}
\qbezier(20,13)(37,18)(26,29)
\qbezier(10,29)(10,35)(26,29)
\qbezier(10,29)(10,23)(26,29)
}

\put(20,20){\circle{40}}

\qbezier(13,21)(20,13)(27,21)
\qbezier(14.5,20)(20,25)(25.5,20)

\put(20,10){\circle{6}}
\multiput(17,10)(6,0){2}{\circle*{1.2}}
\put(20,13){\circle*{1.2}}

\multiput(14,29)(12,0){2}{\circle*{1.2}}

\thinlines
\put(20,10){\circle{5}}
\put(20,10){\circle{4}}
\put(20,10){\circle{3}}
\put(20,10){\circle{2}}
\put(20,10){\circle{1}}

\end{picture}

\end{center}
\caption{The three curves on the left do not represent arcs;
the three curves on the right do.}
\label{fig:arcs}
\end{figure}
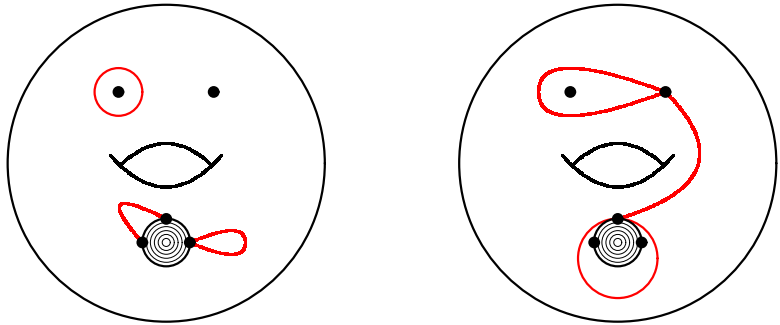

\begin{definition}[\emph{Ordinary arcs}]
An \emph{arc} $\gamma$ in $\SM$ is a curve
in~$\Surf$, considered up to isotopy, such that
\begin{itemize}
\item
the endpoints of $\gamma$ are marked points in~$\Mark$;
\item
$\gamma$ does not intersect itself, except that its endpoints may
  coincide;
\item
except for the endpoints, $\gamma$ is disjoint from~$\Mark$ and
from~$\partial\Surf$; and
\item
$\gamma$~does not cut out
an unpunctured monogon or an unpunctured digon.
\end{itemize}
An arc whose endpoints coincide is called a \emph{loop}.
We denote by $\APSM$ the set of of all arcs in $\SM$.
See Figure~\ref{fig:arcs}.
\end{definition}

\begin{definition}[\emph{Compatibility of ordinary arcs}]
Two arcs are \emph{compatible} if they
(more precisely, some of their isotopic deformations)
do not intersect in the interior of~$\Surf$.
For example, the three arcs shown in Figure~\ref{fig:arcs} on the
right are pairwise compatible.
\end{definition}

\pagebreak[3]

\begin{definition}[\emph{Ideal triangulations}]
\label{def:ideal-triang}
A maximal collection of distinct pairwise compatible arcs forms an
(ordinary) \emph{ideal triangulation}.
The arcs of a triangulation cut~$\Surf$ into \emph{ideal triangles};
note that we do allow \emph{self-folded} triangles, triangles where
two sides are identified.
Each ideal triangulation consists of
\begin{equation}
\label{eq:n=6g+3b+3p+c-6}
n=6g+3b+3p+c-6
\end{equation}
arcs, where
$g$ is the genus of~$\Surf$,
$b$ is the number of boundary components,
$p$ is the number of punctures, and
$c$ is the number of marked points on the boundary~$\partial\Surf$.
\end{definition}

Figure~\ref{fig:triang-4-punct-sphere} shows two triangulations of a
sphere with $4$~punctures.
Each triangulation has $4$ ideal triangles, $2$~of which are
self-folded.

\newsavebox{\digon}
\setlength{\unitlength}{4pt}
\savebox{\digon}(10,10)[bl]{
\qbezier(5,0)(0,5)(5,10)
\qbezier(5,0)(10,5)(5,10)
\put(5,0){\circle*{1}}
\put(5,5){\circle*{1}}
\put(5,10){\circle*{1}}
}

\newsavebox{\lowbar}
\savebox{\lowbar}(10,10)[bl]{
\put(5,0){\line(0,1){5}}
}

\newsavebox{\highbar}
\setlength{\unitlength}{4pt}
\savebox{\highbar}(10,10)[bl]{
\put(5,5){\line(0,1){5}}
}

\newsavebox{\lowloop}
\setlength{\unitlength}{4pt}
\savebox{\lowloop}(10,10)[bl]{
\qbezier(5,0)(2,7)(5,7)
\qbezier(5,0)(8,7)(5,7)
}

\newsavebox{\highloop}
\setlength{\unitlength}{4pt}
\savebox{\highloop}(10,10)[bl]{
\qbezier(5,10)(2,3)(5,3)
\qbezier(5,10)(8,3)(5,3)
}

\newsavebox{\lowtag}
\setlength{\unitlength}{4pt}
\savebox{\lowtag}(10,10)[bl]{
\put(5,3.5){\makebox(0,0){$\atag$}}
}

\newsavebox{\hightag}
\setlength{\unitlength}{4pt}
\savebox{\hightag}(10,10)[bl]{
\put(5,6.5){\makebox(0,0){$\atag$}}
}

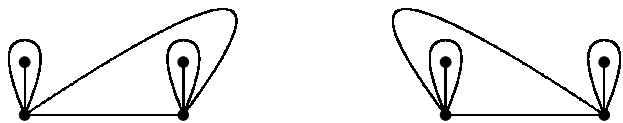
\begin{figure}[htbp]
\begin{center}

\setlength{\unitlength}{4pt}
\begin{picture}(30,12)(0,0)

\thinlines

\put(10,5){\makebox(0,0){\usebox{\lowbar}}}
\put(10,5){\makebox(0,0){\usebox{\lowloop}}}

\put(25,5){\makebox(0,0){\usebox{\lowbar}}}
\put(25,5){\makebox(0,0){\usebox{\lowloop}}}

\put(10,0){\circle*{1}}
\put(10,5){\circle*{1}}
\put(25,0){\circle*{1}}
\put(25,5){\circle*{1}}

\put(10,0){\line(1,0){15}}

\qbezier(10,0)(40,20)(25,0)
\end{picture}
\qquad\quad
\begin{picture}(30,12)(0,0)

\thinlines

\put(10,5){\makebox(0,0){\usebox{\lowbar}}}
\put(10,5){\makebox(0,0){\usebox{\lowloop}}}

\put(25,5){\makebox(0,0){\usebox{\lowbar}}}
\put(25,5){\makebox(0,0){\usebox{\lowloop}}}

\put(10,0){\circle*{1}}
\put(10,5){\circle*{1}}
\put(25,0){\circle*{1}}
\put(25,5){\circle*{1}}

\put(10,0){\line(1,0){15}}

\qbezier(10,0)(-5,20)(25,0)
\end{picture}
\end{center}
\caption{Two triangulations of a sphere with $4$ punctures.
In this example, $n=6$, $g=0$, $b=0$, $p=4$, $c=0$, in the notation
of~\eqref{eq:n=6g+3b+3p+c-6}.} 
\label{fig:triang-4-punct-sphere}
\end{figure}

The assumptions made above ensure that $\SM$
possesses a triangulation without self-folded triangles
(see \cite[Lemma~2.13]{cats1}).

\begin{figure}[htbp]
\begin{center}
\ \hspace{-.5in}
\setlength{\unitlength}{1.5pt}
\begin{picture}(60,60)(0,0)
\thicklines
  \put(0,20){\line(1,2){20}}
  \put(0,20){\line(1,-1){20}}
  \put(20,0){\line(0,1){60}}
  \put(20,0){\line(1,1){40}}
  \put(20,60){\line(2,-1){40}}

  \put(20,0){\circle*{3}}
  \put(20,60){\circle*{3}}
  \put(0,20){\circle*{3}}
  \put(60,40){\circle*{3}}

\put(15,25){\makebox(0,0){$\gamma$}}
\end{picture}
\begin{picture}(40,60)(0,0)
\put(20,30){\makebox(0,0){$\longrightarrow$}}
\end{picture}
\begin{picture}(60,60)(0,0)
\thicklines
  \put(0,20){\line(1,2){20}}
  \put(0,20){\line(1,-1){20}}
  \put(0,20){\line(3,1){60}}
  \put(20,0){\line(1,1){40}}
  \put(20,60){\line(2,-1){40}}

  \put(20,0){\circle*{3}}
  \put(20,60){\circle*{3}}
  \put(0,20){\circle*{3}}
  \put(60,40){\circle*{3}}

\put(26,35){\makebox(0,0){$\gamma'$}}

\end{picture}
\end{center}
\caption{A flip inside a quadrilateral}
\label{fig:flip}
\end{figure}
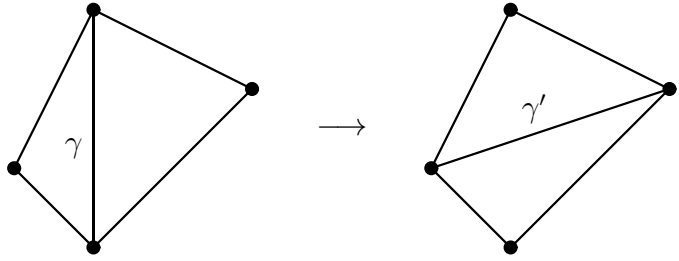

\begin{definition}[\emph{Ordinary flips}]
\label{def:flip}
Ideal triangulations are connected with each other by sequences of
\emph{flips}. Each flip replaces a single arc~$\gamma$ in a
triangulation~$T$ by a (unique) arc $\gamma'\neq \gamma$ that,
together with the remaining arcs in~$T$,
forms a new ideal triangulation.
See Figure~\ref{fig:flip}; also, the two triangulations in
Figure~\ref{fig:triang-4-punct-sphere} are related by a flip.

It is important to note that this operation cannot be applied to an
arc~$\gamma$ that lies inside a self-folded triangle in~$T$.
\end{definition}

\begin{figure}[htbp]
\begin{center}

\setlength{\unitlength}{4pt}
\begin{picture}(40,30)(0,0)

\thinlines

\put(10,5){\makebox(0,0){\usebox{\digon}}}
\put(10,5){\makebox(0,0){\usebox{\lowloop}}}

\put(40,5){\makebox(0,0){\usebox{\digon}}}
\put(40,5){\makebox(0,0){\usebox{\highloop}}}

\put(10,25){\makebox(0,0){\usebox{\digon}}}
\put(10,25){\makebox(0,0){\usebox{\lowbar}}}

\put(40,25){\makebox(0,0){\usebox{\digon}}}
\put(40,25){\makebox(0,0){\usebox{\highbar}}}

\put(15,5){\circle*{1.5}}
\put(35,5){\circle*{1.5}}
\put(15,25){\circle*{1.5}}
\put(35,25){\circle*{1.5}}

\thicklines

\put(15,5){\line(0,1){20}}
\put(35,5){\line(0,1){20}}
\put(15,25){\line(1,0){20}}

\put(25,0){\makebox(0,0){$\DPSM$}}

\end{picture}
\quad\qquad
\setlength{\unitlength}{4pt}
\begin{picture}(40,32)(0,0)

\thinlines


\put(40,15){\makebox(0,0){\usebox{\digon}}}
\put(40,15){\makebox(0,0){\usebox{\highbar}}}
\put(40,15){\makebox(0,0){\usebox{\highloop}}}

\put(10,15){\makebox(0,0){\usebox{\digon}}}
\put(10,15){\makebox(0,0){\usebox{\lowbar}}}
\put(10,15){\makebox(0,0){\usebox{\lowloop}}}

\put(21,28){\makebox(0,0){\usebox{\digon}}}
\put(21,28){\makebox(0,0){\usebox{\lowbar}}}
\put(21,28){\makebox(0,0){\usebox{\highbar}}}

\put(25,25){\circle*{1.5}}
\put(15,15){\circle*{1.5}}
\put(35,15){\circle*{1.5}}

\thicklines

\put(15,15){\line(1,1){10}}
\put(25,25){\line(1,-1){10}}

\put(25,0){\makebox(0,0){$\EPSM$}}

\end{picture}
\qquad\quad

\end{center}
\caption{The arc complex and its dual graph for a once-punctured digon}
\label{fig:arc-cplx+exch-graph-A1A1}
\end{figure}
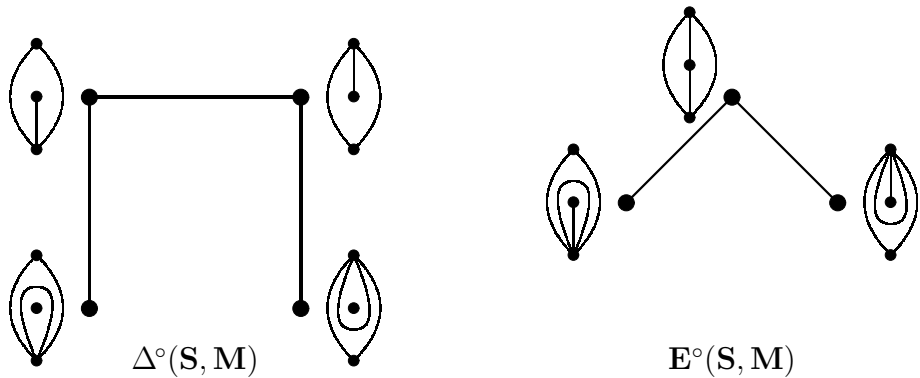

\begin{definition}[\emph{Arc complex and its dual graph}]
The \emph{arc complex} $\DPSM$ is the (possibly infinite) simplicial complex
on the ground set $\APSM$ defined as the clique
complex for the compatibility relation.
In other words, the vertices of $\DPSM$ are the arcs,
and the maximal simplices are the ideal triangulations.
The dual graph of $\DPSM$ is denoted by~$\EPSM$; its vertices are the
triangulations, and its edges correspond to the flips.
See Figure~\ref{fig:arc-cplx+exch-graph-A1A1}.
\end{definition}

In general, the arc complex $\DPSM$ has nonempty boundary since its
dual graph is not $n$-regular: not every arc can be flipped.
In~\cite{cats1}, we suggested a natural way to extend the arc complex
beyond its boundary, obtaining an $n$-regular graph which can be used to build
the desired exchange patterns.
This requires
the concept of a tagged arc.

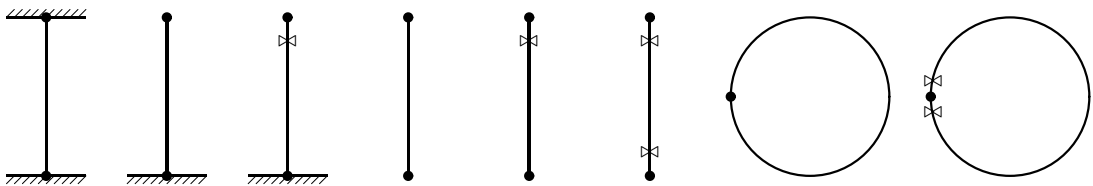
\begin{figure}[htbp]
\begin{center}
\setlength{\unitlength}{3pt}
\begin{picture}(10,20)(0,0)
\thinlines
\multiput(1,0)(1,0){10}{\line(-1,-1){1}}
\multiput(1,21)(1,0){10}{\line(-1,-1){1}}
\thicklines
\put(5,0){\line(0,1){20}}
\put(0,0){\line(1,0){10}}
\put(0,20){\line(1,0){10}}
\multiput(5,0)(0,20){2}{\circle*{1}}
\end{picture}
\quad
\begin{picture}(10,20)(0,0)
\thinlines
\multiput(1,0)(1,0){10}{\line(-1,-1){1}}
\thicklines
\put(5,0){\line(0,1){20}}
\put(0,0){\line(1,0){10}}
\multiput(5,0)(0,20){2}{\circle*{1}}
\end{picture}
\quad
\begin{picture}(10,20)(0,0)
\thinlines
\multiput(1,0)(1,0){10}{\line(-1,-1){1}}
\thicklines
\put(5,0){\line(0,1){20}}
\put(0,0){\line(1,0){10}}
\put(5,17){\makebox(0,0){$\notch$}}
\multiput(5,0)(0,20){2}{\circle*{1}}
\end{picture}
\quad
\begin{picture}(10,20)(0,0)
\thicklines
\put(5,0){\line(0,1){20}}
\multiput(5,0)(0,20){2}{\circle*{1}}
\end{picture}
\quad
\begin{picture}(10,20)(0,0)
\thicklines
\put(5,0){\line(0,1){20}}
\put(5,17){\makebox(0,0){$\notch$}}
\multiput(5,0)(0,20){2}{\circle*{1}}
\end{picture}
\quad
\begin{picture}(10,20)(0,0)
\thicklines
\put(5,0){\line(0,1){20}}
\put(5,17){\makebox(0,0){$\notch$}}
\put(5,3){\makebox(0,0){$\notch$}}
\multiput(5,0)(0,20){2}{\circle*{1}}
\end{picture}
\quad
\begin{picture}(20,20)(-5,0)
\thicklines
\put(5,10){\circle{20}}
\put(-5,10){\circle*{1}}
\end{picture}
\quad
\begin{picture}(20,20)(-5,0)
\thicklines
\put(5,10){\circle{20}}
\put(-4.7,12){\makebox(0,0){$\notch$}}
\put(-4.7,8){\makebox(0,0){$\notch$}}
\put(-5,10){\circle*{1}}
\end{picture}

\end{center}
\caption{Different types of tagged arcs}
\label{fig:tagged-arcs}
\end{figure}

\begin{definition}[\emph{Tagged arcs}]\label{def:tagged-arc}
A \emph{tagged arc} is obtained by taking an arc
that does not cut out a once-punctured monogon and marking (``tagging'')
each of its ends in one of the two ways,
\emph{plain} or \emph{notched}, so that the following conditions are
satisfied:
\begin{itemize}
\item
an endpoint lying on the boundary of $\Surf$ must be tagged plain, and
\item
both ends of a loop must be tagged in the same way.
\end{itemize}
See Figure~\ref{fig:tagged-arcs}.

The set of all tagged arcs in $\SM$ is denoted by~$\ATSM$.
\end{definition}

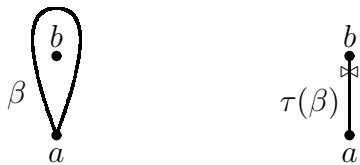
\begin{figure}[htbp]
\begin{center}
\setlength{\unitlength}{1.5pt}
\begin{picture}(40,33)(0,-3)
\thicklines


\qbezier(20,0)(8,32)(20,32)
\qbezier(20,0)(32,32)(20,32)

\multiput(20,0)(0,20){2}{\circle*{2}}
\put(10,10){\makebox(0,0){$\beta$}}
\put(20,-5){\makebox(0,0){$a$}}
\put(20,25){\makebox(0,0){$b$}}
\end{picture}
\qquad\qquad
\begin{picture}(40,33)(0,-3)
\thicklines

\put(20,0){\line(0,1){20}}
\put(20,16){\makebox(0,0){$\notch$}}
\put(10,8){\makebox(0,0){$\tau(\beta)$}}

\multiput(20,0)(0,20){2}{\circle*{2}}
\put(20,-5){\makebox(0,0){$a$}}
\put(20,25){\makebox(0,0){$b$}}
\end{picture}

\end{center}
\caption{\hbox{Representing an arc bounding a punctured monogon by a
    tagged arc}}
\label{fig:arc-as-tagged-arc}
\end{figure}

\begin{definition}[\emph{Representing ordinary arcs by tagged arcs}]
\label{def:arcs-as-tagged-arcs}
Ordinary arcs can be viewed as a special case of tagged arcs,
via the following dictionary.
Let us canonically represent any ordinary (untagged) arc~$\beta$ by a tagged
arc~$\tau(\beta)$ defined as follows.
If $\beta$ does not cut out a once-punctured monogon,
then~$\tau(\beta)$ is simply $\beta$ with both ends tagged plain.
Otherwise, $\beta$ is a loop based at some marked point~$a$ and
cutting out a punctured monogon with the sole puncture~$b$ inside it.
Let $\alpha$ be the unique arc connecting $a$ and $b$ and compatible
with~$\beta$. Then $\tau(\beta)$ is
obtained by tagging $\alpha$ plain at~$a$ and notched at~$b$.
See Figure~\ref{fig:arc-as-tagged-arc}.
\end{definition}

\begin{definition}[\emph{Compatibility of tagged arcs}]
This is an extension of the corresponding notion for ordinary arcs.
Tagged arcs $\alpha$ and $\beta$ are \emph{compatible}
if and only~if
\begin{itemize}
\item
their untagged versions $\alpha^\circ$ and $\beta^{\,\circ}$ are
compatible;
\item
if $\alpha$ and $\beta$ share an endpoint~$a$,
then the ends of $\alpha$ and $\beta$ connecting to~$a$
must be tagged in the same way---unless
$\alpha^\circ=\beta^{\,\circ}$, in which case at least one end of
$\alpha$ must be tagged in the same way as the corresponding end
of~$\beta$.
\end{itemize}
\end{definition}

It is easy to see that the map $\gamma\mapsto\tau(\gamma)$ described in
Definition~\ref{def:arcs-as-tagged-arcs} preserves compatibility.

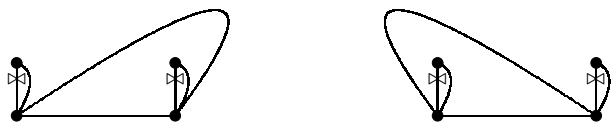
\begin{figure}[htbp]
\begin{center}

\setlength{\unitlength}{4pt}
\begin{picture}(30,10)(0,0)

\thinlines

\put(10,5){\makebox(0,0){\usebox{\lowbar}}}
\put(10,5){\makebox(0,0){\usebox{\lowtag}}}
\qbezier(10,5)(12.5,4)(10,0)

\put(25,5){\makebox(0,0){\usebox{\lowbar}}}
\put(25,5){\makebox(0,0){\usebox{\lowtag}}}
\qbezier(25,5)(27.5,4)(25,0)

\put(10,0){\circle*{1}}
\put(10,5){\circle*{1}}
\put(25,0){\circle*{1}}
\put(25,5){\circle*{1}}

\put(10,0){\line(1,0){15}}

\qbezier(10,0)(40,20)(25,0)
\end{picture}
\qquad\quad
\begin{picture}(30,10)(0,0)

\thinlines

\put(10,5){\makebox(0,0){\usebox{\lowbar}}}
\put(10,5){\makebox(0,0){\usebox{\lowtag}}}
\qbezier(10,5)(12.5,4)(10,0)

\put(25,5){\makebox(0,0){\usebox{\lowbar}}}
\put(25,5){\makebox(0,0){\usebox{\lowtag}}}
\qbezier(25,5)(27.5,4)(25,0)

\put(10,0){\circle*{1}}
\put(10,5){\circle*{1}}
\put(25,0){\circle*{1}}
\put(25,5){\circle*{1}}

\put(10,0){\line(1,0){15}}

\qbezier(10,0)(-5,20)(25,0)
\end{picture}
\end{center}
\caption{Two tagged triangulations of a $4$-punctured sphere.}
\label{fig:tagged-triang-4-punct-sphere}
\end{figure}

\begin{definition}[\emph{Tagged triangulations}]
\label{def:tagged-triang}
A~maximal (by inclusion) collection of pairwise compatible
tagged arcs is called a
\emph{tagged triangulation}.

Each ideal triangulation $T$ can be represented by a tagged
triangulation $\tau(T)$ via the dictionary~$\tau$ described above.
Figure~\ref{fig:tagged-triang-4-punct-sphere} shows two tagged
triangulations obtained by applying~$\tau$ to the
triangulations in Figure~\ref{fig:triang-4-punct-sphere}.
\end{definition}

All tagged
triangulations have the same cardinality~$n$
given by~\eqref{eq:n=6g+3b+3p+c-6} \cite[Theorem~7.9]{cats1}.

\begin{figure}[ht]
\begin{center}
\setlength{\unitlength}{4pt}
\begin{picture}(40,35)(5,-5)

\thinlines

\put(10,5){\makebox(0,0){\usebox{\digon}}}
\put(10,5){\makebox(0,0){\usebox{\lowbar}}}
\put(10,5){\makebox(0,0){\usebox{\lowtag}}}

\put(40,5){\makebox(0,0){\usebox{\digon}}}
\put(40,5){\makebox(0,0){\usebox{\highbar}}}
\put(40,5){\makebox(0,0){\usebox{\hightag}}}

\put(10,25){\makebox(0,0){\usebox{\digon}}}
\put(10,25){\makebox(0,0){\usebox{\lowbar}}}

\put(40,25){\makebox(0,0){\usebox{\digon}}}
\put(40,25){\makebox(0,0){\usebox{\highbar}}}

\put(15,5){\circle*{1.5}}
\put(35,5){\circle*{1.5}}
\put(15,25){\circle*{1.5}}
\put(35,25){\circle*{1.5}}

\thicklines

\put(15,5){\line(0,1){20}}
\put(35,5){\line(0,1){20}}
\put(15,25){\line(1,0){20}}
\put(15,5){\line(1,0){20}}

\put(25,-4){\makebox(0,0){$\DTSM$}}

\end{picture}
\quad\qquad
\begin{picture}(40,39)(5,-5)

\thinlines

\put(29,1){\makebox(0,0){\usebox{\digon}}}
\put(29,1){\makebox(0,0){\usebox{\lowbar}}}
\put(29,1){\makebox(0,0){\usebox{\lowtag}}}
\put(29,1){\makebox(0,0){\usebox{\highbar}}}
\put(29,1){\makebox(0,0){\usebox{\hightag}}}

\put(40,15){\makebox(0,0){\usebox{\digon}}}
\qbezier(40,15)(37.5,16)(40,20)
\put(40,15){\makebox(0,0){\usebox{\highbar}}}
\put(40,15){\makebox(0,0){\usebox{\hightag}}}

\put(10,15){\makebox(0,0){\usebox{\digon}}}
\qbezier(10,15)(12.5,14)(10,10)
\put(10,15){\makebox(0,0){\usebox{\lowbar}}}
\put(10,15){\makebox(0,0){\usebox{\lowtag}}}

\put(21,29){\makebox(0,0){\usebox{\digon}}}
\put(21,29){\makebox(0,0){\usebox{\lowbar}}}
\put(21,29){\makebox(0,0){\usebox{\highbar}}}

\put(25,5){\circle*{1.5}}
\put(25,25){\circle*{1.5}}
\put(15,15){\circle*{1.5}}
\put(35,15){\circle*{1.5}}

\thicklines

\put(25,5){\line(1,1){10}}
\put(25,5){\line(-1,1){10}}
\put(15,15){\line(1,1){10}}
\put(25,25){\line(1,-1){10}}

\put(15,0){\makebox(0,0){$\ETSM$}}

\end{picture}
\end{center}
\caption{Tagged arc complex and its dual graph for a once-punctured
  digon}
\label{fig:tagged-cplx-once-punct-digon}
\end{figure}
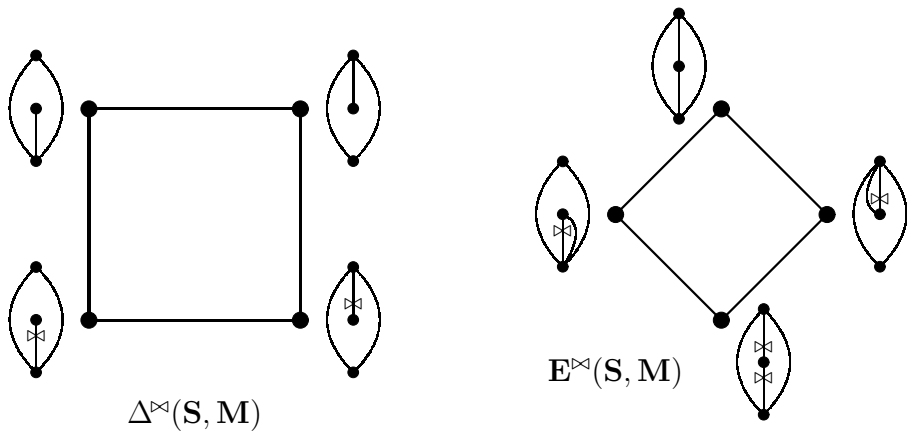

\begin{definition}[\emph{Tagged arc complex}]
The \emph{tagged arc complex} $\DTSM$ is the simplicial complex whose vertices
are tagged arcs and whose simplices are collections of pairwise compatible
tagged arcs.
See Figure~\ref{fig:tagged-cplx-once-punct-digon} on the left.

The ordinary arc complex $\DPSM$ can be viewed a subcomplex of~$\DTSM$
(via the map~$\tau$); cf.\ Figures~\ref{fig:arc-cplx+exch-graph-A1A1}
and~\ref{fig:tagged-cplx-once-punct-digon}.
\end{definition}

The maximal simplices of $\DTSM$ are the
tagged triangulations, so $\DTSM$ is \emph{pure} of dimension $n-1$.
Furthermore, $\DTSM$ is a
\emph{pseudomanifold}, i.e., each simplex of codimension~1 is
contained in precisely two maximal simplices.
To rephrase, for every tagged arc in an arbitrary tagged
triangulation, we can apply a \emph{tagged flip}
(replace it by a different tagged arc) in a unique way to
produce another tagged triangulation.

\begin{definition}[\emph{Dual graph of the tagged arc complex}]
The dual graph $\ETSM$  of the pseudomanifold $\DTSM$ has
tagged triangulations as its vertices. Two such vertices are connected
by an edge if these tagged triangulations are related by a tagged flip.
Thus, $\ETSM$ is a (possibly infinite) $n$-regular graph.
An example is shown in Figure~\ref{fig:tagged-cplx-once-punct-digon}
on the right.
\end{definition}

\begin{remark}[\emph{Nomenclature of tagged triangulations and tagged
    flips}]
\label{rem:two-types-of-flips}
Compatibility
of tagged arcs is invariant with respect to a simultaneous change of
all tags at a given puncture.
Let us take any tagged triangulation~$T'$ and perform such a change
at every puncture where all ends of~$T'$ are notched.
It is not hard to verify that the resulting tagged triangulation~$T''$
represents some ideal triangulation~$T$
(possibly containing self-folded triangles): $T''=\tau(T)$.
In other words, each tagged triangulation can be obtained from an
ordinary one by applying $\tau$ and then placing notches on all arcs
around some punctures.

This observation can be used to give a concrete description of all
possible tagged flips.
Let $T'$ be a tagged triangulation, and let $T''$ and~$T$ be as
above. Then each of the
$n$ tagged flips out of~$T'$ is of one of the two kinds:
\begin{enumerate}
\item[(\ref*{sec:cats1-recap}.D)]\label{item:digon-flip}
a flip performed inside a once-punctured digon,
as represented by one of the $4$~edges of the graph $\ETSM$ in
Figure~\ref{fig:tagged-cplx-once-punct-digon}.
The tagging at the vertices of the digon does not change; or
\item[(\ref*{sec:cats1-recap}.Q)]\label{item:quad-flip}
an ordinary flip inside a quadrilateral in~$T$ (cf.\
Definition~\ref{def:flip}).
As before, the sides of the quadrilateral do not have to be distinct.
Moreover, those sides (stripped of their tagging) should be arcs
of~$T$ but not necessarily of~$T'$:
specifically,
such a side can be a loop in~$T$ enclosing a once-punctured monogon.
The tagging at each vertex of the quadrilateral remains the same.
\end{enumerate}
See \cite[Section~9.2]{cats1} for more discussion and proofs.
For example, the two tagged triangulations in
Figure~\ref{fig:tagged-triang-4-punct-sphere} are related by a tagged
flip of type~(\ref{sec:cats1-recap}.Q).
\end{remark}

By \cite[Proposition~7.10]{cats1},
the dual graph $\ETSM$  (and therefore the complex $\DTSM$)
is connected---i.e., any two
tagged triangulations can be connected by a sequence of flips---unless
$\SM$ is a surface with no boundary and a single puncture, in which case
$\ETSM$ (resp., $\DTSM$)
consists of two isomorphic connected components, one in which
the ends of all tagged arcs are plain, and another in which they are
all notched.

\begin{definition}[\emph{Exchange graph of tagged triangulations}]
\label{def:ESM}
We denote by
$\ESM$ a connected component of $\ETSM$.
More precisely, we set $\ESM=\ETSM$ unless $\SM$ has no boundary and a single
puncture; in the latter case, $\ESM$ is the connected component of
$\ETSM$ in which all arcs are plain.
\end{definition}

As shown in \cite[Theorem~7.11]{cats1}, there is a natural class of
normalized exchange patterns (equivalently, cluster
algebras) whose underlying graph is $\mathbf{E}=\ESM$.
(Strictly speaking, this result was obtained under the assumption that
$\SM$ is not a closed surface with two punctures. We are not going to
make this assumption herein.)
In such a pattern, the cluster variables are labeled by
the tagged arcs while clusters correspond to tagged triangulations.
We review this construction.

\pagebreak[3]

\begin{figure}[ht]
\begin{center}
\setlength{\unitlength}{4pt}
\begin{picture}(40,41)(25,0)
\thinlines
  \put(10,0){\line(-1,2){10}}
  \put(40,20){\line(-1,2){10}}
  \put(10,0){\line(1,0){20}}
  \put(10,40){\line(1,0){20}}
  \put(0,20){\line(1,2){10}}
  \put(30,0){\line(1,2){10}}
  \put(10,0){\line(1,2){10}}
  \put(0,20){\line(1,0){40}}
  \put(10,40){\line(1,-2){10}}
  \put(20,20){\line(1,-2){10}}
  \put(10,40){\line(3,-2){30}}

  \put(0,20){\circle*{1}}
  \put(10,0){\circle*{1}}
  \put(10,40){\circle*{1}}
  \put(30,0){\circle*{1}}
  \put(30,40){\circle*{1}}
  \put(40,20){\circle*{1}}
  \put(20,20){\circle*{1}}

\thicklines

\put(23,10){\makebox(0,0){$2$}}
\put(13,10){\makebox(0,0){$3$}}
\put(10,18){\makebox(0,0){$4$}}
\put(30,18){\makebox(0,0){$1$}}
\put(27,31){\makebox(0,0){$6$}}
\put(13,31){\makebox(0,0){$5$}}

\put(45,19){
$B(T)=
\left[\,\,
\begin{matrix}
0 & -1& 0 & 0 & 1 &-1\\[.05in]
1 & 0 &-1 & 0 & 0 & 0\\[.05in]
0 & 1 & 0 &-1& 0 & 0\\[.05in]
0 & 0 & 1 & 0&-1 & 0\\[.05in]
-1& 0 & 0 & 1 & 0 & 1\\[.05in]
1 & 0 & 0 & 0 & -1& 0 
\end{matrix}
\,\,\right]
$
}
\end{picture}
\qquad\qquad

\end{center}
\caption{The signed adjacency matrix for a triangulation of a 
once-punctured hexagon}
\label{fig:B(T)-hexagon}
\end{figure}
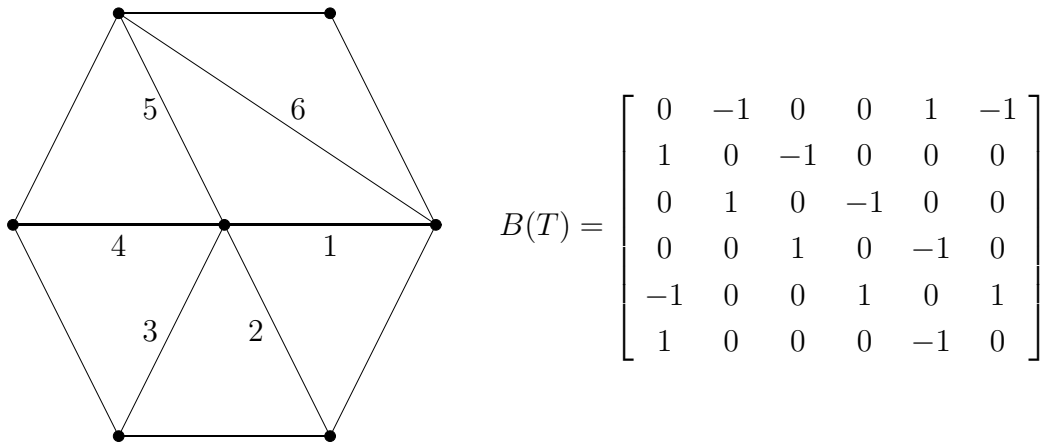

\begin{definition}[\emph{Signed adjacency matrix}]
\label{def:signed-adjacency-matrix}
The key ingredient in building an exchange pattern on $\ESM$ is a rule that
associates with each tagged triangulation $T$ a skew-symmetric
exchange matrix~$B(T)$ called a \emph{signed adjacency matrix}
of~$T$. The rows and columns of $B(T)$ are labeled by the arcs in~$T$.
The direct definition of $B(T)$ is fairly technical, and we refer the
reader to \cite[Definitions 4.1, 9.18]{cats1} for those
technicalities.
For the immediate purposes of this review, we make the following shortcut.
Let us start by defining $B(T)$ for an ordinary ideal triangulation~$T$
without self-folded triangles; this was first done
in~\cite{fock-goncharov1, fock-goncharov2, gsv2}.
Under that assumption, one sets 
\begin{equation}
\label{eq:B(T)=sum}
B(T)=\sum_\Delta B^\Delta\,,
\end{equation}
the sum over all ideal triangles $\Delta$ in~$T$ of the $n \times n$
matrices $B^\Delta=(b^\Delta_{ij})$ given by
\begin{equation}
\label{eq:B-proper}
b^\Delta_{ij}=
\begin{cases}
1 & \text{if $\Delta$ has sides $i$ and $j$, with $j$
  following~$i$ in the clockwise order;}\\
-1 & \text{if the same holds, with the counterclockwise order;}\\
0 & \text{otherwise.}
\end{cases}
\end{equation}
See Figure~\ref{fig:B(T)-hexagon} for an example.
\end{definition}

One can then extend the definition of $B(T)$ to arbitrary tagged
triangulations by requiring that 
\begin{itemize}
\item
whenever $T$ and $\overline T$ are
related by a flip of a tagged arc~$k$, the associated signed adjacency
matrices are related by the corresponding mutation:
\begin{equation}
\label{eq:B-mut-flip}
B(\overline T)=\mu_k(B(T)); 
\end{equation}
\item
if $T$ is a triangulation without self-folded triangles, and $T'$ is
obtained from $T$ by putting a notch at each end of each arc in~$T$
adjacent to a puncture, then $B(T')=B(T)$. 
\end{itemize}
It follows from \cite[Proposition~7.10, Definition~9.6,
  Lemma~9.7]{cats1} that this definition is
consistent, that is, there exists a unique 
collection of matrices $B(T)$ satisfying~\eqref{eq:B-mut-flip}. 

Each matrix $B(T)$ is skew-symmetric,
with entries equal to $0$, $\pm1$, or~$\pm2$.

It turns out that replacing ideal triangulations by (more general)
tagged triangulations does not extend the class of the associated
exchange matrices~$B(T)$ \cite[Proposition~12.3]{cats1}:
each matrix $B(T)$ corresponding to a tagged triangulation is
identical, up to simultaneous permutations of rows and columns, to a
matrix corresponding to an ordinary ideal triangulation.

\begin{remark}
\label{rem:trop-coeffs-from-boundary}
There are various ways to extend the matrices $B(T)$ to rectangular
matrices~$\tilde B(T)$ (cf.\ Definition~\ref{def:B-tilde}),
creating coefficient systems of geometric type.
A~very general construction of this kind will be discussed in
Chapter~\ref{sec:shear}.
Here we briefly discuss an easy special case that was already pointed
out in~\cite{fock-goncharov1, gsv2}.

Let $\BSM$ denote the set of boundary segments between adjacent marked
points on~$\partial\Surf$.
Consider the tropical coefficient semifield
\[
\PP=\Trop (q_\beta:\beta\in\BSM)
\]
generated by the variables $q_\beta$
labeled by such boundary segments.
For an ideal triangulation $T$ without self-folded triangles,
define the $(|\BSM|+n)\times n$ matrices $\tilde B(T)$ by the same
equations \eqref{eq:B(T)=sum}--\eqref{eq:B-proper} as before, but with
the understanding that the row index $i$ can now be either an arc or a
boundary segment.
The matrices $\tilde B(T)$ still satisfy the mutation
rule  $\tilde B(\overline T)=\mu_k(\tilde B(T))$.
This can be deduced
from \eqref{eq:B-mut-flip} by gluing a triangle on the other side of
each boundary segment (thus making it into a legitimate arc), then
``freezing'' all these arcs (i.e., not allowing to flip them).
The same argument allows us to extend the definition of $\tilde B(t)$
to arbitrary tagged triangulations, as was done for $B(t)$'s,
resulting in a well defined tropical coefficient system.
\end{remark}

\chapter{Structural results}
\label{sec:structural}

In this chapter, we formulate those of our results whose statements do
not require any references to Teichm\"uller theory or hyperbolic
geometry---even though their proofs will rely on geometric arguments.
These results concern structural properties of
exchange patterns (or cluster algebras) whose exchange matrices can be
described as signed adjacency matrices of triangulations of a bordered
surface.
Most crucially, we show
that, for any choice of (normalized) coefficients,
there is an exchange pattern on the $n$-regular graph
$\mathbf{E}=\ESM$ (see Definition~\ref{def:ESM})
whose exchange matrices are the signed adjacency matrices~$B(T)$.
More precisely, we have the following theorem.

\begin{theorem}
\label{th:patterns-on-SM}
Let $T_\circ$ be a tagged triangulation consisting of $n$ tagged arcs
in~$\SM$.
Let $\Sigma_\circ\!=\!(\xx(T_\circ),\pp(T_\circ),B(T_\circ))$
be a (normalized) seed as in Definitions~\ref{def:seed}
and~\ref{def:n-exch-pattern}, that is: 
\begin{itemize}
\item
$\xx(T_\circ)$ is an $n$-tuple of formal variables labeled by
the arcs in~$T_\circ$;
\item
$\pp(T_\circ)$ is a $2n$-tuple of elements of a semifield~$\PP$ satisfying~\eqref{eq:normalization};
\item
$B(T_\circ)$ is the signed adjacency matrix
of~$T_\circ$.
\end{itemize}
Then there is a unique exchange pattern $(\Sigma_T)$
on $\ESM$ such that $\Sigma_{T_\circ}=\Sigma_\circ$.

More precisely, let $\mathcal{F}$ be the field of rational functions
in the variables $\xx(T_\circ)$ with coefficients in~$\ZZ\PP$.
Then there exist unique elements $x_\gamma(T)\in\mathcal{F}$ and
$p^\pm_\gamma(T)\in\PP$  labeled by the tagged triangulations $T\in\ESM$
and the tagged arcs $\gamma\in T$ such that
\begin{itemize}
\item
every triple $\Sigma_T=(\xx(T),\pp(T),B(T))$ is a seed, where
$\xx(T)=(x_\gamma(T))_{\gamma\in T}$,
$\pp(T)=(p^\pm_\gamma(T))_{\gamma\in T}$,
and $B(T)$ is the signed adjacency matrix of~$T$;
\item
each cluster variable $x_\gamma=x_\gamma(T)$ does not depend on~$T$;
\item
if a tagged triangulation $T'$ is obtained from $T$ by flipping
a tagged arc $\gamma\in T$, then $\Sigma_T$ is obtained from
$\Sigma_{T'}$ by the seed mutation
replacing $x_\gamma$ by~$x_{\gamma'}$.
\end{itemize}
Furthermore, all cluster variables $x_\gamma$
(hence all seeds~$\Sigma_T$) are distinct.
\end{theorem}

In the terminology of~\cite[Section~7]{ca1}, the last statement means
that $\ESM$ is the \emph{exchange graph} of the exchange
pattern~$(\Sigma_T)$.

Theorem~\ref{th:patterns-on-SM} implies that the structural results
obtained in \cite[Theorem~5.6]{cats1} hold in full generality, for
arbitrary bordered surfaces with marked points:

\begin{corollary}
\label{cor:structural-clust}
Let $\Acal$ be a cluster algebra whose exchange matrices arise from
triangulations of a surface~$\SM$.
Then each seed in $\Acal$ is uniquely determined by its
cluster; the cluster complex (see \cite[Definition~5.4]{cats1}) and
the exchange graph $\Exch$ of $\Acal$ do not depend on the choice of
coefficients in~$\Acal$;
the seeds containing a given cluster variable form a connected
subgraph of\/~$\Exch$;
and several cluster variables appear together in the same cluster if
and only if every pair among them~does.
The cluster complex is the complex of tagged arcs, as in
\cite[Theorem~7.11]{cats1};
it is the clique complex for its $1$-skeleton, and
is a connected pseudomanifold.
\end{corollary}

Theorem~\ref{th:patterns-on-SM} and
Corollary~\ref{cor:structural-clust} are proved in
Chapter~\ref{sec:laminated-teich}, using results and
constructions from the intervening chapters.
The proof is based on interpreting the
cluster variables $x_\gamma$ as generalized \emph{lambda lengths},
which are particular functions on the appropriately defined extension
of the Teichm\"uller
space of~$\SM$.

\chapter{Lambda lengths on bordered surfaces with punctures}
\label{sec:hyperbolic-cluster}

The machinery of lambda lengths was introduced and developed by
R.~Penner~\cite{penner-decorated, penner-lambda} in his study of
decorated Teichm\"uller spaces.
In this chapter, we adapt Penner's definitions to the case at hand, and
give a couple of useful geometric lemmas.

Throughout the paper, $\SM$ is a bordered surface with marked
points as described at the beginning of
Chapter~\ref{sec:cats1-recap}.
The (cusped) \emph{Teichm\"uller space} $\Teich\SM$ consists of all
complete finite-area hyperbolic structures with constant curvature
$-1$ on $\Surf\setminus \Mark$, with geodesic boundary at
$\partial\Surf\setminus\Mark$, considered up to $\Diff_0\SM$,
diffeomorphisms of $\Surf$ fixing $\Mark$ that are homotopic to the
identity.
(Thus there is a cusp at each point of $\Mark$.)  Our assumptions on
$\SM$ guarantee that $\Teich\SM$ is non-empty.
In fact, it is a manifold of dimension $n-p=6g+3b+2p+c-6$ in the notation
of Definition~\ref{def:ideal-triang}.

For a given hyperbolic structure in~$\Teich\SM$,
each arc can be represented by a unique geodesic.
Since there are cusps at the marked points, such a geodesic segment is of
infinite length.
So if we want to measure the
``length'' of a geodesic arc between two marked points, we need to
renormalize. This is done as follows.

\begin{definition}[\emph{Decorated Teichm\"uller space}
\cite{penner-decorated, penner-bordered, penner-lambda}]
\label{def:decor-teich}
A~point in a decorated Teichm\"uller space $\dTeich\SM$ is a
hyperbolic structure as above 
together with a collection of horocycles $h_p$, one around each
cusp corresponding to a marked point $p \in \Mark$.
\end{definition}

Appropriately interpreted, a \emph{horocycle} around a cusp~$p$ is the set
of points at an equal distance from~$p$: although the cusp is
infinitely far away from any point in the surface, there is still a
well-defined way to compare the distance to~$p$ from two different points
in the surface.  A horocycle can also be characterized as a curve
perpendicular to every geodesic to~$p$.

\begin{remark}
Our definition is a common generalization of those explicitly
given by Penner in \emph{loc.\ cit.}, as we simultaneously decorate
both the punctures and the marked points on the boundary. The
possibility of extending his theory to
this generality was already mentioned by Penner \cite[comments following
  Theorem~5.10]{penner-lambda}.
\end{remark}

Recall that $\BSM$ denotes the set of segments of the
boundary~$\partial\Surf$ between two adjacent marked points.
The cardinality of $\BSM$ is thus equal to~$c$,
the number of marked points on~$\partial\Surf$.

\begin{definition}[\emph{Lambda lengths}
\cite{penner-decorated, penner-bordered, penner-lambda}]
\label{def:lambda-length}
Fix $\sigma\in\dTeich\SM$. Let $\gamma$ be an arc in~$\APSM$,
or a boundary segment in~$\BSM$.
We will use the notation $\gamma_\sigma$ for the geodesic
representative of~$\gamma$ (relative to~$\sigma$).

Let $l(\gamma)=l_\sigma(\gamma)$ be the
  signed distance along~$\gamma_\sigma$ between the horocycles at either end
  of~$\gamma$ (positive if the two horocycles do not intersect, negative if
  they do intersect).
  The \emph{lambda length}~$\lambda(\gamma)=\lambda_\sigma(\gamma)$
  of~$\gamma$ is defined
by\footnote{This definition coincides with the one in
    \cite[Section~4]{penner-lambda} (or~\cite{gsv2}), and
differs by a factor of~$\sqrt{2}$ from the
    definition in~\cite{penner-decorated}.
The choice made here makes Lemma~\ref{lem:horocyclic-segment} below work with no
    factors.}
\begin{equation}
  \lambda(\gamma) = \exp(l(\gamma)/2).
\end{equation}
\end{definition}

Definition~\ref{def:lambda-length} can also be interpreted in terms of
a certain
dot product between two null vectors corresponding to the two
endpoints.   See~\cite{penner-decorated} and
Remark~\ref{rem:geom-grassmannian}.

For a given $\gamma\in\APSM\cup\BSM$, one can view the lambda length
\[
\lambda(\gamma):\sigma\mapsto\lambda_\sigma(\gamma)
\]
as a function on  the decorated Teichm\"uller space~$\dTeich\SM$.
Penner shows that such lambda lengths can be used to coordinatize
$\dTeich\SM$, as follows.

\begin{theorem}
\label{th:decorated-coord}
  For any triangulation $T$ of $\SM$, the map
  \[
  \prod_{\gamma \in T\cup\BSM} \lambda(\gamma) : \dTeich\SM \to
  \RR_{>0}^{n+c}
  \]
  is a homeomorphism.
\end{theorem}
(Recall from~\eqref{eq:n=6g+3b+3p+c-6} that $n$
is the total number of  arcs in~$T$, and $c$ is the number of marked
points on the boundary.)

\begin{remark}
\label{rem:penner-5.10}
  The first version of this theorem was proved by Penner
  \cite[Theorem~3.1]{penner-decorated}, which treats the case of
  closed surfaces with punctures.  This was later extended
  \cite{penner-bordered,penner-lambda}; the most relevant statement
  for us is \cite[Theorem~5.10]{penner-lambda}, which is not quite the
  statement above, since there the punctures in the interior are not
  decorated by horocycles and therefore are treated differently.
  Theorem~\ref{th:decorated-coord} follows easily from the same
  arguments, for instance using the doubling argument of
  \cite[Theorem~5.10]{penner-lambda} to reduce it to the case of a
  closed surface.
\end{remark}

For our purposes, the crucial property of lambda lengths
is the ``Ptolemy relation,''
the basic prototype of an exchange relation
in cluster algebras.

\begin{figure}[ht]
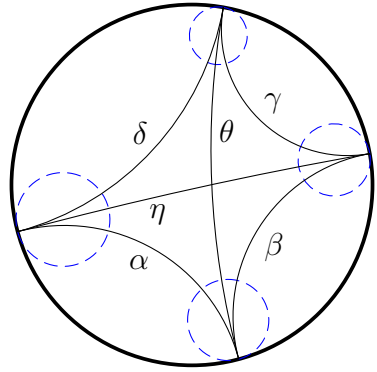



\[
\mfig{draws/hyp-disk.0}
\]
\caption{Sides and diagonals in a hyperbolic ideal quadrilateral,
  drawn in a lift to the hyperbolic plane}
\label{fig:ptolemy}
\end{figure}

\begin{proposition}[{Ptolemy relations}
{\cite[Proposition 2.6(a)]{penner-decorated}}]
\label{pr:geometric-flip}
Let
\[
\alpha,\beta,\gamma,\delta\in\APSM\cup\BSM
\]
be arcs or boundary segments (not necessarily distinct) that cut out a
quadrilateral in~$S$; we assume that the sides of the quadrilateral, listed in
cyclic order, are $\alpha,\beta,\gamma,\delta$.
Let $\eta$ and $\theta$ be the two diagonals of this quadrilateral;
see Figure~\ref{fig:ptolemy}.
Then the corresponding lambda lengths satisfy the Ptolemy relation
  \begin{equation}
\label{eq:ptolemy-lambda}
  \lambda(\eta) \lambda(\theta) = \lambda(\alpha)\lambda(\gamma)
  + \lambda(\beta)\lambda(\delta).
  \end{equation}
\end{proposition}

There is a Ptolemy relation \eqref{eq:ptolemy-lambda}
associated to each ordinary flip in an
ideal triangulation (cf.\ Definition~\ref{def:flip}).
Note that some sides of the relevant quadrilateral may be glued to
each other, changing the appearance of the relation.
See for example Figure~\ref{fig:ptolemy-4-punct-sphere}.

\begin{figure}[htbp]
\begin{center}

\setlength{\unitlength}{8pt}
\begin{picture}(30,10.5)(2.5,0)

\thinlines

\qbezier(10,0)(7,7)(10,7)
\qbezier(10,0)(13,7)(10,7)

\qbezier(25,0)(22,7)(25,7)
\qbezier(25,0)(28,7)(25,7)

\put(10,0){\circle*{0.5}}
\put(10,5){\circle*{0.5}}
\put(25,0){\circle*{0.5}}
\put(25,5){\circle*{0.5}}

\put(10,0){\line(1,0){15}}

\qbezier(10,0)(40,20)(25,0)
\qbezier(10,0)(-5,20)(25,0)

\put(8,7){\makebox(0,0){$\alpha$}}
\put(17.5,1){\makebox(0,0){$\beta$}}
\put(27,7){\makebox(0,0){$\gamma$}}
\put(21,8){\makebox(0,0){$\theta$}}
\put(14,8){\makebox(0,0){$\eta$}}

\end{picture}
\end{center}

\[
  \lambda(\eta) \lambda(\theta) = \lambda(\alpha)\lambda(\gamma)
  + \lambda(\beta)^2
\]
\caption{Ptolemy relation on a $4$-punctured sphere.
Cf.\ Figure~\ref{fig:triang-4-punct-sphere}.}
\label{fig:ptolemy-4-punct-sphere}
\end{figure}

By Theorem~\ref{th:decorated-coord},
each triangulation provides a set of coordinates on $\dTeich\SM$, while
Proposition~\ref{pr:geometric-flip} allows us to relate the
coordinatizations corresponding to different triangulations.
In the absence of punctures, this leads to an exchange pattern
(with a special choice of coefficients) in
which the lambda lengths play the role of cluster variables;
cf.\ \cite{fg-dual-teich, gsv2}.
For a punctured surface, the situation is more delicate, for
reasons both geometric and combinatorial:
as we know, not every arc can be flipped without leaving the realm of
ordinary triangulations.
There is also a (related) algebraic reason, provided by the following
lemma, a special case of Proposition~\ref{pr:geometric-flip}.

\begin{figure}[htbp]
\begin{center}
\setlength{\unitlength}{2pt}
\begin{picture}(80,48)(0,-4)
\thinlines
\put(40,20){\circle*{2}}
\put(0,20){\circle*{2}}
\put(80,20){\circle*{2}}

\put(0,20){\line(1,0){80}}

\qbezier(25,20)(25,-5)(80,20)
\qbezier(25,20)(25,45)(80,20)

\qbezier(0,20)(20,40)(40,40)
\qbezier(0,20)(20,0)(40,0)
\qbezier(40,0)(60,0)(80,20)
\qbezier(40,40)(60,40)(80,20)

\put(40,43){\makebox(0,0){$\alpha$}}
\put(40,-4){\makebox(0,0){$\beta$}}
\put(15,23){\makebox(0,0){$\theta$}}
\put(23,12){\makebox(0,0){$\eta$}}
\put(50,23){\makebox(0,0){$\gamma$}}

\put(40,16){\makebox(0,0){$p$}}

\end{picture}
\end{center}
\caption{Arcs in a punctured digon.}
\label{fig:lambda-digon}
\end{figure}
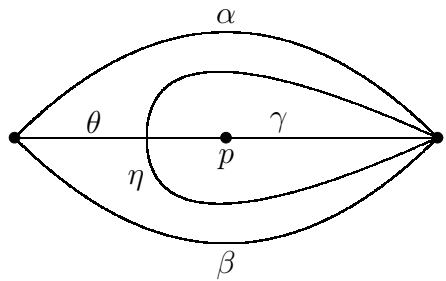

\begin{corollary}
\label{cor:digon-problematic}
Let $\alpha,\beta,\gamma,\eta,\theta\in\APSM\cup\BSM$ be as shown in
Figure~\ref{fig:lambda-digon}, that~is:
$\alpha$ and $\beta$ bound a digon with a sole puncture~$p$ inside it;
$\theta$ and $\gamma$ connect~$p$ to the vertices of the
digon; $\eta$~is the loop enclosing~$\gamma$.
Then
\begin{equation}
\label{eq:ptolemy-problematic}
\lambda(\eta) \lambda(\theta)
= \lambda(\alpha)\lambda(\gamma) + \lambda(\beta)\lambda(\gamma).
\end{equation}
\end{corollary}

The Ptolemy relation \eqref{eq:ptolemy-problematic}
cannot be an instance of a cluster exchange~\eqref{eq:exchange-rel-xx}
since the two terms on the right-hand side
of~\eqref{eq:ptolemy-problematic}
have a common factor~$\lambda(\gamma)$.
Thus Corollary~\ref{cor:digon-problematic} shows that
in the punctured case,
complexities associated with setting up a cluster algebra structure
already arise for ordinary flips,
namely those that create self-folded triangles.

This issue can be resolved by introducing tagged arcs and their lambda
lengths, and by extending the
Ptolemy relations to the case of tagged flips.
In the case of a tagged arc with a notched end,
our definition of a lambda length will require the notion of
(the hyperbolic distance from) a \emph{conjugate horocycle}.

\begin{definition}
For an horocycle~$h$ around a puncture in the interior of~$\Surf$, we
denote by~$L(h)$  the length of~$h$ as a (non-geodesic) curve in the
hyperbolic metric.
Two horocycles~$h$ and~$\bar h$ around the same interior marked point
are called \emph{conjugate} if $L(h) L(\bar h) = 1$.
\end{definition}

\pagebreak[3]

\begin{lemma}[{\cite[Lemma~4.4]{penner-lambda}}]
\label{lem:horocyclic-segment}
Fix a decorated hyperbolic structure in $\dTeich\SM$.
Consider a triangle in $\Surf\!\setminus\!\Mark$ with vertices
$p,q,r\!\in\!\Mark$
whose sides have lambda lengths $\lambda_{pq}$, $\lambda_{pr}$, and
$\lambda_{qr}$.
Then the length~$L_r$ of the horocyclic
  segment cut out by the triangle at vertex $r$ is given by
  \[
  L_r = \frac{\lambda_{pq}}{\lambda_{pr}\lambda_{qr}}.
  \]
\end{lemma}

\begin{figure}[ht]
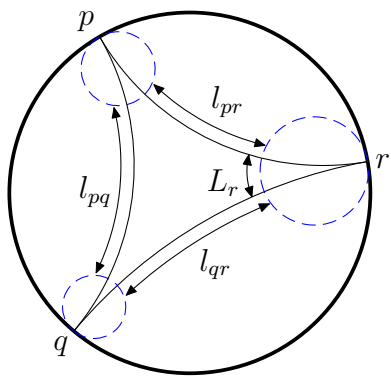

  \[
  \mfigb{draws/hyp-disk.10}
  \]
  \caption{The lengths in the statement of Lemma~\ref{lem:horocyclic-segment}.}
  \label{fig:horo-seg-statement}
\end{figure}

\begin{figure}[ht]
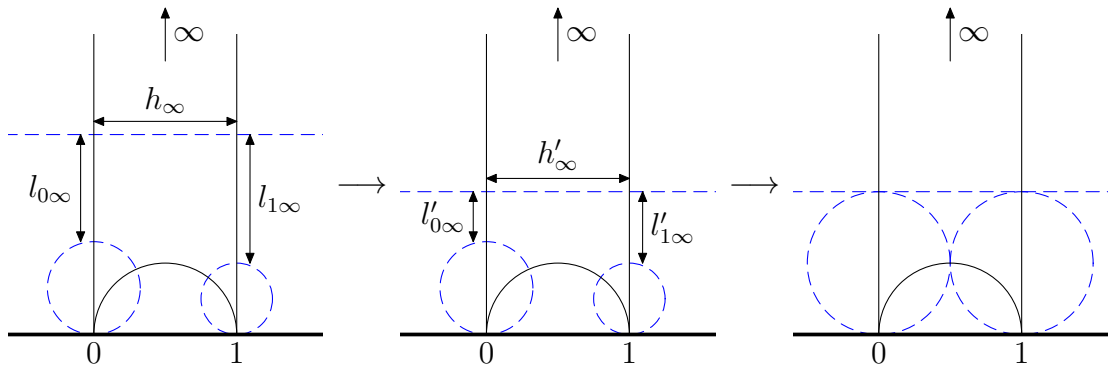

\[
\mfigb{draws/hyp.0}
    \longrightarrow\mfigb{draws/hyp.1}
    \longrightarrow\mfigb{draws/hyp.2}
\]
  \caption{Moving a horocycle in the proof of
    Lemma~\ref{lem:horocyclic-segment}.
At the last step, all the hyperbolic lengths
    $l_{xy}$ become~0.}
  \label{fig:horocyclic-segment}
\end{figure}

\begin{proof}
  First let us see how the side lengths~$l_{ij}$, the
  lambda lengths~$\lambda_{ij}$, and the horocyclic lengths~$L_i$
  change when we change the choice of horocycle.  For convenience, we
  work in the upper-half-plane model for the hyperbolic plane, with
 the metric
  \[
  ds^2 = \frac{dx^2 + dy^2}{y^2},
  \]
  and put the three vertices at $0$, $1$, and $\infty$, as in
  Figure~\ref{fig:horocyclic-segment} on the left.  Let us move the
  horocycle around~$\infty$ from an initial Euclidean height of $y$ to
  a height of $y'$.  By elementary integration, the new lengths are
 $L_\infty' = L_\infty (y/y')$,
 $l_{i\infty}' = l_{i\infty} + \ln(y'/y)$, and
  $\lambda_{i\infty}' = \lambda_{i\infty}\sqrt{y'/y}$,
  for $i \in \{0,1\}$.  The other lambda lengths and horocyclic
  lengths are unchanged.  By symmetry, similar statements are true with
  $\{0,1,\infty\}$ permuted.

  From this we see that, up to scale,
 $\dfrac{\lambda_{01}}{\lambda_{0\infty}\lambda_{1\infty}}$
  is the unique expression in the $\lambda_{ij}$ that is covariant
in the same way as $L_\infty$
  with respect to the $\RR_{>0}^3$-action associated with moving the
 three horocycles.  To fix the scale,
  consider the case where all three horocycles just touch, as
in Figure~\ref{fig:horocyclic-segment} on the right.  In this
  case, $L_\infty=1$, as in the statement of the lemma.
\end{proof}

\begin{figure}[ht]
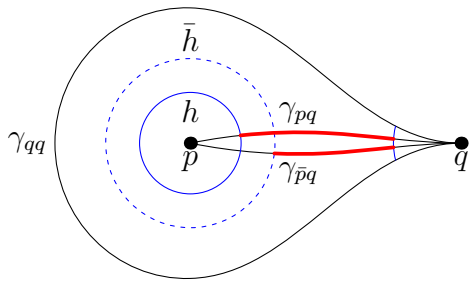

\begin{center}
$\mfigb{draws/geom-cluster.20}$
  \caption{The punctured monogon in
    Lemma~\ref{lem:conj-lambdas}.
The geodesic segments
whose lengths are measured are shown in red.}
\label{fig:conj-lambdas}
\end{center}
\end{figure}

\begin{lemma}
\label{lem:conj-lambdas}
Consider a punctured monogon with the vertex $q\in\partial\Surf$
and a sole puncture~$p$ in the interior.
Choose a horocycle around~$q$, and a horocycle $h$ around~$p$.
Let $\lambda_{qq}$ and $\lambda_{pq}$
be the corresponding lambda lengths for the boundary
of the monogon and the arc $\gamma_{pq}$ connecting $p$ and $q$ inside it,
respectively.
Let $\bar h$ be the horocycle around $p$ conjugate to~$h$,
and let $\lambda_{\bar p q}$ be the corresponding lambda length
of~$\gamma_{pq}$;
see Figure~\ref{fig:conj-lambdas}.
Then $\lambda_{qq} = \lambda_{pq} \lambda_{\bar p q}$.
\end{lemma}

\begin{proof}
By Lemma~\ref{lem:horocyclic-segment} applied to the self-folded
triangle with sides $\gamma_{pq}$ (twice) and $\gamma_{qq}$, we have
$L(h)= \dfrac{\lambda_{qq}}{\lambda_{pq}^2}$; similarly,
$L(\bar h)=\dfrac{\lambda_{qq}}{\lambda_{\bar pq}^2}$.
Since $h$ and $\bar h$ are conjugate, we obtain
$1=L(h)L(\bar h)=\dfrac{\lambda_{qq}^2}{\lambda_{pq}^2\lambda_{\bar
    pq}^2}$, and the claim follows.
\end{proof}

\chapter{Lambda lengths of tagged arcs}
\label{sec:lambda-tagged}

Lemma~\ref{lem:conj-lambdas} can be used to define a cluster algebra
structure associated with the decorated Teichm\"uller space
$\dTeich\SM$ of a general bordered surface with punctures, extending
the construction in~\cite{fock-goncharov1, fock-goncharov2, gsv2}.
As mentioned earlier,
the key idea is to interpret a notched end of a tagged arc as an
indication that in defining the corresponding lambda length,
we should take the distance to the \emph{conjugate} horocycle.

\begin{definition}[\emph{Lambda lengths of tagged arcs}]
\label{def:lambda-length-tagged}
Fix a decorated hyperbolic structure $\sigma\in\dTeich\SM$.
The \emph{lambda length~$\lambda(\gamma)=\lambda_\sigma(\gamma)$}
of a tagged arc $\gamma\in\ATSM$ is defined as follows.
If both ends of $\gamma$ are tagged plain, then the definition of
$\lambda(\gamma)$ given in Definition~\ref{def:lambda-length} stands.
Otherwise, the definition should be adjusted by replacing each
horocycle $h_p$ at a notched end~$p$ of $\gamma$ by the corresponding
conjugate horocycle~$\bar h_p$.
\end{definition}

In order to write relations among these lambda lengths,
we will need the following lemma.

\begin{figure}[htbp]
\begin{center}
\setlength{\unitlength}{2.5pt}
\begin{picture}(40,32)(0,0)
\thinlines

\qbezier(20,0)(-10,32)(20,32)
\qbezier(20,0)(50,32)(20,32)

\qbezier(20,0)(25,18)(20,20)
\qbezier(20,0)(15,18)(20,20)
\put(22.4,16){\makebox(0,0){$\notch$}}

\multiput(20,0)(0,20){2}{\circle*{2}}
\put(2,24){\makebox(0,0){$\eta$}}
\put(20,-4){\makebox(0,0){$q$}}
\put(20,24){\makebox(0,0){$p$}}
\put(13.7,11){$\gamma$}
\put(23.4,11){$\gamma'$}
\end{picture}
\end{center}
\caption{Arcs $\gamma, \gamma'$ and the enclosing loop~$\eta$
in Lemma~\ref{lem:conj-lambdas-2}}
\label{fig:lambda-conj}
\end{figure}
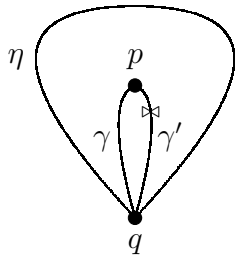

\begin{lemma}
\label{lem:conj-lambdas-2}
Let $\gamma$ and $\gamma'$ be two tagged arcs
connecting marked points $p,q\in\Mark$.
Assume that the untagged versions of $\gamma$ and $\gamma'$ coincide.
Also assume that $\gamma$ and $\gamma'$ have identical tags at~$q$,
and different tags at~$p$. (See Figure~\ref{fig:lambda-conj}.)
Let $\eta$ be the loop based at~$q$ wrapping around~$p$, so that
$\eta$ encloses a monogon with a sole puncture~$p$ inside it.
Then $\lambda(\eta)=\lambda(\gamma)\lambda(\gamma')$, where we compute
$\lambda(\eta)$ using $h_q$ or $\bar h_q$ according to whether
$\gamma$ and~$\gamma'$ are plain or notched at~$q$, respectively.
\end{lemma}

\begin{proof}
In view of Definition~\ref{def:lambda-length-tagged},
Lemma~\ref{lem:conj-lambdas-2} is a restatement of
Lemma~\ref{lem:conj-lambdas}.
\end{proof}

\begin{remark}
One delicate aspect associated with Lemma~\ref{lem:conj-lambdas-2} is
that $\eta$ itself is \emph{not} a legal tagged arc since it encloses a
once-punctured monogon.
Assume furthermore that both ends of $\eta$ are tagged plain.
Then $\eta$ is an arc in $\APSM$; as such, it is
represented by $\tau(\eta)=\gamma'$
(cf.\ Figure~\ref{fig:arc-as-tagged-arc}).
However, $\lambda(\eta)$ is not the same as~$\lambda(\gamma')$.
\end{remark}

Lemma~\ref{lem:conj-lambdas-2} allows us to write the exchange
relations associated with the tagged flips
of types (\ref{sec:cats1-recap}.D) and (\ref{sec:cats1-recap}.Q)
described in Remark~\ref{rem:two-types-of-flips}.

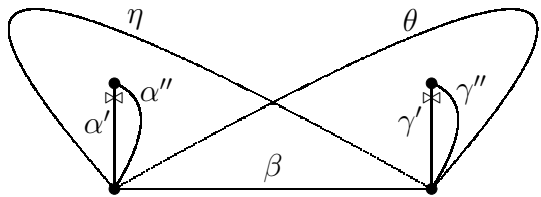
\begin{figure}[htbp]
\begin{center}

\setlength{\unitlength}{8pt}
\begin{picture}(30,8.5)(2.5,0)

\thinlines

\put(10,0){\line(0,1){5}}
\qbezier(10,5)(12.5,4)(10,0)
\put(10,4.3){\makebox(0,0){$\atag$}}

\put(25,0){\line(0,1){5}}
\qbezier(25,5)(27.5,4)(25,0)
\put(25,4.3){\makebox(0,0){$\atag$}}

\put(10,0){\circle*{0.5}}
\put(10,5){\circle*{0.5}}
\put(25,0){\circle*{0.5}}
\put(25,5){\circle*{0.5}}

\put(10,0){\line(1,0){15}}

\qbezier(10,0)(40,17)(25,0)
\qbezier(10,0)(-5,17)(25,0)
\put(9.2,3.2){\makebox(0,0){$\alpha'$}}
\put(12,4.7){\makebox(0,0){$\alpha''$}}
\put(17.5,1){\makebox(0,0){$\beta$}}
\put(24,3.2){\makebox(0,0){$\gamma'$}}
\put(26.9,4.7){\makebox(0,0){$\gamma''$}}
\put(24,8){\makebox(0,0){$\theta$}}
\put(11,8){\makebox(0,0){$\eta$}}

\end{picture}
\end{center}

\[
  \lambda(\eta) \lambda(\theta)
= \lambda(\alpha')\lambda(\alpha'')\lambda(\gamma')\lambda(\gamma'')
  + \lambda(\beta)^2
\]
\caption{Ptolemy relation for a tagged flip in a $4$-punctured
  sphere.
Cf.\ Figures~\ref{fig:tagged-triang-4-punct-sphere}
and~\ref{fig:ptolemy-4-punct-sphere}.}
\label{fig:ptolemy-tagged-4-punct-sphere}
\end{figure}

\begin{definition}[\emph{Ptolemy relations for tagged arcs}]
\label{def:ptolemy-tagged}
If two tagged triangulations $T_1$ and $T_2$ are related by a flip
of type~(\ref{sec:cats1-recap}.Q),
then the corresponding lambda
lengths are related by an appropriate specialization
of the equation~\eqref{eq:ptolemy-lambda}.
We will continue to refer to such relations among
lambda lengths of tagged arcs as (generalized) \emph{Ptolemy relations}.
Note that these relations can be more complicated than their
counterparts for the ordinary arcs:
some of the arcs $\alpha, \beta,\gamma,\delta$
appearing in~\eqref{eq:ptolemy-lambda} may bound a once-punctured
monogon and so might not be present in $T_1$ and~$T_2$.
In such a case, following Lemma~\ref{lem:conj-lambdas-2} we should
replace the lambda length of each such loop by the product of
lambda lengths of the two tagged arcs in $T_1$ (equivalently,~$T_2$)
that it encloses.
See Figure~\ref{fig:ptolemy-tagged-4-punct-sphere} for an example.
\end{definition}

\begin{definition}[\emph{Digon relations}]
  \label{def:digon-rel}
For a tagged flip of type~(\ref{sec:cats1-recap}.D),
consider the punctured digon as in
Corollary~\ref{cor:digon-problematic} and Figure~\ref{fig:lambda-digon}.
Making use of
Lemma~\ref{lem:conj-lambdas-2},
we can rewrite the relation~\eqref{eq:ptolemy-problematic}
in the form
\begin{equation}
\label{eq:exchange-digon}
\lambda(\gamma')\lambda(\theta) = \lambda(\alpha) + \lambda(\beta),
\end{equation}
where $\gamma'$ denotes~$\gamma$ with a notch at~$p$, as in
Figure~\ref{fig:lambda-conj}.  We refer to this as a \emph{digon
  relation}.  As in the case of Ptolemy relations, if one of the sides
of the digon ($\alpha$ or~$\beta$) bounds a punctured monogon, we
replace the lambda length of the corresponding loop by the product of
the lambda lengths of the two tagged arcs it encloses.
\end{definition}

Our next goal is to show that lambda lengths of tagged arcs on
a given bordered surface with marked points naturally form a normalized
exchange pattern whose exchange relations are
the relations of
Definitions~\ref{def:ptolemy-tagged} and~\ref{def:digon-rel}.
Making these statements precise will require a bit of preparation.

First, the coefficient semifield $\PP$ is going to be the tropical
semifield (see Definition~\ref{def:tropical})
generated by the lambda lengths of the boundary segments
(cf.\ Remark~\ref{rem:trop-coeffs-from-boundary}):
\begin{equation}
\label{eq:P-trop-boundary}
\PP=\Trop(\lambda(\gamma):\gamma\in\BSM).
\end{equation}
(If $\Surf$ is closed, then $\PP=\{1\}$ is the trivial one-element
semifield.)
We note that formula~\eqref{eq:P-trop-boundary} makes sense in view of
Theorem~\ref{th:decorated-coord}, which enables us to treat these
lambda lengths as independent variables.

Second, let us describe the clusters.
As in the case of ordinary arcs,
each lambda length of a tagged arc $\gamma\in\ATSM$ can be viewed as a
function $\sigma\mapsto\lambda_\sigma(\gamma)$ on the
decorated Teichm\"uller space $\dTeich\SM$.
For a tagged triangulation $T$ of~$\SM$,~let
\begin{equation}
\label{eq:xx=lambda}
\xx(T)=\{\lambda(\gamma) : \gamma \in T\}
\end{equation}
denote the collection of lambda lengths of the tagged arcs
in~$T$.

Third, the ambient field.
Let us pick an ordinary triangulation $T_\circ$ without self-folded
triangles.
It follows from Theorem~\ref{th:decorated-coord} that the lambda
lengths in $\xx(T_\circ)$ are algebraically independent over the field of
fractions of~$\PP$.
Let $\Fcal=\Fcal(T_\circ)$ be the field generated (say over~$\RR$)
by these lambda lengths.
That is, $\Fcal$ consists of all functions defined
on a dense subset of $\dTeich\SM$ which can be written as a rational
expression (say with real coefficients) in the lambda lengths of the
arcs in~$T_\circ$.

\begin{theorem}
\label{th:cluster-lambda}
There exists a unique normalized exchange pattern $(\Sigma_T)$, here
identified with its positive realization (see
Definition~\ref{def:positive-realization} and
Proposition~\ref{prop:positive-realization}) with the following
properties:
\begin{itemize}
\item
the coefficient semifield $\PP$ is the tropical semifield
 generated by the lambda lengths of boundary
segments, as in~\eqref{eq:P-trop-boundary};
\item
the ambient field $\Fcal=\Fcal(T_\circ)$ is generated over $\PP$ by
the lambda lengths of a given triangulation~$T_\circ$ with no
self-folded triangles;
\item
the underlying $n$-regular graph $\Exch\!=\!\ESM$ is the exchange graph
of~tagged triangulations (see Definition~\ref{def:ESM}), 
and the seeds 
$\Sigma_T\!=\!(\xx(T),\pp(T),B(T))$
are labeled by the vertices of~$\ESM$,
as in Theorem~\ref{th:patterns-on-SM};
\item
each cluster $\xx(T)$ consists of the lambda lengths of the tagged
arcs in~$T$, as in~\eqref{eq:xx=lambda};
\item
each exchange matrix $B(T)$ is the signed adjacency matrix of~$T$;
\item
the exchange relations out of each seed $\Sigma_T$
are the Ptolemy relations (see Definition~\ref{def:ptolemy-tagged})
and the digon relations (see Definition~\ref{def:digon-rel})
associated with the two respective types of tagged flips from~$T$
(cf.\ Remark~\ref{rem:two-types-of-flips}).
\end{itemize}
Neither the ambient field $\Fcal(T_\circ)$
nor the entire exchange pattern $(\Sigma_T)$
depend on the choice of the initial triangulation~$T_\circ$.
\end{theorem}

\begin{proof}
We know that the signed adjacency matrices $B(T)$ associated
with tagged triangulations~$T$ satisfy the mutation
rule~\eqref{eq:B-mut-flip}, as required in the definition of an
exchange pattern.
We also know from Theorem~\ref{th:decorated-coord} that the lambda
lengths forming the initial cluster $\xx(T_\circ)$ are algebraically
independent.
It remains to verify that
\begin{itemize}
\item[(i)]
the relations \eqref{eq:ptolemy-lambda}
and~\eqref{eq:exchange-digon}
associated with arbitrary tagged flips
can be viewed as exchange relations
\eqref{eq:exchange-rel-E} for the signed adjacency matrices~$B(T)$
of tagged triangulations~$T$,
\item[(ii)]
the coefficients appearing in these relations satisfy the mutation rules
\eqref{eq:p-mutation1-E}--\eqref{eq:p-mutation2-E}, and
\item[(iii)]
the normalization condition \eqref{eq:normalization} holds in the
tropical semifield~$\PP$.
\end{itemize}
Straightforward albeit somewhat tedious details of these verifications
are omitted.
It helps to note that a statement essentially equivalent to claim (ii) has
been already checked in Remark~\ref{rem:trop-coeffs-from-boundary}.
\end{proof}

\begin{remark}
It is tempting to try to deduce Theorem~\ref{th:patterns-on-SM} from
Theorem~\ref{th:cluster-lambda} by expressing cluster variables for
any exchange pattern with exchange matrices $B(T)$ as lambda
lengths of tagged arcs, perhaps rescaled to get different
coefficients. It turns out however that this simplistic approach does
not produce the most general coefficient patterns, as required for
Theorem~\ref{th:patterns-on-SM}.
Instead, we will need to develop, in subsequent chapters,
a more complicated concept of generalized lambda lengths for
\emph{laminated Teichm\"uller spaces} associated with
\emph{opened surfaces}.
\end{remark}

\begin{remark}
\label{rem:param-teich-tagged}
We note that Theorem~\ref{th:cluster-lambda} implies that
lambda lengths of tagged arcs in any cluster (i.e., in any tagged
triangulation) parametrize the decorated Teichm\"uller space
$\dTeich\SM$, extending Theorem~\ref{th:decorated-coord} verbatim to
the case of tagged triangulations. 
\end{remark}

\begin{example}
Let $\SM$ be a once-punctured digon, with notation as in
Figure~\ref{fig:lambda-digon}.
The coefficient semifield is $\PP=\Trop(\lambda(\alpha),\lambda(\beta))$.
The four cluster variables
$\lambda(\gamma)$, $\lambda(\theta)$,
$\lambda(\gamma')$, $\lambda(\theta')$,
are labeled by the tagged arcs in~$\SM$.
The four clusters correspond to the four tagged triangulations, cf.\
Figure~\ref{fig:tagged-cplx-once-punct-digon}.
The two exchange relations have the
form~\eqref{eq:exchange-digon}.
The resulting exchange pattern has finite type $A_1\times
A_1$, in the nomenclature of~\cite{ca2}.
\end{example}

In the case of surfaces with no punctures, there is no tagging, and
Theorem~\ref{th:cluster-lambda} specializes to its counterparts given
by V.~Fock and A.~Goncharov~\cite{fock-goncharov1, fock-goncharov2}
and by M.~Gekhtman, M.~Shapiro, and A.~Vainshtein~\cite{gsv2}.
The case of an unpunctured disk discussed in
Example~\ref{ex:unpunctured-disk} below was already treated
in~\cite[Section~12.2]{ca2}, without the hyperbolic geometry
interpretation.

\begin{example}
\label{ex:unpunctured-disk}
Let $\SM$ be an unpunctured $(n+3)$-gon with the vertices
$v_1,\dots,v_{n+3}$, labeled counterclockwise.
For $1\!\le\! i\!<\!j\!\le\! n\!+\!3$,
let $\gamma_{ij}$ denote the arc or boundary segment connecting $v_i$
and~$v_j$, that is, a diagonal or a side of the $(n+3)$-gon.
Denote $\lambda_{ij}=\lambda(\gamma_{ij})$.
Applying  the construction in Theorem~\ref{th:cluster-lambda}
to this special case, we get the coefficient semifield
\[
\PP=\Trop(\lambda_{12},\lambda_{23},\dots,\lambda_{n+3,1})
\]
generated by the lambda lengths of the sides of the $(n+3)$-gon;
the cluster variables are the lambda lengths of diagonals.
The corresponding cluster algebra (of type~$A_n$) can be interpreted (see
\cite[Proposition~12.7]{ca2}) as a homogeneous coordinate ring of the
Grassmannian $\operatorname{Gr}_{2,n+3}$ of $2$-dimensional subspaces
in~$\CC^{n+3}$.
See Example~\ref{example:G2n+3} for a more detailed treatment.
\end{example}

In the case of a once-punctured disk, we recover a particular cluster
algebra of type~$D_n$ that has been described (from a different
perspective) in \cite[Section~12.4]{ca2}.

\begin{example}
\label{ex:once-punctured-disk}
Let $\SM$ be an $n$-gon ($n\ge 3$) with vertices $v_1,\dots,v_n$
(labeled counterclockwise) and a single puncture~$p$ inside it. 
For $1\le i<j\le n$, there are two arcs or boundary
segments connecting $v_i$ and~$v_j$, depending on which side of the
curve the puncture~$p$ is on.
Let $\gamma_{ij}$ (resp.,  $\gamma_{i\bar j}=\gamma_{\bar j i}$)
denote the curve that has $p$ on the left (resp., right)
as we move from $v_i$ to~$v_j$.
There are also plain arcs $\gamma_{i\bar i}$ connecting $v_i$ to~$p$,
and tagged arcs $\tilde\gamma_{i\bar i}$ that have a
notched end at~$p$.
Replace $\gamma$'s with $\lambda$'s to denote the corresponding lambda
lengths. Then $\lambda_{12},\dots,\lambda_{n,1}$ generate the tropical
semifield of coefficients;
the remaining $\lambda$'s are cluster variables. As
always, clusters correspond to tagged triangulations.
The resulting cluster algebra coincides with the cluster algebra
$\mathcal{A}_\circ$ described in \cite[Example~12.15]{ca2},
and identified in \cite[Proposition~12.16]{ca2} with the
coordinate ring of the affine cone over the Schubert divisor in
$\operatorname{Gr}_{2,n+2}$.
\end{example}

\chapter{Opened surfaces}
\label{sec:opened-surfaces}

As shown in Chapter~\ref{sec:lambda-tagged}, the lambda lengths
of tagged arcs form an exchange pattern.
It is important to note that the coefficients in such an exchange
pattern are of a very special kind: they are monomials in the
lambda lengths of the boundary segments.
(For example, in the case of a closed surface with punctures, the
coefficients are trivial.)
In order to construct exchange patterns with general coefficients (as
in Theorem~\ref{th:patterns-on-SM}), we will need to modify our
geometric setting, extending the Teichm\"uller space from surfaces
with cusps at marked points to \emph{opened surfaces}.

It is not unusual in Teichm\"uller theory to
allow both cusped surfaces and surfaces with
geodesic boundary in the same moduli and Teichm\"uller spaces.
One standard model is the space of all complex structures on the
complement of the marked points.  A complex structure on the
neighborhood of a singularity can have two possible behaviours: it
can have a removable singularity at the marked point, corresponding to
a cusp in the hyperbolic structure after uniformization; or it can be
equivalent to the complex plane minus a closed disk, corresponding to
a non-finite volume end after uniformization.
It will be more convenient for us to use a different
(equivalent) model:
instead of complex structures, we will work with hyperbolic metrics,
truncated so that they have geodesic boundary and finite volume.
In addition, we add an orientation on each geodesic boundary component.

On the combinatorial/topological level, our construction will be based
on the following concept.

\begin{definition}[\emph{Opening of a surface}]
\label{def:opened-surface}
Let $\barM = \Mark \setminus\partial\Surf$ denote the set of
punctures of~$\Surf$.
For a subset $P \subset \barM$,
the corresponding \emph{opened surface}~$\Surf_P$
is obtained from $\Surf$ by removing a small open disk
  around each point in~$P$.  For $p \in P$, let $C_p$ be the
  boundary component of $\Surf_P$ created in this way.
We then introduce a new marked point $M_p$ on each component~$C_p$, and set
\[
\Mark_P=(\Mark\setminus P) \cup \{M_p\}_{p\in P}\,,
\]
creating a new bordered surface with marked points
$(\Surf_P,\Mark_P)$.
The sets of marked points $\Mark_P$ and $\Mark$ can be identified with
each other in a natural way.
See Figure~\ref{fig:opening}.
\end{definition}

\begin{figure}[htbp]
\begin{center}
\setlength{\unitlength}{1.8pt}
\begin{picture}(50,50)(-15,-10)
\thicklines
  \put(0,10){\circle*{3}}
  \put(40,10){\circle*{3}}
  \put(-20,10){\circle*{3}}
  \put(60,10){\circle*{3}}

\put(-5,10){\makebox(0,0){$p$}}
\put(45,10){\makebox(0,0){$q$}}
\put(-25,10){\makebox(0,0){$a$}}
\put(65,10){\makebox(0,0){$b$}}

\qbezier(-20,10)(-20,30)(20,30)
\qbezier(-20,10)(-20,-10)(20,-10)
\qbezier(60,10)(60,30)(20,30)
\qbezier(60,10)(60,-10)(20,-10)

\put(20,37){\makebox(0,0){$\SM$}}

\end{picture}
\qquad\qquad\qquad\qquad\qquad\qquad
\begin{picture}(50,50)(0,-10)
\thicklines
  \put(0,10){\circle{20}}
\put(40,10){\circle{20}}
  \put(0,20){\circle*{3}}
  \put(40,20){\circle*{3}}
  \put(-20,10){\circle*{3}}
  \put(60,10){\circle*{3}}

\qbezier(-20,10)(-20,30)(20,30)
\qbezier(-20,10)(-20,-10)(20,-10)
\qbezier(60,10)(60,30)(20,30)
\qbezier(60,10)(60,-10)(20,-10)

\put(6,24){\makebox(0,0){$M_p$}}
\put(34,24){\makebox(0,0){$M_q$}}
\put(12,1){\makebox(0,0){$C_p$}}
\put(29,1){\makebox(0,0){$C_q$}}

\put(-25,10){\makebox(0,0){$a$}}
\put(65,10){\makebox(0,0){$b$}}

\put(20,37){\makebox(0,0){$(\Surf_P,\Mark_P)$}}
\end{picture}
\end{center}
\caption{Opening of a surface.
Here $\SM$ is a twice-punc\-tured digon,
$\Mark\!=\!\{a,b,p,q\}$,
$P\!=\!\barM\!=\!\{p,q\}$, $\Mark_P\!=\!\{a,b,M_p,M_q\}$.}
\label{fig:opening}
\end{figure}

\begin{remark}
  Be careful to distinguish a surface with an opening from a surface
  with an extra boundary component with one marked
  point.  We will consider different Teichm\"uller spaces in the two cases,
  and treat them rather differently.
\end{remark}

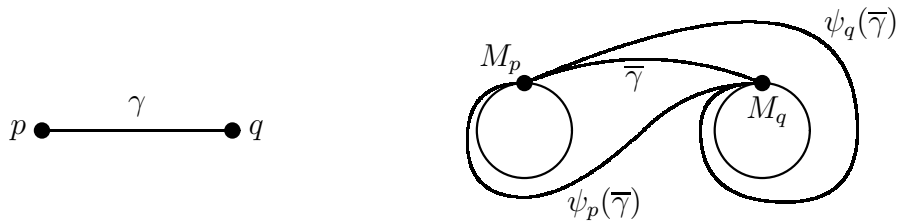
\begin{figure}[htbp]
\begin{center}
\setlength{\unitlength}{1.8pt}
\begin{picture}(60,20)(0,-5)
\thicklines
  \put(0,10){\line(1,0){40}}
  \put(0,10){\circle*{3}}
  \put(40,10){\circle*{3}}

\put(-5,10){\makebox(0,0){$p$}}
\put(45,10){\makebox(0,0){$q$}}
\put(20,15){\makebox(0,0){$\gamma$}}
\end{picture}
\qquad\qquad\qquad
\begin{picture}(60,40)(0,-5)
\thicklines
  \put(0,10){\circle{20}}
\put(50,10){\circle{20}}
  \put(0,20){\circle*{3}}
  \put(50,20){\circle*{3}}

\qbezier(0,20)(25,30)(50,20)

\qbezier(0,20)(70,50)(70,10)
\qbezier(55,-5)(70,-5)(70,10)
\qbezier(37,10)(37,-5)(55,-5)
\qbezier(37,10)(37,20)(50,20)

\qbezier(0,20)(-12,20)(-12,8)
\qbezier(-12,8)(-12,-4)(0,-4)
\qbezier(0,-4)(10,-4)(25,10)
\qbezier(25,10)(35,20)(50,20)

\put(-5,25){\makebox(0,0){$M_p$}}
\put(51,14){\makebox(0,0){$M_q$}}

\put(23,21){\makebox(0,0){$\bargamma$}}
\put(9,-7){{$\psi_p(\bargamma)$}}
\put(63,31){{$\psi_q(\bargamma)$}}

\end{picture}
\end{center}
\caption{An arc connecting two punctures, and three of its lifts}
\label{fig:lift}
\end{figure}

There is a natural ``projection'' map
\begin{equation}
\label{eq:collapse-arcs}
\varkappa_P:\AP(\Surf_P,\Mark_P)\to\APSM
\end{equation}
(surjective but not injective)
that corresponds to collapsing the new boundary components~$C_p$.
We will refer to any $\bargamma\in\AP(\Surf_P,\Mark_P)$ that projects
onto a given arc $\gamma \in \APSM$
as a \emph{lift} of~$\gamma$.
To describe these lifts, we introduce, for every $p \in P$, the map
\begin{equation}
\label{eq:clockwise-twist}
\psi_p:\AP(\Surf_P,\Mark_P) \rightarrow
\AP(\Surf_P,\Mark_P)
\end{equation}
that takes each arc ending at $M_p$ and twists it once clockwise
around~$C_p$ (with a negative Dehn twist).
Then, for example, an arc~$\gamma\in\APSM$ connecting two distinct punctures
$p, q \in P$ has the lifts
\begin{equation}
\label{eq:lifts-pq}
\varkappa_P^{-1}(\gamma)=\{(\psi_p)^n(\psi_q)^m\bargamma\}_{n, m \in \ZZ}\,,
\end{equation}
where $\bargamma$ is some particular lift of~$\gamma$.
See Figure~\ref{fig:lift}.
If~$\gamma\in\APSM$ goes from a puncture $p \in P$ back to itself,
then $\varkappa_P^{-1}(\gamma)$ consists of
two orbits under the action of~$\psi_p\,$.

Opening all the punctures in $\barM$ results in the ``largest'' opened
surface
\begin{equation}
\label{eq:barSM}
\barSM=(\Surf_{\barM}, \Mark_{\barM}).
\end{equation}
Its arc complex
\begin{equation}
\label{eq:largest-arcs}
\barAPSM\stackrel{\rm def}{=}\AP(\Surf_{\barM}, \Mark_{\barM})
\end{equation}
naturally projects onto all the other arc complexes
$\AP(\Surf_P,\Mark_P)$;
that is, the map $\varkappa_{\barM}$ factors through every other map
$\varkappa_P$, for $P\subset\barM$.

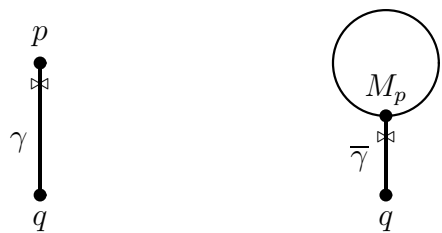
\begin{figure}[htbp]
\begin{center}
\setlength{\unitlength}{2pt}
\begin{picture}(40,35)(0,-3)
\thicklines
\put(20,0){\line(0,1){25}}
\put(20,21){\makebox(0,0){$\notch$}}
\put(16,10){\makebox(0,0){$\gamma$}}

\multiput(20,0)(0,25){2}{\circle*{2}}
\put(20,-5){\makebox(0,0){$q$}}
\put(20,30){\makebox(0,0){$p$}}
\end{picture}
\qquad\qquad
\begin{picture}(40,35)(0,-3)
\thicklines
\put(20,0){\line(0,1){15}}
\put(20,11){\makebox(0,0){$\notch$}}
\put(15,7){\makebox(0,0){$\bargamma$}}

\multiput(20,0)(0,15){2}{\circle*{2}}
\put(20,25){\circle{20}}
\put(20,-5){\makebox(0,0){$q$}}
\put(20,20){\makebox(0,0){$M_p$}}
\end{picture}

\end{center}
\caption{Lift of a tagged arc}
\label{fig:lifting-tagged-arc}
\end{figure}

\begin{definition}[\emph{Lifts of tagged arcs}]
\label{def:lift-tagged-arc}
To lift a tagged arc $\gamma\in\ATSM$ to an opened surface
$(\Surf_P,\Mark_P)$ (in particular, to $\barSM$),
we simply lift the untagged version of~$\gamma$
to $\AP(\Surf_P,\Mark_P)$ (resp., $\barAPSM$), and then affix
the same tags as the ones used at the corresponding ends of~$\gamma$.
See Figure~\ref{fig:lifting-tagged-arc}.
Thus, the lifted tagged arc $\bargamma$ may have a notched end at an
unopened point $p\notin P$, or at a marked point $M_p$ ($p\in P$).
We denote by $\AT(\Surf_P,\Mark_P)$ (resp., $\barATSM$)
the set of all such tagged arcs $\bargamma$ on $(\Surf_P,\Mark_P)$
(resp., $\barSM$).  Note that as in $\AT\SM$, tagged arcs
that enclose a monogon containing a single puncture~$p$ are forbidden,
whether or not $p \in P$.
\end{definition}

\chapter{Lambda lengths on opened surfaces}
\label{sec:lambda-opened}

We are now prepared to describe our main Teichm\"uller-theoretic
construction. This will be done in two steps,
Definitions~\ref{def:partial-Teich} and~\ref{def:complete-decorated}.

\begin{definition}
  \label{def:dec-punct}
  A \emph{decorated set of marked points} $\dPunct$ is a subset $P \subset
  \barM$ of the punctures, together with a choice of
  orientation on $C_p$ for each $p \in P$; this orientation can be
  clockwise or counterclockwise.
\end{definition}

\begin{definition}
\label{def:partial-Teich}
  Fix a decorated set of marked points $\dPunct$.
  We define the \emph{partially opened Teichm\"uller
    space} $\Teich_{\dPunct}(\Surf_P,\Mark_P)$ as the space of all
  finite-volume, complete
  hyperbolic metrics on $\Surf_P \setminus (\Mark\setminus P)$ with geodesic
  boundary, modulo isotopy.
  For a decorated set of marked points $\dPunct$, 
  the \emph{decorated partially opened
    Teichm\"uller space $\dTeich_{\dPunct}(\Surf_P,\Mark_P)$} is the same set of
  metrics as in $\Teich_{\dPunct}(\Surf_P,\Mark_P)$, modulo isotopy relative
  to~$\{M_p\}_{p\in P}$
and with a choice of horocycle around each point in $\Mark
  \setminus P$.
\end{definition}

(The orientations on the
boundary will be used shortly.)

That is, a hyperbolic structure in $\Teich_{\dPunct}(\Surf_P,\Mark_P)$ has
a cusp at each point in~$\Mark\setminus P$
(i.e., at each original marked point on~$\partial\Surf$ and at each
puncture in~$\barM\setminus P$),
and a new circular geodesic boundary component $C_p$ arising from each point
$p\in P$.  The boundary is otherwise geodesic.  Note in particular
that there are no cusps at the points $M_p$, which are
introduced merely to help parametrize the new boundary.

\begin{remark}
\label{rem:dTeich-lift}
The 
space $\dTeich_{\dPunct}(\Surf_P,\Mark_P)$
is a fibration over $\Teich_{\dPunct}(\Surf_P,\Mark_P)$ with
fibers $\mathbb{R}^{\Mark}$.  The decorations look different depending on the
marked point: for points not in $P$, the decoration is a choice of
horocycle as in the previous chapters, while at the new geodesic boundary
the decoration comes from restricting the isotopies to those that
leave the boundary components $C_p$ ($p\in P$) fixed.  The difference
(isotopies of $\Surf_P$ that are isotopic to the identity, but not
while fixing the~$C_p$) is isotopies that twist the surface around the
$C_p$.  Such isotopies have a single real parameter for each $p\in P$,
namely the amount of twisting.
\end{remark}

\begin{figure}[ht]
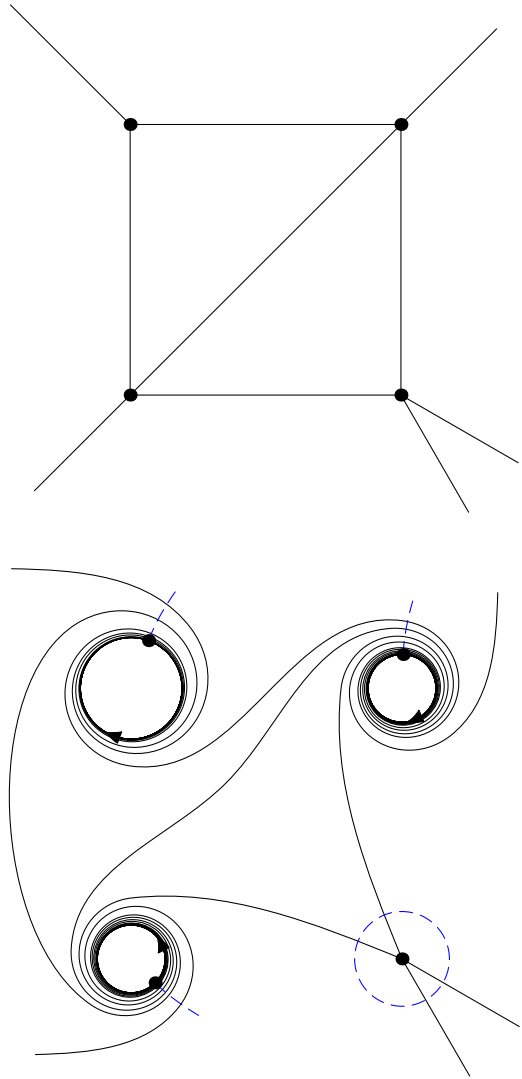

\[
\begin{array}{c}
  \mfigb{draws/geom-cluster.0}
\\[1.5in]
\mfigb{draws/geom-cluster.1}
\end{array}
  \]
  \caption{Representing arcs by geodesics on an opened surface.
Shown 
above
is a portion of the original surface; 
below is
a particular opening of the surface, endowed with a hyperbolic
structure and an orientation on the opened boundary components.
The lower right marked point, the only one not in~$P$,
has been left as an interior
    cusp; the remaining three have been opened into circular geodesic
    boundary components.
    The lower left component is oriented counterclockwise and the other
    two are oriented clockwise.
}
\label{fig:opening-surface}
\end{figure}

Given a decorated set of marked points $\dPunct$ and a
geometric structure
$\sigma\in\Teich_{\dPunct}(\Surf_P,\allowbreak\Mark_P)$ to each
arc $\gamma\in \APSM$, we can associate a unique infinite,
non-self\/intersecting geodesic $\gamma_\sigma$ on~$\Surf_P$
(geodesic with respect to~$\sigma$):
at endpoints of $\gamma$ that
are not opened,
the geodesic $\gamma_\sigma$ runs out to the cusp, while at
endpoints that are in~$P$, it spirals around $C_p$ in the
chosen direction.
See Figure~\ref{fig:opening-surface}.

For an arc $\bargamma \in \AP(\Surf_P, \Mark_P)$ on an opened
surface, we set $\bargamma_\sigma=\gamma_\sigma$,
where $\gamma=\varkappa_P(\bargamma)$ is obtained from $\bargamma$ by 
the collapsing map $\varkappa_P$  of~\eqref{eq:collapse-arcs}.
Thus the geodesic representative $\bargamma_\sigma$ does not depend on
how much $\bargamma$ winds around its opened ends.

We now coordinatize the Teichm\"uller spaces
$\dTeich_{\dPunct}(\Surf_P,\Mark_P)$ by introducing appropriate generalizations
of Penner's lambda lengths. This will require the following notion.

\pagebreak[3]

\begin{definition}
Fix a decorated set of marked points $\dPunct$ and a geometric structure
$\sigma\in\dTeich_{\dPunct}(\Surf_P,\Mark_P)$.
For each $p \in P$,
there is a
\emph{perpendicular horocyclic segment}~$h_p$ near~$C_p$: we take a
(short) segment of the horocycle from $M_p \in C_p$ which is
perpendicular to~$C_p$ and to all geodesics $\gamma_\sigma$ that spiral
to~$C_p$ in the direction given by the chosen orientation.  In
Figure~\ref{fig:opening-surface}, the horocycle segments $h_p$ are drawn as
dashed curves, as is the horocycle decorating the marked point
not in~$P$.
\end{definition}

\begin{remark}
  The perpendicular horocyclic segment can be obtained by following
  the \emph{horocyclic flow} from $M_p$ perpendicular to the boundary.
  (The horocyclic flow is similar to the geodesic flow, but follows
  the unique horocycle through a given point in a given direction.)
  It can be thought of as the set of points an equal distance to the
  ideal point obtained by following the boundary~$C_p$ infinitely far
  in the direction of its orientation.  See, e.g.,
  \cite{hedlund-fuchsian} for more on the horocyclic flow.
\end{remark}

\begin{definition}[\emph{Lambda lengths on an opened surface}]
\label{def:lambda-opened}
We next define lamb\-da lengths~$\lambda_\sigma(\gamma)$
for arcs $\gamma \in \AP(\Surf_P, \Mark_P)$ and for 
$\sigma\in\dTeich_P(\Surf_P,\Mark_P)$.
As in Definition~\ref{def:lambda-length},
we set
\begin{equation}
\label{eq:lambda-opened}
\lambda(\gamma)=\lambda_\sigma(\gamma)=e^{l(\gamma)/2},
\end{equation}
where $l(\gamma)=l_\sigma(\gamma)$ is the distance between appropriate
intersections of the geodesic~$\gamma_\sigma$ with the horocycles at
its two ends.
At ends of $\gamma_\sigma$ that spiral around one of the openings~$C_p$,
there will be many intersections between $\gamma_\sigma$ and the
horocyclic segment~$h_p$, and we need to pick one of them.
Assume that $\gamma$ connects two ends $M_p$ and~$M_q$,
with both $p$ and $q$ in~$P$ (this is the most complicated case).
Suppose furthermore that $\gamma$ twists sufficiently far in the
direction of the orientation of the boundary.
Then there are unique intersections between $\gamma_\sigma$ and
each of $h_p$ and~$h_q$
such that the path that runs
\begin{itemize}
\item along $h_p$ from $M_p$ to one intersection, then
\item along $\gamma_\sigma$ to the other intersection, then
\item along $h_q$ to the other endpoint $M_q$
\end{itemize}
is homotopic to the original arc~$\gamma$, as shown in
Figure~\ref{fig:horocycle-intersection}.  If one or both of the ends
of $\gamma$ are not in $P$, we leave $\gamma_\sigma$ unmodified at that end and
pick the unique intersection between $\gamma_\sigma$ and the corresponding
horocycle.  In either case, $l(\gamma)$ is the (signed) distance along
$\gamma_\sigma$ between the
chosen intersections with the two horocycles.

\begin{figure}[ht]
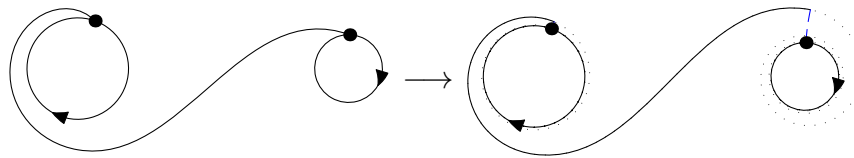

  \[
  \mfigb{draws/geom-cluster.5} \longrightarrow \mfigb{draws/geom-cluster.6}
  \]
  \caption{Finding the correct intersection with the perpendicular
    horocycles.}
  \label{fig:horocycle-intersection}
\end{figure}

\pagebreak[3]

In order to extend the definition to \emph{all} arcs
$\gamma\in\AP(\Surf_P,\Mark_P)$, not
just those that twist sufficiently much,
we postulate how $l(\gamma)$ and $\lambda(\gamma)$ change when we
twist $\gamma$ around the boundary.
Specifically, we mandate that
\begin{equation}
\label{eq:l-twist}
l(\psi_p\gamma) = n_p(\gamma)\,l(p) + l(\gamma)
\end{equation}
where 
\begin{align}
\label{eq:psi_p}
&\text{$\psi_p$ is the clockwise twist defined
by~\eqref{eq:clockwise-twist},}
\\
\label{eq:n_p}
&\text{$n_p(\gamma)$ is the number of ends of $\gamma$ that
  touch~$M_p$, and} 
\\
\label{eq:l(p)}
&\text{$l(p) =
\begin{cases}
  -\text{length of $C_p$}&\text{if $p \in P$ and $C_p$ is oriented counterclockwise;}\\
  0& \text{if $p \not\in P$;}\\
  \text{length of $C_p$}&\text{if $p \in P$ and $C_p$ is oriented clockwise.}
\end{cases}
$}
\end{align}
Accordingly (cf.~\eqref{eq:lambda-opened}), we have
\begin{equation}
\label{eq:lambda-twist}
\lambda(\psi_p\gamma) = \lambda(p)^{n_p(\gamma)}\,\lambda(\gamma), 
\end{equation}
where
\begin{equation}
\label{eq:lambda(p)}
\lambda(p) = \lambda_\sigma(p) = e^{l(p)/2}.
\end{equation}
In order for $l(\gamma)$ and $\lambda(\gamma)$ to be well defined, we
need of course to check that the requirements
\eqref{eq:l-twist}--\eqref{eq:lambda-twist} are consistent
with the earlier definitions given in the case where $\gamma$ twists
sufficiently much.  This follows from the following lemma.
\end{definition}

\begin{lemma}
  The distance along a geodesic $\gamma_\sigma$ between successive
  intersections with the horocycle $h_p$ is always equal to $|l(p)|$.
\end{lemma}

\begin{figure}[ht]
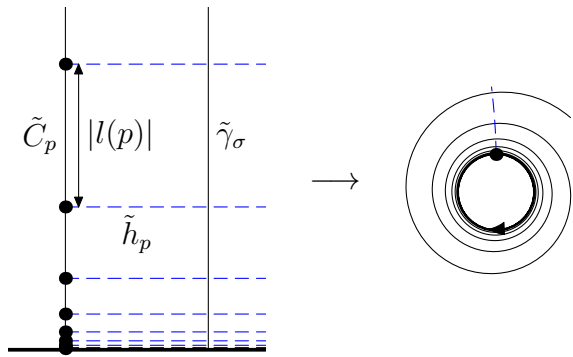

  \[
  \mfigb{draws/hyp.30} \quad\longrightarrow\quad \mfigb{draws/geom-cluster.30}
  \]
  \caption{Shown on the right is a geodesic $\gamma_\sigma$ spiraling to a boundary
    component~$C_p$, and the corresponding horocycle~$h_p$. 
On the left, their respective universal covers $\tilde\gamma_\sigma$, $\tilde
C_p$, and~$\tilde h_p$, in the upper-half-space model of $\HH^2$.
The distance between horocycles is the length of~$C_p$.}
  \label{fig:twist-difference}
\end{figure}

\begin{proof}
  Consider a segment~$s$ of $\gamma_\sigma$ between successive
  intersections with $h_p$ which is very close to $C_p$.  Because of
  the spiraling nature of $\gamma_\sigma$, as the distance from $s$ to
  $C_p$ approaches~$0$, the length of~$s$
  approaches the length of~$C_p$.  But since $\gamma_\sigma$ is part
  of a family of geodesics perpendicular to $h_p$, we can move~$s$
  within the family of geodesics without changing the length.
  Now move $s$ out from $C_p$ until it coincides with the
  desired segment.

  Alternatively, a purely geometric proof is sketched in
  Figure~\ref{fig:twist-difference}.  Consider the universal cover
  $\HH^2$ of $\SM$ and one lift $\tilde C_p$ of $C_p$ within it.
  Place the endpoint of $\tilde C_p$ to which $\gamma_\sigma$ is
  spiraling at infinity in the upper-half-space model of~$\HH^2$.
  Then the lifts $\tilde h_p$ of $h_p$ appear as straight lines
  parallel to the real axis, and a lift $\tilde \gamma_\sigma$ of
  $\gamma_\sigma$ appears as a line parallel to the imaginary axis
  (and $\tilde C_p$).  The distance between successive intersections
  on $\tilde \gamma_\sigma$ is independent of the left-right position
  of $\tilde \gamma_\sigma$, and in particular it agrees with the
  distance along $\tilde C_p$, namely $|l(p)|$.
\end{proof}

\begin{definition}
\label{def:complete-decorated}
The \emph{complete decorated Teichm\"uller space} $\TTSM$ is the
disjoint union over all $3^{|\barM|}$ possible choices of decorated sets of
marked points $\dPunct$ of $\dTeich_{\dPunct}(\Surf_P,\Mark_P)$.

It remains to describe the topology on $\TTSM$.
For an arc $\gamma\in\barAPSM$ (see Definition~\ref{def:lift-tagged-arc}),
define the lambda length
\[
\lambda(\gamma):\TTSM\to\mathbb{R}
\]
on each stratum of $\TTSM$
by projecting $\gamma$ to the appropriate set
$\Arcs^{\circ}(\Surf_P,\Mark_P)$ and using the
construction above.
  The topology on $\TTSM$ (making it into a connected
  space) is the weakest in which $\lambda(\bargamma)$ is
  continuous for all lifted arcs $\bargamma\in\barAPSM$.
\end{definition}

\begin{lemma}
  \label{lem:ptolemy-opened}
  Inside a quadrilateral in $\barSM$ with sides $\baralpha$,
  $\barbeta$, $\bargamma$, and~$\bardelta$ and diagonals $\bareta$
  and~$\bartheta$ as in
  Proposition~\ref{pr:geometric-flip}, we have the Ptolemy relation
  \begin{equation}
    \label{eq:ptolemy-opened}
    \lambda(\bareta)\lambda(\bartheta)=\lambda(\baralpha)\lambda(\bargamma)+\lambda(\barbeta)\lambda(\bardelta).
  \end{equation}
\end{lemma}

Note that the arcs $\baralpha,\dots,\bartheta$ have to form a
quadrilateral in $\barSM$; it is not enough for their projections to
$\SM$ to form a quadrilateral.  In particular, if an arc appears twice
on the boundary of a quadrilateral in $\SM$, we may have to take two
different lifts of it to $\barSM$ in order for
Lemma~\ref{lem:ptolemy-opened} to apply.

\begin{proof}
  This is equivalent to Proposition~\ref{pr:geometric-flip}:
the geometry is identical to what we had before, once a lift
to the universal cover is made.
\end{proof}

\begin{proposition}
\label{pr:param-teich-opened}
  Let $T$ be an ideal triangulation of $\SM$ without self-folded triangles.
For each $\gamma \in T$, fix an arc
  $\bargamma \in \barAPSM$ that projects to~$\gamma$.
Then the map
  \[
\Phi = \Biggl(\prod_{p \in \barM} \lambda(p)\Biggr)
\times
\Biggl(\prod_{\beta\in\BSM} \lambda(\beta)\Biggr)
\times
\Biggl(\prod_{\gamma \in T} \lambda(\bargamma)\Biggr)
: \TTSM \to \RR_{>0}^{n+\abs{\Mark}}
  \]
  is a homeomorphism, where
  $n$ is the number of arcs in~$T$ as in
  formula~\eqref{eq:n=6g+3b+3p+c-6}.
\end{proposition}

\begin{proof}
  For any vector $\Lambda$ in $\RR_{>0}^{n+\abs{\Mark}}$, we can
  construct $\Phi^{-1}(\Lambda)$, the unique geometric structure in
  $\TTSM$ with the corresponding set of lambda lengths, as follows.
  First note that if the arc~$\gamma$ has distinct endpoints $p$ and~$q$
  in~$\barM$, then 
  for any alternate lift $\bargamma'$ of~$\gamma$ there are $n,m \in
  \ZZ$ so that
  \[
  \bargamma' = \psi_p^n\psi_q^m(\bargamma)
  \]
(cf.~\eqref{eq:lifts-pq}). 
  Then by equation~\eqref{eq:lambda-twist}, for any hyperbolic
  structure with these lambda coordinates, we have 
\begin{equation}
\label{eq:lambda-bargamma'}
  \lambda(\bargamma') = \lambda(p)^{n}\,\lambda(q)^{m}\,\lambda(\bargamma).
\end{equation}
  This equation and similar ones can then be used 
to compute lambda lengths of all lifts of the arcs in~$T$.

  Now for each ideal triangle in $T$ with sides $\gamma_1$,
  $\gamma_2$, $\gamma_3$, pick lifts $\bargamma_1'$, $\bargamma_2'$,
  $\bargamma_3'$ that form a triangle in $\barSM$.
  Then
  take a decorated ideal hyperbolic triangle (i.e., a triangle with
  choice of horocycles around each cusp) so that the lambda lengths 
  of the sides of the triangle match with the $\lambda(\bargamma_i')$
(computed using~\eqref{eq:lambda-bargamma'}),
  as in Penner's proof of Theorem~\ref{th:decorated-coord}.  (There is
  a unique decorated ideal hyperbolic triangle with given
  lambda lengths.)

  We next need to glue these triangles together to form a hyperbolic
  surface.  For each arc $\gamma \in T$, we have two different lifts
  $\bargamma'$, $\bargamma''$ of $\gamma$ coming from the two
  different triangles that have this arc as a side (or the two
  different sides of the same triangle, in case $\gamma$ is the
  repeated edge of a self-folded triangle).  Suppose that $\bargamma''
  = \psi_p^n \psi_q^m(\bargamma')$, where $p$ and $q$ are the
  endpoints of~$\gamma$ as before.  Then glue the two hyperbolic triangles
  so that the the horocycles around the vertex corresponding to~$p$
  are offset by $n\cdot l(p)$ and the horocycles around the vertex~$q$ are
  offset by $m\cdot l(q)$.  It is then elementary to verify that the
  resulting glued surface has a metric completion which is a surface
  with the desired lambda lengths, proving the
  surjectivity of~$\Phi$. Conversely, since each decorated ideal
  triangle is determined by its lambda lengths and the gluings between
  adjacent triangles are determined by the data, $\Phi$ is injective.

  By definition of the topology on $\TTSM$, the map $\Phi$ is continuous.  It
  remains to show that $\Phi^{-1}$ is continuous.  To do this, we must
  show that for an arbitrary arc $\baralpha\in\barAPSM$, the
  lambda length $\lambda(\baralpha)$ is a continuous function of the
  given coordinates.  Let $\alpha = \varkappa(\baralpha)$ be the arc
  in $\APSM$ corresponding to~$\baralpha$.  We can move from $T$ to a
  triangulation that contains $\alpha$ by a series of edge flips in
  quadrilaterals.  For each such flip, Lemma~\ref{lem:ptolemy-opened}
  lets us write the lambda length of one
  lift of the new diagonal in terms of lambda lengths of lifts of the
  old arcs.  (If we can write one lift of a given arc in terms of the
  given lambda coordinates, we can write all lifts in terms of these
  lambda coordinates by multiplying by appropriate powers of the
  $\lambda(p)$ for $p \in \barM$.)  We end up inductively writing
  $\lambda(\baralpha)$ as an algebraic function with non-zero denominator in
  terms of the original coordinates.  Thus each $\lambda(\baralpha)$
  is a continuous function when pulled back to $\RR^{n+|\Mark|}$, so
  by definition of the topology on $\TTSM$ it follows that $\Phi^{-1}$ is
  continuous.
\end{proof}

We next wish to extend Proposition~\ref{pr:param-teich-opened} to the
case of tagged triangulations, as in
Remark~\ref{rem:param-teich-tagged}.

\begin{definition}[\emph{Lambda lengths of tagged arcs on an opened surface}]
  \label{def:lambda-opened-tagged}
  For $\sigma\in\TTSM$ and $\gamma\in\ATSM$, define $\gamma_\sigma$ to
  be the unique infinite, non-self\/inter\-secting geodesic which
  at each notched end spirals \emph{against} the orientation chosen on
  $C_p$ and is otherwise as before.  For $p\in\barM$, set
  \begin{equation}
    \nu(p)=2\ln\abs{\lambda(p)-\lambda(p)^{-1}}.
  \end{equation}
Let $\overline{M_p}$ be the point on $C_p$ a (signed) distance of
  $\nu(p)$ from $M_p$ in the direction against the orientation
  of~$C_p$.  Define the \emph{conjugate perpendicular horocycle}
  $\overline h_p$ to be the horocycle passing through $\overline{M_p}$
  and perpendicular to $C_p$ and to all geodesics spiraling against the
  orientation on~$C_p$.  Finally, for $\bargamma\in\barATSM$, define
  $l(\bargamma)$ to be the length between intersections
  with horocycles as before, using the conjugate perpendicular
  horocycle for notched ends that meet an opened puncture.
  Specifically, if $\bargamma$ is notched at~$p$ and plain at~$q$,
  there is a unique path that is homotopic to $\bargamma$ and runs
  \begin{itemize}
  \item along $C_p$ from $M_p$ to $\overline{M_p}$ a distance of
    $\nu(p)$ against the orientation of~$C_p$, then
  \item along $\overline h_p$ from $\overline{M_p}$ to an intersection
    with $\bargamma_\sigma$, then
  \item along $\bargamma_\sigma$ to an intersection with $h_q$, then
  \item along $h_q$ to $M_q$.
  \end{itemize}
  There is a similar path if $\bargamma$ is notched at both ends.
  Set $\lambda(\bargamma)=e^{l(\bargamma)/2}$ as before.

  We will also allow the obvious extensions of $l(\bargamma)$ and
  $\lambda(\bargamma)$ to a version of tagged arcs on $\barSM$ which
  enclose punctured monogons (so are not in $\barATSM$).
\end{definition}

\begin{figure}[ht]
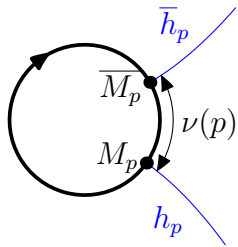

  \[
  \mfig{draws/geom-cluster.50}
  \]
  \caption{The effect of $\nu(p)$ on the horocycles.  Shown are the
    horocycle $h_p$ and the conjugate horocycle $\overline{h}_p$ in
    the case when $l(p) > 0$ (so $C_p$ is oriented clockwise), and
    $l(p)$ is large enough so that $\nu(p) > 0$.  (Here $l(p) \approx
    1.05$ so $\nu(p) \approx 0.2$.)}
  \label{fig:effect-nu}
\end{figure}

\begin{remark}
  \label{rem:nu-asymptotics}
  The correction term $\nu(p)$ is chosen so that
  Lemma~\ref{lem:holed-monogon} below comes out with no correction
  factors, which in turn implies that for each lifted arc
  $\bargamma\in\barATSM$, 
the lambda length $\lambda(\bargamma)$ is a continuous function on $\TTSM$.

  As $l(p)$ approaches $\pm \infty$, $\nu(p)$ is asymptotic to
  $\abs{l(p)}$, which amounts to saying that in the limit as $l(p)$
  gets large,
  $\overline{M_p}$ differs from $M_p$ by a full turn against the
  orientation on~$C_p$.  On the other hand, for $l(p)$ close to zero
  (when the boundary is close to a cusp),
  $\nu(p)$ is asymptotic to $2 \ln \abs{l(p)}$, which is large
  and negative.
\end{remark}

\begin{question}
  Is there a more geometrically natural way to define the conjugate
  perpendicular horocycle (Definition~\ref{def:lambda-opened-tagged})?
\end{question}

We next investigate the properties of the lambda lengths of tagged arcs.
In order to complete our construction of exchange patterns associated
with opened surfaces, we need to define exchange relations
involving
\begin{itemize}
\item[(i)]
the lambda lengths
$\lambda(\bargamma)$, for $\bargamma\in \barATSM$,
\item[(ii)]
the lambda lengths $\lambda(\beta)$, for $\beta\in\BSM$, and
\item[(iii)]
the lambda lengths $\lambda(p)$, for $p \in \barM$.
\end{itemize}
Some of these relations are easy to obtain.  We see right away that
the lambda lengths of types (i) and~(ii) still obey the Ptolemy
relation (Lemma~\ref{lem:ptolemy-opened}), provided the arcs form a
quadrilateral on $\barSM$ as before and the tags of the three arcs
meeting at each vertex of the quadrilateral agree with each other.

Next, there is change of the lift.  Equation~\eqref{eq:lambda-twist}
holds as before, with the convention that $n_p(\gamma)$
(cf.~\eqref{eq:n_p}) is a
\emph{signed} count: a plain end of $\gamma$ at~$p$ contributes~$+1$,
a notched end contributes~$-1$:
\begin{gather}
  \label{eq:signed-n_p}
  \text{$n_p(\gamma)$ is the signed number of ends of $\gamma$ that
    touch~$M_p$} \\
  \label{eq:signed-lambda-twist}
  \lambda(\psi_p\gamma) = \lambda(p)^{n_p(\gamma)}\,\lambda(\gamma)
\end{gather}

We will also need a relation associated with a tagged flip inside
an opened monogon, an analogue of Lemma~\ref{lem:conj-lambdas}.

\begin{figure}[ht]
\begin{center}
\setlength{\unitlength}{2.5pt}
\begin{picture}(60,30)(-20,8)
\thicklines
\put(0,10){\circle{20}}
\put(10,10){\circle*{2}}
\put(50,10){\circle*{2}}

\qbezier(-18,10)(-18,-35)(50,10)
\qbezier(-18,10)(-18,55)(50,10)

\put(-10,11.5){\makebox{\vector(0,1){0}}}

\thinlines
\qbezier(10,10)(25,-5)(50,10)
\qbezier(10,10)(25,25)(50,10)

\put(5,10){\makebox(0,0){$M_p$}}
\put(-7.5,21){\makebox(0,0){$C_p$}}
\put(53,10){\makebox(0,0){$q$}}


\put(12,12){\rotatebox{-50}{\makebox(0,0){$\notch$}}}

\put(20,1){\makebox(0,0){$\delta$}}
\put(20,19){\makebox(0,0){$\varrho$}}

\put(-21,10){\makebox(0,0){$\eta$}}

\end{picture}
\qquad\qquad\qquad
$\mfigb{draws/geom-cluster.40}$
\end{center}
  \caption{An opened monogon}
  \label{fig:opened-monogon}
\end{figure}
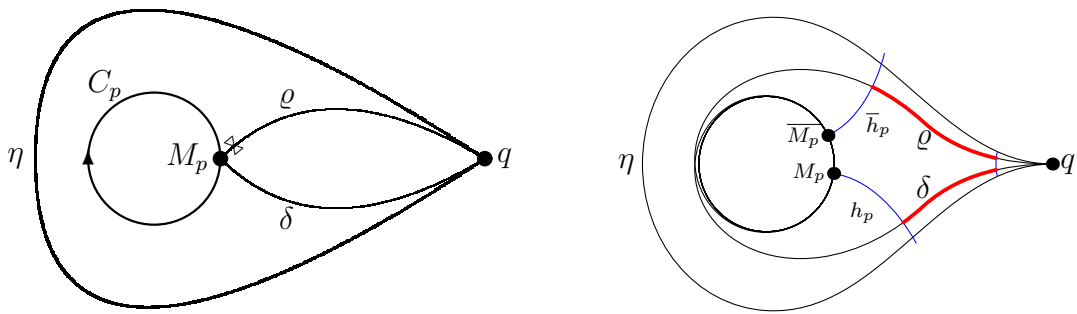

\begin{lemma}
\label{lem:holed-monogon}
Inside an opened surface, consider a monogon with a marked
vertex~$q$ and a single boundary component~$C_p$ in
  the interior.
Let $\delta$ and $\varrho$ be two compatible parallel 
tagged arcs in $\barATSM$ connecting $q$ and $M_p$, with $\delta$
plain and $\varrho$ notched at $M_p$, as shown in
Figure~\ref{fig:opened-monogon} on the left,
and let $\eta$ be the outer boundary of the monogon, tagged like
$\delta$ and $\varrho$ at~$q$.
  Then
\begin{equation}
\label{eq:holed-monogon}
\lambda(\delta)\lambda(\varrho) =
     \lambda(\eta).
\end{equation}
\end{lemma}

\pagebreak[3]

\begin{figure}[ht]
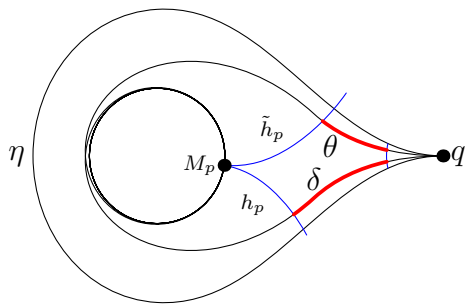

  \begin{center}
    $\mfigb{draws/geom-cluster.41}$
  \end{center}
  \caption{The opened monogon with alternate horocycles}
  \label{fig:opened-monogon-2}
\end{figure}

\begin{proof}
Let $\tilde h_p$ be the horocycle which is like $\overline
h_p$ but perpendicular to $C_p$ at $M_p$ instead of $\overline{M_p}$
as in Figure~\ref{fig:opened-monogon-2},
and let $\theta$ be the tagged arc like $\varrho$ but with lambda
length measured with respect to~$\tilde h_p$. Then
\eqref{eq:holed-monogon} is equivalent to
\begin{equation}
  \label{eq:holed-monogon-2}
  \lambda(\delta)\lambda(\theta) =
     \frac{\lambda(\eta)}{\abs{\lambda(p)-\lambda(p)^{-1}}}.
\end{equation}
Strange as it may seem at first glance, \eqref{eq:holed-monogon-2} is yet
another instance of the same Ptolemy relation.
To see that, suppose first that $C_p$ is oriented clockwise (so
that $\lambda(p)>1$, cf.\ \eqref{eq:l(p)},~\eqref{eq:lambda(p)}), and
consider Figure~\ref{fig:cover-monogon}, which on the top
shows lifts of the arcs $\delta$, $\theta$, and $\eta$ to the
universal cover of the monogon.  The bottom of
Figure~\ref{fig:cover-monogon} shows a different triple of lifts together
with lifts of the arcs $\delta' = \psi_p \delta$ and $\theta' =
\psi_p^{-1} \theta$.  

\begin{figure}[ht]
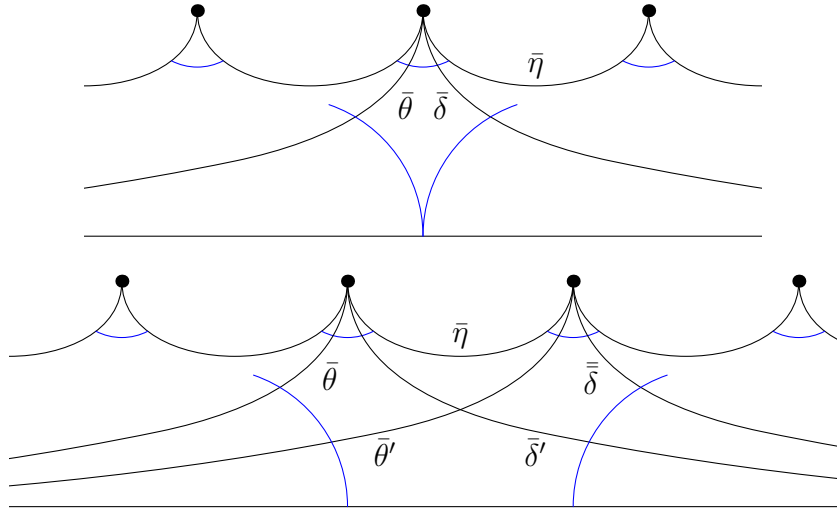

  \begin{center}
    $\mfigb{draws/geom-cluster.45}$\\[10pt]
    $\mfigb{draws/geom-cluster.46}$
  \end{center}
  \caption{The universal cover of the opened monogon, with different
    choices of lifts.}
  \label{fig:cover-monogon}
\end{figure}

\noindent
Applying the Ptolemy relation to the
quadrilateral with diagonals $\bar\delta'$ and~$\bar\theta'$, 
we~get
\begin{equation*}
\lambda(\delta')\lambda(\theta') = \lambda(\delta)\lambda(\theta) +
\lambda(p)\lambda(\eta).
\end{equation*}
Combining this with
\begin{align*}
  \lambda(\delta') &= \lambda(p)\,\lambda(\delta)\\
  \lambda(\theta') &= \lambda(p)\,\lambda(\theta)
\end{align*}
(from Equation~(\ref{eq:signed-lambda-twist})) we deduce
\[
(\lambda(p)^2-1)\lambda(\delta)\lambda(\theta) = \lambda(p)\,\lambda(\eta)
\]
as desired.
If $C_p$ is oriented counterclockwise instead, $\delta$ spirals
counterclockwise and $\varrho$ and $\theta$ spiral clockwise.  In this
case define
$\delta'=\psi_p^{-1} \delta$ and $\theta' = \psi_p \theta$.  Then
\begin{equation*}
\lambda(\delta')\lambda(\theta') = \lambda(\delta)\lambda(\theta) +
\lambda(\eta)\lambda(p)^{-1}
\end{equation*}
and we again deduce~\eqref{eq:holed-monogon-2}.
\end{proof}

Lemma~\ref{lem:holed-monogon} lets us find the relations
associated with the tagged flips of 
type~(\ref{sec:cats1-recap}.D).

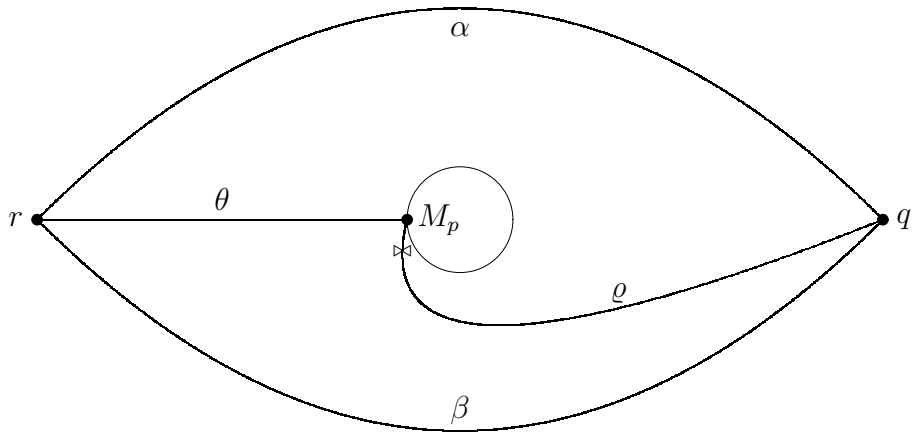
\begin{figure}[htbp]
\begin{center}
\setlength{\unitlength}{4pt}
\begin{picture}(80,40)(0,0)
\thinlines
\put(40,20){\circle{10}}
\put(35,20){\circle*{1}}
\put(0,20){\circle*{1}}
\put(80,20){\circle*{1}}

\put(0,20){\line(1,0){35}}
\qbezier(35,20)(30,0)(80,20)

\qbezier(0,20)(20,40)(40,40)
\qbezier(0,20)(20,0)(40,0)
\qbezier(40,0)(60,0)(80,20)
\qbezier(40,40)(60,40)(80,20)

\put(38,20){\makebox(0,0){$M_p$}}
\put(-2,20){\makebox(0,0){$r$}}
\put(82,20){\makebox(0,0){$q$}}
\put(34.55,17){\makebox(0,0){$\notch$}}

\put(40,38){\makebox(0,0){$\alpha$}}
\put(40,2){\makebox(0,0){$\beta$}}
\put(17.5,22){\makebox(0,0){$\theta$}}
\put(55,13){\makebox(0,0){$\varrho$}}


\end{picture}
\end{center}
\caption{An exchange relation in an opened digon.}
\label{fig:opened-digon}
\end{figure}

\begin{lemma}
\label{lem:exch-opened-digon}
Consider an opened digon with vertices $r$ and~$q$ and an 
opening $C_p$ with a marked point~$M_p$.
Let $\alpha$, $\beta$, $\varrho$, and~$\theta$ be the tagged arcs
shown in Figure~\ref{fig:opened-digon}.
(Possible tags at $r$ and $q$ have been suppressed in the picture.)
We assume that the arcs in $\{\alpha,\beta,\varrho,\theta\}$
are tagged so that any two of them are compatible,
with the exception of the pair $(\varrho,\theta)$.
Then
\begin{equation}
\label{eq:exch-opened-digon}
\lambda(\varrho)\, \lambda(\theta)=\lambda(\alpha)+\lambda(p)^{-1}\,\lambda(\beta).
\end{equation}
\end{lemma}

\begin{figure}[htbp]
\begin{center}
\setlength{\unitlength}{4pt}
\begin{picture}(80,40)(0,0)
\thinlines
\put(40,20){\circle{10}}
\put(35,20){\circle*{1}}
\put(0,20){\circle*{1}}
\put(80,20){\circle*{1}}

\put(0,20){\line(1,0){35}}
\qbezier(35,20)(30,0)(80,20)
\qbezier(35,20)(30,40)(80,20)

\qbezier(25,20)(25,-5)(80,20)
\qbezier(25,20)(25,45)(80,20)

\qbezier(0,20)(20,40)(40,40)
\qbezier(0,20)(20,0)(40,0)
\qbezier(40,0)(60,0)(80,20)
\qbezier(40,40)(60,40)(80,20)

\put(38,20){\makebox(0,0){$M_p$}}
\put(-2,20){\makebox(0,0){$r$}}
\put(82,20){\makebox(0,0){$q$}}

\put(40,38){\makebox(0,0){$\alpha$}}
\put(40,2){\makebox(0,0){$\beta$}}
\put(15,22){\makebox(0,0){$\theta$}}
\put(25,12){\makebox(0,0){$\eta$}}
\put(55,14){\makebox(0,0){$\gamma$}}
\put(55,26){\makebox(0,0){$\delta$}}


\end{picture}
\end{center}
\caption{Proof of Lemma~\ref{lem:exch-opened-digon}.}
\label{fig:opened-digon-proof}
\end{figure}
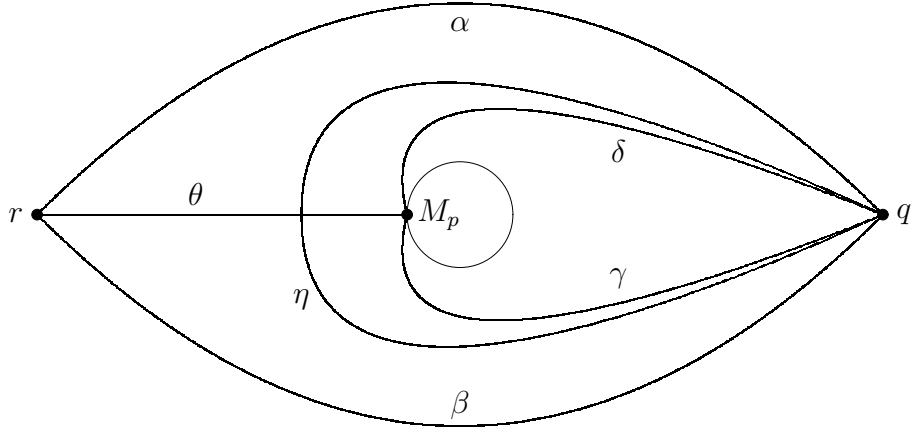

\begin{proof}
Let us introduce the arcs $\gamma$, $\delta$, and~$\eta$ as in
Figure~\ref{fig:opened-digon-proof} (all tagged plain at~$p$).
The lambda lengths of the six arcs in
Figure~\ref{fig:opened-digon-proof} satisfy the Ptolemy
relation~\eqref{eq:ptolemy-lambda}:
\[
\lambda(\theta)\lambda(\eta)=\lambda(\alpha)\lambda(\gamma)+\lambda(\beta)\lambda(\delta).
\]
We also have
\begin{equation}
\lambda(\delta)=\lambda(\gamma)\,\lambda(p)^{-1}
\end{equation}
(by~\eqref{eq:signed-lambda-twist}) and
\begin{equation}
\lambda(\eta)=\lambda(\gamma)\lambda(\varrho)
\end{equation}
(by Lemma~\ref{lem:holed-monogon}).
Putting everything together, we
obtain~\eqref{eq:exch-opened-digon}.
\end{proof}

\begin{corollary}
  \label{cor:param-teich-opened-top}
  Let $T$ be a tagged triangulation of $\SM$.
  For each $\gamma \in T$, fix an arc
  $\bargamma \in \barATSM$ that projects to~$\gamma$.
Then the map $\Phi$ as defined in Proposition~\ref{pr:param-teich-opened}
  is a homeomorphism.
\end{corollary}

\begin{proof}
  The triangulation~$T$ can be connected to an ideal
  triangulation $T'$ of $\SM$ with no notched arcs or self-folded
  triangles by a sequence of
  flips in quadrilaterals and digons.  At each step,
  Lemmas~\ref{lem:ptolemy-opened} and~\ref{lem:holed-monogon} let us
  express the lambda lengths after the flip in terms of those before
  the flip, so the map $\Phi$ above and the analogue defined with
  respect to $T'$ are related by a homeomorphism on the target.  But
  the latter map is a homeomorphism by
  Proposition~\ref{pr:param-teich-opened}.
\end{proof}

\begin{remark}
We note that while the definitions of lambda lengths on the opened
surface depend in an essential way on the chosen orientations of the
boundaries~$C_p\,$, the relations \eqref{eq:ptolemy-opened},
\eqref{eq:holed-monogon}, and~\eqref{eq:exch-opened-digon} that they
satisfy have the same form irrespective of the choices of orientations. 
\end{remark}

\chapter{Non-normalized exchange patterns from surfaces}
\label{sec:non-normalized-patterns-from-surfaces}

In this chapter, we describe a construction of a non-normalized exchange
pattern on the exchange graph $\ESM$ of tagged triangulations of the
original surface~$\SM$.
This construction is \emph{different} from the one given in
Chapter~\ref{sec:lambda-tagged}: although it is more complicated, it is
eventually going to provide---after proper rescaling---a more general
class of coefficients.
Here is the basic idea: rather than designating the lambda length of a
tagged arc $\gamma\in\ATSM$ as the corresponding
cluster variable,
we take the lambda length of an arbitrary lift of~$\gamma$
to the opened surface~$\barSM$ (see~\eqref{eq:barSM} and
Figure~\ref{fig:lift}).
It turns out that we do not have to coordinate these lifts:
the corresponding lambda lengths will always
form an exchange pattern.
In contrast to the simpler construction in
Chapter~\ref{sec:lambda-tagged}, this new exchange pattern will not be
normalized.

We begin by setting up the coefficient group $\PP=\PP\SM$
as the (free) abelian multiplicative group generated by the set
\begin{equation}
\label{eq:coeff-group-non-normalized}
\{\lambda(p): p \in \barM\} \cup \{\lambda(\beta):\beta\in\BSM\}
\end{equation}
of lambda lengths of boundary components $\beta$ and
opened circular components~$C_p$.
By Proposition~\ref{pr:param-teich-opened}, we can view (and treat)
these lambda lengths either as functions on the complete decorated
Teichm\"uller space $\TTSM$ or as formal variables (=coordinate
functions).

For each tagged arc $\gamma\in\ATSM$,
let us fix an arbitrary lift
$\bargamma\in\barATSM$ (see Definition~\ref{def:lift-tagged-arc}),
and set $x(\gamma)=\lambda(\bargamma)$.
Then, for each tagged triangulation $T\in\ESM$, define
\begin{equation}
\label{eq:xx(gamma)}
\xx(T)=\{\,x(\gamma):\gamma\in T\,\}.
\end{equation}
In view of Corollary~\ref{cor:param-teich-opened-top},
the rescaled lambda lengths in $\xx(T)$ can be treated as
formal variables algebraically independent
over the field of fractions of~$\PP\SM$.

We are now ready to state our next theorem:
the lambda lengths of lifts of tagged arcs form
a non-normalized exchange pattern.

\begin{theorem}
\label{th:cluster-lambda-opened}
For an arbitrary choice of lifts $\bargamma$ of tagged arcs
$\gamma\in\ATSM$, there exists a (unique) non-normalized exchange
pattern on $\ESM$ with the following properties:
\begin{itemize}
\item
the coefficient group is $\PP=\PP\SM$;
\item
the cluster variables are the lambda lengths $\lambda(\bargamma)$;
\item
the cluster $\xx(T)$ at a vertex $T\in\ESM$ is given by
\eqref{eq:xx(gamma)};
\item
the ambient field is generated over $\PP$ by some
(equivalently, any) cluster~$\xx(T)$;
\item
the exchange matrices are the signed adjacency matrices~$B(T)$; and
\item
the exchange relations out of each seed are the
relations \eqref{eq:ptolemy-opened} and
\eqref{eq:exch-opened-digon}
associated with the corresponding tagged flips,
properly rescaled via \eqref{eq:signed-lambda-twist} to reflect the
choices of lifts.
\end{itemize}
\end{theorem}

To be more accurate, the description of cluster variables above
should refer to a ``positive realization'' of the exchange pattern in
question, in the spirit of Definition~\ref{def:positive-realization}.
Even though this pattern is not of geometric type (as it is not
normalized), the corresponding
notions still have clear meaning, and the analogue of
Proposition~\ref{prop:positive-realization} holds.

\begin{proof}
The proof is similar to the proof of Theorem~\ref{th:cluster-lambda}.
As before, the real issue is coefficients:
we need to demonstrate that they satisfy the requisite mutation rules
\eqref{eq:p-mutation1}--\eqref{eq:p-mutation2}.
It is straightforward to check that these rules hold for each triple
of flips/mutations
\[
T_1 \stackrel{\mu_x}\longleftrightarrow
T_2 \stackrel{\mu_z}\longleftrightarrow
T_3 \stackrel{\mu_x}\longleftrightarrow
T_4\,,
\]
if the lifts of the arcs involved
are chosen in a coordinated way (this is essentially the same
verification as before)---and
therefore this rule would hold for any lifts, by
Proposition~\ref{pr:cv-rescaling}.
\end{proof}

\chapter{Laminations and shear coordinates}
\label{sec:shear}

In this chapter, we briefly review a small fragment---as this is all we
need---of W.~Thurston's theory of measured
laminations~\cite{thurston-travaux, thurston-bulletin},
and its relationship with matrix mutations.
Our exposition is an abridged adaptation of the one given by V.~Fock
and A.~Goncharov \cite[Section~3]{fg-dual-teich}.
An interested reader is referred to the cited sources for further details.

In the next chapter, we will extend these constructions to
the tagged setting.

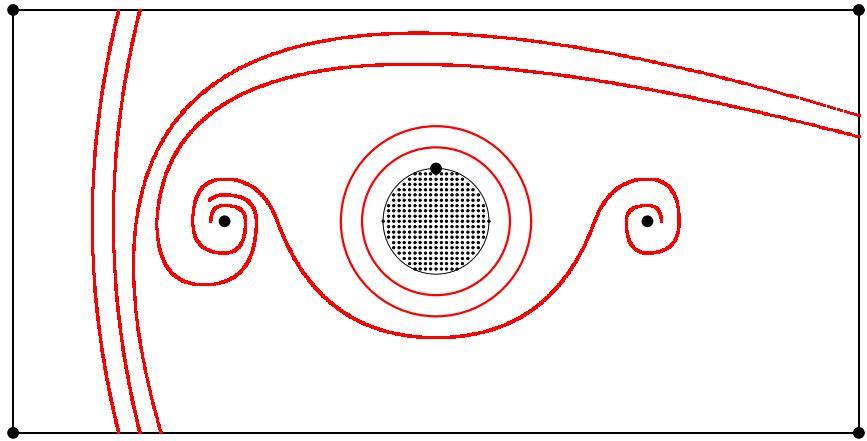
\begin{figure}[htbp]
\begin{center}
\setlength{\unitlength}{4pt}
\begin{picture}(80,40)(0,0)
\thinlines
\put(0,0){\line(1,0){80}}
\put(0,40){\line(1,0){80}}
\put(0,0){\line(0,1){40}}
\put(80,0){\line(0,1){40}}

  \put(0,0){\circle*{1}}
  \put(0,40){\circle*{1}}
  \put(80,0){\circle*{1}}
  \put(80,40){\circle*{1}}

  \put(40,20){\circle{10}}
\multiput(35,20)(0.5,0){21}{\circle*{0.2}}
\multiput(35.5,20.5)(0.5,0){19}{\circle*{0.2}}
\multiput(35.5,21)(0.5,0){19}{\circle*{0.2}}
\multiput(35.5,21.5)(0.5,0){19}{\circle*{0.2}}
\multiput(36,22)(0.5,0){17}{\circle*{0.2}}
\multiput(36,22.5)(0.5,0){17}{\circle*{0.2}}
\multiput(36.5,23)(0.5,0){15}{\circle*{0.2}}
\multiput(37,23.5)(0.5,0){13}{\circle*{0.2}}
\multiput(37.5,24)(0.5,0){11}{\circle*{0.2}}
\multiput(38,24.5)(0.5,0){9}{\circle*{0.2}}

\multiput(35.5,19.5)(0.5,0){19}{\circle*{0.2}}
\multiput(35.5,19)(0.5,0){19}{\circle*{0.2}}
\multiput(35.5,18.5)(0.5,0){19}{\circle*{0.2}}
\multiput(36,18)(0.5,0){17}{\circle*{0.2}}
\multiput(36,17.5)(0.5,0){17}{\circle*{0.2}}
\multiput(36.5,17)(0.5,0){15}{\circle*{0.2}}
\multiput(37,16.5)(0.5,0){13}{\circle*{0.2}}
\multiput(37.5,16)(0.5,0){11}{\circle*{0.2}}
\multiput(38,15.5)(0.5,0){9}{\circle*{0.2}}

  \put(40,25){\circle*{1}}

  \put(20,20){\circle*{1}}
  \put(60,20){\circle*{1}}

\thicklines
\darkred{\qbezier(10,0)(5,20)(10,40)}
\darkred{\qbezier(12,0)(7,20)(12,40)}
\darkred{\qbezier(14,0)(-2,55)(80,30)}
\darkred{\qbezier(13.6,20)(15,45)(80,28)}

\darkred{\qbezier(40,9)(51,9)(55,20)}
\darkred{\qbezier(40,9)(29,9)(25,20)}
\darkred{\qbezier(55,20)(56.5,24)(60,24)}
\darkred{\qbezier(25,20)(23.5,24)(20,24)}
\darkred{\qbezier(63,20)(63,24)(60,24)}
\darkred{\qbezier(17,20)(17,24)(20,24)}
\darkred{\qbezier(63,20)(63,17)(60,17)}
\darkred{\qbezier(17,20)(17,17)(20,17)}
\darkred{\qbezier(13.6,20)(13.6,14)(18,14)}
\darkred{\qbezier(20,17)(22,17)(22,20)}
\darkred{\qbezier(60,17)(58,17)(58,20)}
\darkred{\qbezier(18,14)(23,14)(23,20)}
\darkred{\qbezier(23,20)(23,22.5)(20,22.5)}
\darkred{\qbezier(20,21.5)(22,21.5)(22,20)}
\darkred{\qbezier(60,21.5)(58,21.5)(58,20)}
\darkred{\qbezier(20,21.5)(18.7,21.5)(18.7,20)}
\darkred{\qbezier(60,21.5)(61.3,21.5)(61.3,20)}
\darkred{\qbezier(20,22.5)(19,22.5)(18.6,22)}

  \darkred{\put(40,20){\circle{14}}}
  \darkred{\put(40,20){\circle{18}}}

\end{picture}
\end{center}
\caption{A lamination in a twice-punctured annulus with a total of 7
  marked points.}
\label{fig:lamin}
\end{figure}

\begin{figure}[htbp]
\begin{center}
\setlength{\unitlength}{4pt}
\begin{picture}(40,20)(0,0)
\thinlines

  \put(20,20){\circle*{1}}
  \put(15,10){\circle*{1}}

  \put(20,10){\circle{20}}

\thicklines
\lightblue{\qbezier(15,18.4)(20,15)(25,18.4)}
\lightblue{\qbezier(15,1.6)(20,5)(25,1.6)}

\lightblue{\put(15,10){\circle{6}}}
\lightblue{\put(25,10){\circle{6}}}

\end{picture}
\hspace{-0.6in}
\begin{picture}(40,20)(0,0)
\thinlines

  \put(20,20){\circle*{1}}
  \put(20,10){\circle*{1}}

  \put(20,10){\circle{20}}

\thicklines

\lightblue{\qbezier(17,4)(17,6)(20,6)}
\lightblue{\qbezier(17,4)(17,2)(20,2)}

\lightblue{\qbezier(24,10)(24,6)(20,6)}
\lightblue{\qbezier(26,8)(26,2)(20,2)}

\lightblue{\qbezier(24,10)(24,13.5)(20,13.5)}
\lightblue{\qbezier(26,8)(26,14.5)(20,14.5)}

\lightblue{\qbezier(17,10)(17,13.5)(20,13.5)}
\lightblue{\qbezier(16,10)(16,14.5)(20,14.5)}

\lightblue{\qbezier(17,10)(17,8)(20,8)}
\lightblue{\qbezier(16,10)(16,7)(20,7)}

\lightblue{\qbezier(20,8)(22,8)(22,10)}
\lightblue{\qbezier(20,7)(23,7)(23,10)}

\lightblue{\qbezier(23,10)(23,12.5)(20,12.5)}
\lightblue{\qbezier(20,11.8)(22,11.8)(22,10)}

\lightblue{\qbezier(20,11.8)(19,11.8)(18.8,11)}
\lightblue{\qbezier(20,12.5)(18.5,12.5)(18,11)}

\end{picture}
\end{center}
\caption{Curves that are not allowed in a lamination.}
\label{fig:not-lamin}
\end{figure}
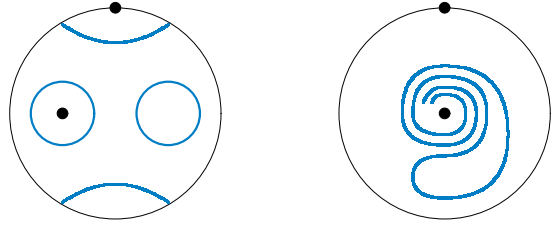

\begin{definition}
An \emph{integral unbounded measured lamination}---in this paper,
frequently just a \emph{lamination}---on a marked
surface~$\SM$ is a finite collection of non-self\/intersecting and
pairwise non\hyp intersecting curves in~$\Surf$, modulo isotopy relative
to~$\Mark$, subject to the restrictions specified below.
Each curve must be one of the following:
\begin{itemize}
\item
a closed curve (an embedded circle);
\item
a curve connecting two unmarked points on the boundary of~$\Surf$;
\item
a curve starting at an unmarked point on the boundary and, at its
other end, spiraling into a puncture (either clockwise or
counterclockwise); or
\item
a curve both of whose ends spiral into punctures (not necessarily
distinct).
\end{itemize}
See Figure~\ref{fig:lamin}.
Also, the following types of curves are not allowed 
(see Figure~\ref{fig:not-lamin}):
\begin{itemize}
\item
a curve that bounds an unpunctured or once-punctured disk;
\item
a curve with two endpoints on the boundary of~$\Surf$ which is
isotopic to a piece of boundary containing no marked
points, or a single marked point; and
\item a curve with two ends spiraling into the same puncture in the
  same direction
  without enclosing anything else.
\end{itemize}
\end{definition}

Note that a curve with two ends spiraling into the same puncture,
with the two ends spiraling in opposite directions, is excluded since
it is self\/intersecting.

Laminations on a marked surface $\SM$
can be coordinatized using W. Thurston's shear
coordinates.

\pagebreak[3]

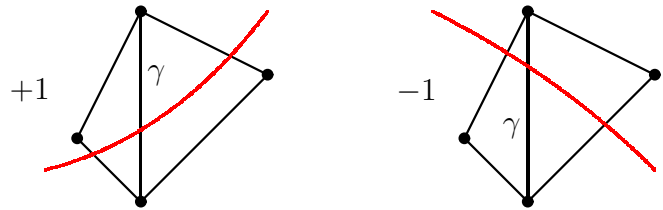
\begin{figure}[htbp]
\begin{center}
\setlength{\unitlength}{1.2pt}
\begin{picture}(60,60)(0,0)
\thicklines
  \put(0,20){\line(1,2){20}}
  \put(0,20){\line(1,-1){20}}
  \put(20,0){\line(0,1){60}}
  \put(20,0){\line(1,1){40}}
  \put(20,60){\line(2,-1){40}}

  \put(20,0){\circle*{3}}
  \put(20,60){\circle*{3}}
  \put(0,20){\circle*{3}}
  \put(60,40){\circle*{3}}

\put(25,40){\makebox(0,0){$\gamma$}}
\put(-15,35){\makebox(0,0){$+1$}}

\darkred{\qbezier(-10,10)(30,20)(60,60)}

\end{picture}
\qquad\qquad\qquad
\begin{picture}(60,60)(0,0)
\thicklines
  \put(0,20){\line(1,2){20}}
  \put(0,20){\line(1,-1){20}}
  \put(20,0){\line(0,1){60}}
  \put(20,0){\line(1,1){40}}
  \put(20,60){\line(2,-1){40}}

  \put(20,0){\circle*{3}}
  \put(20,60){\circle*{3}}
  \put(0,20){\circle*{3}}
  \put(60,40){\circle*{3}}

\put(-15,35){\makebox(0,0){$-1$}}
\put(15,23){\makebox(0,0){$\gamma$}}

\darkred{\qbezier(-10,60)(30,40)(60,10)}

\end{picture}
\end{center}
\caption{Defining the shear coordinate $b_\gamma(T,L)$.
The curve on the left contributes~$+1$, the one on the right contributes~$-1$.}
\label{fig:shear-def}
\end{figure}

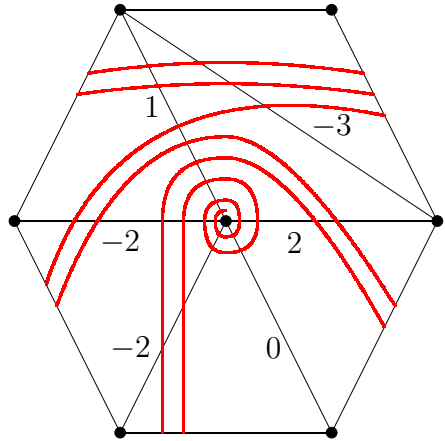
\begin{figure}[htbp]
\begin{center}
\setlength{\unitlength}{4pt}
\begin{picture}(40,40)(0,0)
\thinlines
  \put(10,0){\line(-1,2){10}}
  \put(40,20){\line(-1,2){10}}
  \put(10,0){\line(1,0){20}}
  \put(10,40){\line(1,0){20}}
  \put(0,20){\line(1,2){10}}
  \put(30,0){\line(1,2){10}}
  \put(10,0){\line(1,2){10}}
  \put(0,20){\line(1,0){40}}
  \put(10,40){\line(1,-2){10}}
  \put(20,20){\line(1,-2){10}}
  \put(10,40){\line(3,-2){30}}

  \put(0,20){\circle*{1}}
  \put(10,0){\circle*{1}}
  \put(10,40){\circle*{1}}
  \put(30,0){\circle*{1}}
  \put(30,40){\circle*{1}}
  \put(40,20){\circle*{1}}
  \put(20,20){\circle*{1}}

\thicklines

\darkred{\qbezier(16,0)(16,10)(16,20)}
\darkred{\qbezier(16,20)(16,24)(20,24)}
\darkred{\qbezier(23,20)(23,24)(20,24)}
\darkred{\qbezier(23,20)(23,17)(20,17)}
\darkred{\qbezier(20,17)(18,17)(18,20)}
\darkred{\qbezier(20,22)(18,22)(18,20)}
\darkred{\qbezier(20,22)(21.3,22)(21.3,20)}
\darkred{\qbezier(20,18.5)(21.3,18.5)(21.3,20)}
\darkred{\qbezier(20,18.5)(19,18.5)(19,20)}
\darkred{\qbezier(20,21)(19,21)(19,20)}

\darkred{\qbezier(14,0)(14,10)(14,20)}
\darkred{\qbezier(14,20)(14,26)(20,26)}

\darkred{\qbezier(35,10)(26,26)(20,26)}

\darkred{\qbezier(36,12)(26,28)(20,28)}
\darkred{\qbezier(4,12)(10,28)(20,28)}
\darkred{\qbezier(3,14)(10,35)(35,30)}

\darkred{\qbezier(6,32)(20,34)(34,32)}
\darkred{\qbezier(7,34)(20,36)(33,34)}

\put(24.5,8){\makebox(0,0){$0$}}
\put(11,8){\makebox(0,0){$-2$}}
\put(10,18){\makebox(0,0){$-2$}}
\put(26.5,18){\makebox(0,0){$2$}}
\put(30,29){\makebox(0,0){$-3$}}
\put(13,30.8){\makebox(0,0){$1$}}

\end{picture}
\quad

\end{center}
\caption{Shear coordinates of a lamination $L$ with respect to an
  ordinary triangulation~$T$.
}
\label{fig:shear-coord}
\end{figure}

\begin{definition}[\emph{Shear coordinates}]
\label{def:shear}
Let $L$ be an integral unbounded measured lamination.
Let $T$ be a triangulation
without self-folded triangles.
For each arc $\gamma$ in~$T$, the corresponding \emph{shear
  coordinate} of~$L$ with respect to the triangulation~$T$,
denoted by~$b_\gamma(T,L)$, is
defined as a sum of contributions from all intersections of curves in~$L$
with the arc~$\gamma$.
Specifically, such an intersection contributes $+1$ (resp.,~$-1$)
to $b_\gamma(T,L)$ if the corresponding segment of a curve in~$L$ cuts
through the quadrilateral surrounding $\gamma$ cutting through edges
in the shape of an~`S' (resp., in the shape of a~`Z'), as shown
in Figure~\ref{fig:shear-def} on the left (resp., on the right).
Note that at most one of these two types of intersection can occur.
Note also that even though a spiraling curve can intersect an arc
infinitely many times, the number of intersections that contribute to
the computation of $b_\gamma(T,L)$ is always finite.

An example is shown in Figure~\ref{fig:shear-coord}.
\end{definition}

An alternative (more conceptual) definition of shear coordinates can
be given using the notion of tropical lambda lengths, 
cf.~\eqref{eq:b-quadrilateral-sum}.

\begin{theorem}[W.~Thurston]
\label{th:thurston-1}
For a fixed triangulation~$T$ without self-folded triangles, the map
\[
L\mapsto (b_\gamma(T,L))_{\gamma\in T}
\]
is a bijection between integral unbounded measured
laminations and~$\ZZ^n$.
\end{theorem}

\begin{figure}[htbp]
\begin{center}
\begin{tabular}{ccc}
\setlength{\unitlength}{4pt}
\begin{picture}(25,15)(10,14)
\put(10,27){$L_1$}
\thinlines
  \put(10,20){\line(1,0){20}}
  \multiput(10,20)(10,0){3}{\circle*{1}}
\qbezier(10,20)(20,35)(30,20)
\qbezier(10,20)(20,5)(30,20)

\thicklines

\darkred{\qbezier(16.5,13.6)(16.5,18)(16.5,20)}
\darkred{\qbezier(16.5,20)(16.5,24)(20,24)}
\darkred{\qbezier(23,20)(23,24)(20,24)}
\darkred{\qbezier(23,20)(23,17)(20,17)}
\darkred{\qbezier(20,17)(18,17)(18,20)}
\darkred{\qbezier(20,22)(18,22)(18,20)}
\darkred{\qbezier(20,22)(21.3,22)(21.3,20)}
\darkred{\qbezier(20,18.5)(21.3,18.5)(21.3,20)}
\darkred{\qbezier(20,18.5)(19,18.5)(19,20)}
\darkred{\qbezier(20,21)(19,21)(19,20)}

\put(26,18.5){\makebox(0,0){$0$}}
\put(14,18.5){\makebox(0,0){$-1$}}

\end{picture}
&
\begin{picture}(25,15)(10,14)
\put(10,27){$L_2$}
\thinlines
  \put(10,20){\line(1,0){20}}
  \multiput(10,20)(10,0){3}{\circle*{1}}
\qbezier(10,20)(20,35)(30,20)
\qbezier(10,20)(20,5)(30,20)

\thicklines

\darkred{\qbezier(24,20)(24,24.5)(20,24.5)}
\darkred{\qbezier(16,20)(16,24.5)(20,24.5)}
\darkred{\qbezier(16,13.8)(16,18)(16,20)}
\darkred{\qbezier(24,13.8)(24,18)(24,20)}

\put(26,18.5){\makebox(0,0){$1$}}
\put(14,18.5){\makebox(0,0){$-1$}}

\end{picture}
&
\begin{picture}(25,15)(10,14)
\put(10,27){$L_3$}
\thinlines
  \put(10,20){\line(1,0){20}}
  \multiput(10,20)(10,0){3}{\circle*{1}}
\qbezier(10,20)(20,35)(30,20)
\qbezier(10,20)(20,5)(30,20)

\thicklines

\darkred{\qbezier(23.5,13.6)(23.5,18)(23.5,20)}
\darkred{\qbezier(23.5,20)(23.5,24)(20,24)}
\darkred{\qbezier(17,20)(17,24)(20,24)}
\darkred{\qbezier(17,20)(17,17)(20,17)}
\darkred{\qbezier(20,17)(22,17)(22,20)}
\darkred{\qbezier(20,22)(22,22)(22,20)}
\darkred{\qbezier(20,22)(18.7,22)(18.7,20)}
\darkred{\qbezier(20,18.5)(18.7,18.5)(18.7,20)}
\darkred{\qbezier(20,18.5)(21,18.5)(21,20)}
\darkred{\qbezier(20,21)(21,21)(21,20)}

\put(26,18.5){\makebox(0,0){$1$}}
\put(14,18.5){\makebox(0,0){$0$}}

\end{picture}
\\[.3in] 

\begin{picture}(25,15)(10,14)
\put(10,27){$L_6$}
\thinlines
  \put(10,20){\line(1,0){20}}
  \multiput(10,20)(10,0){3}{\circle*{1}}
\qbezier(10,20)(20,35)(30,20)
\qbezier(10,20)(20,5)(30,20)

\thicklines

\darkred{\qbezier(23.5,26.4)(23.5,22)(23.5,20)}
\darkred{\qbezier(23.5,20)(23.5,16)(20,16)}
\darkred{\qbezier(17,20)(17,16)(20,16)}
\darkred{\qbezier(17,20)(17,23)(20,23)}
\darkred{\qbezier(20,23)(22,23)(22,20)}
\darkred{\qbezier(20,18)(22,18)(22,20)}
\darkred{\qbezier(20,18)(18.7,18)(18.7,20)}
\darkred{\qbezier(20,21.5)(18.7,21.5)(18.7,20)}
\darkred{\qbezier(20,21.5)(21,21.5)(21,20)}
\darkred{\qbezier(20,19)(21,19)(21,20)}

\put(26,18.5){\makebox(0,0){$-1$}}
\put(14,18.5){\makebox(0,0){$0$}}

\end{picture}
&
\begin{picture}(25,15)(10,14)
\put(10,27){$L_5$}
\thinlines
  \put(10,20){\line(1,0){20}}
  \multiput(10,20)(10,0){3}{\circle*{1}}
\qbezier(10,20)(20,35)(30,20)
\qbezier(10,20)(20,5)(30,20)

\thicklines

\darkred{\qbezier(24,20)(24,15.5)(20,15.5)}
\darkred{\qbezier(16,20)(16,15.5)(20,15.5)}
\darkred{\qbezier(16,26.2)(16,22)(16,20)}
\darkred{\qbezier(24,26.2)(24,22)(24,20)}

\put(26,18.5){\makebox(0,0){$-1$}}
\put(14,18.5){\makebox(0,0){$1$}}

\end{picture}
&
\begin{picture}(25,15)(10,14)
\put(10,27){$L_4$}
\thinlines
  \put(10,20){\line(1,0){20}}
  \multiput(10,20)(10,0){3}{\circle*{1}}
\qbezier(10,20)(20,35)(30,20)
\qbezier(10,20)(20,5)(30,20)

\thicklines

\darkred{\qbezier(16.5,26.4)(16.5,22)(16.5,20)}
\darkred{\qbezier(16.5,20)(16.5,16)(20,16)}
\darkred{\qbezier(23,20)(23,16)(20,16)}
\darkred{\qbezier(23,20)(23,23)(20,23)}
\darkred{\qbezier(20,23)(18,23)(18,20)}
\darkred{\qbezier(20,18)(18,18)(18,20)}
\darkred{\qbezier(20,18)(21.3,18)(21.3,20)}
\darkred{\qbezier(20,21.5)(21.3,21.5)(21.3,20)}
\darkred{\qbezier(20,21.5)(19,21.5)(19,20)}
\darkred{\qbezier(20,19)(19,19)(19,20)}

\put(26,18.5){\makebox(0,0){$0$}}
\put(14,18.5){\makebox(0,0){$1$}}

\end{picture}
\end{tabular}
\setlength{\unitlength}{4pt}
\begin{picture}(20,20)(1,10)

\put(10,10){\circle*{1}}

\thinlines
\put(10,10){\vector(1,0){10}}
\put(10,10){\vector(0,1){10}}
\put(10,10){\vector(-1,1){10}}
\put(10,10){\vector(-1,0){10}}
\put(10,10){\vector(0,-1){10}}
\put(10,10){\vector(1,-1){10}}

\put(2,8){\makebox(0,0){$\mathbf{y}_1$}}
\put(2,16){\makebox(0,0){$\mathbf{y}_2$}}
\put(12,19){\makebox(0,0){$\mathbf{y}_3$}}
\put(19,12){\makebox(0,0){$\mathbf{y}_4$}}
\put(17,1){\makebox(0,0){$\mathbf{y}_5$}}
\put(8,1){\makebox(0,0){$\mathbf{y}_6$}}

\end{picture}
\end{center}
\caption{Six ``elementary'' laminations of a once-punctured digon,
and the corresponding vectors of shear coordinates.
}
\label{fig:lamin-digon}
\end{figure}
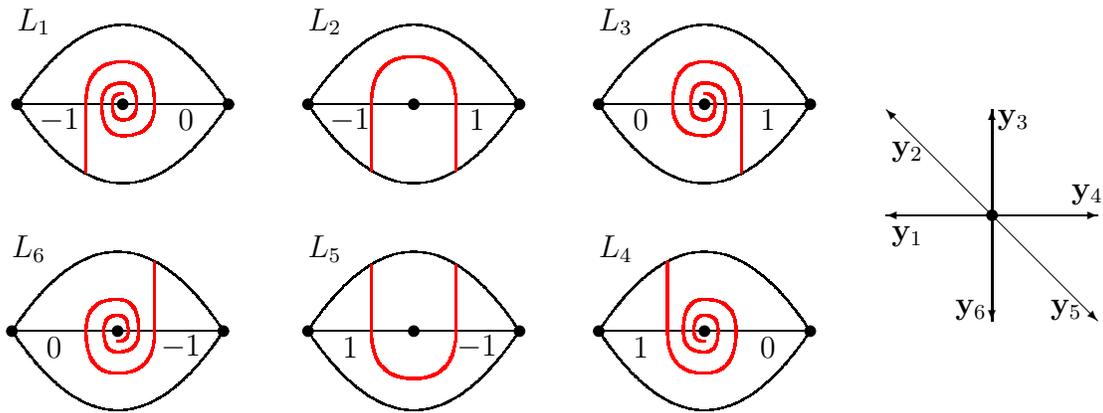

\begin{example}
\label{example:elem-lamin-punctured-digon}
Figure~\ref{fig:lamin-digon} shows the six ``elementary'' laminations
$L_1,\dots,L_6$ of
a once-punctured digon (each lamination $L_i$ consisting of a single curve),
and their shear coordinates with
respect to a particular triangulation~$T$.
It is easy to see that any integral lamination of the digon consists
of several (possibly none) curves homotopic to some~$L_i$, together
with several (possibly none) curves homotopic to $L_{(i+1)\bmod
  6}$.
It~is also easy to see that, in agreement with
Theorem~\ref{th:thurston-1},  each vector in $\ZZ^2$ can be uniquely
written as a non-negative integer linear combination
\[
p \mathbf{y}_i + q \mathbf{y}_{(i+1)\bmod 6} \quad (p,q\in\ZZ_{\ge
  0}),
\]
where the six vectors
\[
 \mathbf{y}_1=[-1,0],\ \  \mathbf{y}_2=[-1,1], \ \  \mathbf{y}_3=[0,1], \ \
 \mathbf{y}_4=[1,0], \ \  \mathbf{y}_5=[1,-1], \ \  \mathbf{y}_6=[0,-1]
\]
represent the shear coordinates of $L_1,\dots,L_6$, respectively.
\end{example}

\begin{definition}[\emph{Multi-laminations and associated extended
    exchange matrices}]
\label{def:B-tilde-multi-lam}
A~\emph{multi-lamination} is simply a finite family of
laminations.
Let us fix a multi-lamination
\[
\LL=(L_{n+1},\dots, L_{m})
\]
of size $m-n$;
this choice of indexing will be convenient in the sequel.
For a triangulation~$T$ of $\SM$ without self-folded triangles, define an
$m\times n$ matrix
\[
\tilde B=\tilde B(T,\LL)=(b_{ij})
\]
(cf.\ Definition~\ref{def:B-tilde}) as follows.
The top $n\times n$ part of $\tilde B$ is the~signed
adjacency matrix $B(T)=(b_{ij})_{1\le i,j\le n}$ (cf.\
Definition~\ref{def:signed-adjacency-matrix}),
whereas the bottom $m-n$ rows are formed by the shear coordinates of the
laminations $L_i$ with respect to the triangulation~$T$:
\begin{equation}
\label{eq:bijTL}
b_{ij}=b_j(T,L_i)\quad \text{if} \quad  n<i\le m.
\end{equation}
\end{definition}

\begin{figure}[htbp]
\begin{center}
\setlength{\unitlength}{4pt}
\begin{picture}(40,41)(25,0)
\thinlines
  \put(10,0){\line(-1,2){10}}
  \put(40,20){\line(-1,2){10}}
  \put(10,0){\line(1,0){20}}
  \put(10,40){\line(1,0){20}}
  \put(0,20){\line(1,2){10}}
  \put(30,0){\line(1,2){10}}
  \put(10,0){\line(1,2){10}}
  \put(0,20){\line(1,0){40}}
  \put(10,40){\line(1,-2){10}}
  \put(20,20){\line(1,-2){10}}
  \put(10,40){\line(3,-2){30}}

  \put(0,20){\circle*{1}}
  \put(10,0){\circle*{1}}
  \put(10,40){\circle*{1}}
  \put(30,0){\circle*{1}}
  \put(30,40){\circle*{1}}
  \put(40,20){\circle*{1}}
  \put(20,20){\circle*{1}}

\thicklines

\darkred{\qbezier(16,0)(16,10)(16,20)}
\darkred{\qbezier(16,20)(16,24)(20,24)}
\darkred{\qbezier(23,20)(23,24)(20,24)}
\darkred{\qbezier(23,20)(23,17)(20,17)}
\darkred{\qbezier(20,17)(18,17)(18,20)}
\darkred{\qbezier(20,22)(18,22)(18,20)}
\darkred{\qbezier(20,22)(21.3,22)(21.3,20)}
\darkred{\qbezier(20,18.5)(21.3,18.5)(21.3,20)}
\darkred{\qbezier(20,18.5)(19,18.5)(19,20)}
\darkred{\qbezier(20,21)(19,21)(19,20)}

\darkred{\qbezier(14,0)(14,10)(14,20)}
\darkred{\qbezier(14,20)(14,26)(20,26)}

\darkred{\qbezier(35,10)(26,26)(20,26)}

\darkred{\qbezier(36,12)(26,28)(20,28)}
\darkred{\qbezier(4,12)(10,28)(20,28)}
\darkred{\qbezier(3,14)(10,35)(35,30)}

\darkred{\qbezier(6,32)(20,34)(34,32)}
\darkred{\qbezier(7,34)(20,36)(33,34)}


\put(45,19){
$\tilde B=
\left[\,\,
\begin{matrix}
0 & -1& 0 & 0 & 1 &-1\\[.05in]
1 & 0 &-1 & 0 & 0 & 0\\[.05in]
0 & 1 & 0 &-1& 0 & 0\\[.05in]
0 & 0 & 1 & 0&-1 & 0\\[.05in]
-1& 0 & 0 & 1 & 0 & 1\\[.05in]
1 & 0 & 0 & 0 & -1& 0\\[.05in]
\hline\\[-.15in]
2 & 0 &-2 &-2 & 1 & -3
\end{matrix}
\,\,\right]
$
}
\end{picture}
\quad

\end{center}
\caption{The matrix $\tilde B=\tilde B(T,\LL)$ for the
example in Figure~\ref{fig:shear-coord}.
The multi-lamination $\LL=(L)$ consists of a single
  lamination~$L$.
We use the labeling of arcs in~$T$ shown in Figure~\ref{fig:B(T)-hexagon}.
}
\label{fig:shear-coord-matrix}
\end{figure}
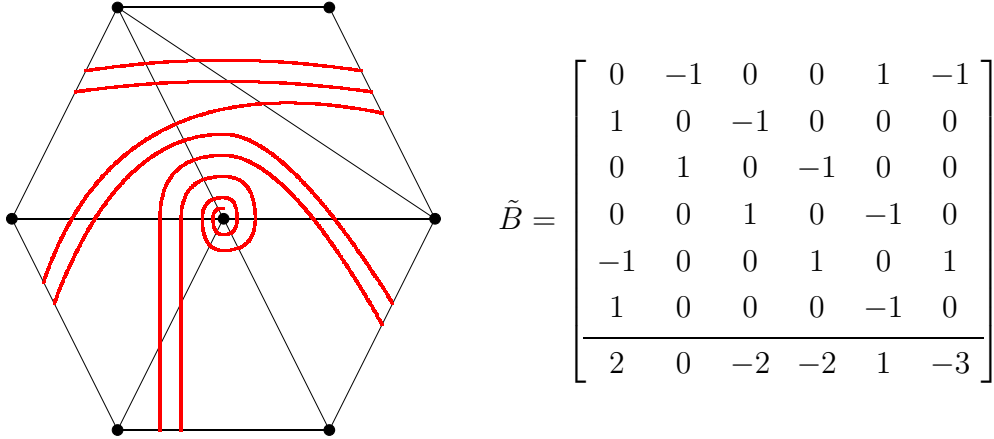

The key observation is that, under ordinary quadrilateral flips,
the matrices $\tilde B(T)$ transform according to the mutation rules.

\begin{theorem}[W.~Thurston--V.~Fock--A.~Goncharov]
\label{th:thurston-2}
Let us fix a multi-lamination~$\LL$.
If~triangulations $T$ and $T'$ without self-folded triangles are
related by a flip of an arc~$k$,
then the corresponding matrices $\tilde B(T,\LL)$ and $\tilde
B(T',\LL)$ are related by a mutation in direction~$k$.
\end{theorem}

Although the reader will not find the exact statement of
Theorem~\ref{th:thurston-2} in the work of the above authors,
it can be seen to be a version of the results contained therein.
More specifically, applying the definition of a matrix mutation to the
case under consideration results in identities for the shear
coordinates (with respect to $T$ and~$T'$) that appear, for example,
at the end of \cite[Section~3.1]{fg-dual-teich}.
It is elementary to check these identities directly; in
Chapter~\ref{sec:tropical-lambda} we will give a more
conceptual proof.

\chapter{Shear coordinates with respect to tagged triangulations}
\label{sec:shear-tagged}

In this chapter, we define shear coordinates for triangulations with
self-folded triangles and, more generally, for tagged
triangulations.
We then obtain the appropriate analogues of
Theorems~\ref{th:thurston-1} and~\ref{th:thurston-2}.

\begin{definition}[\emph{Shear coordinates with respect to a tagged
    triangulation}]
\label{def:shear-tagged}
We extend Definition~\ref{def:shear} by defining the \emph{shear
  coordinates}~$b_\gamma(T,L)$ of an integral unbounded measured
lamination~$L$ with respect to an arbitrary tagged triangulation~$T$.
(Here $\gamma$ runs over the tagged arcs in~$T$.)
These coordinates are uniquely defined by the following rules:
\begin{enumerate}
\item[(i)]
Suppose that tagged triangulations $T_1$ and $T_2$ coincide 
except that at a particular puncture~$p$, the tags of the arcs
in~$T_1$ are all different from the tags of their counterparts in~$T_2$.
Suppose that laminations $L_1$ and $L_2$ coincide 
except that each curve in $L_1$ that spirals into~$p$ has been
replaced in $L_2$ by a curve that spirals in the opposite direction.
Then $b_{\gamma_1}(T_1,L_1)=b_{\gamma_2}(T_2,L_2)$ for each tagged arc
$\gamma_1\in T_1$ and its counterpart $\gamma_2\in T_2$.
\item[(ii)]
By performing tag-changing transformations $L_1\leadsto L_2$ with
$L_1$ and $L_2$ as above, we can convert any tagged triangulation into
a tagged triangulation~$T$ that does not contain any notches except
possibly inside once-punctured digons.
Let $T^\circ$ denote the ideal triangulation that is represented
by such~$T$; that is, $T=\tau(T^\circ)$ in the notation of
Definitions~\ref{def:arcs-as-tagged-arcs} and~\ref{def:tagged-triang}.
Let $\gamma^\circ$ be an arc in~$T^\circ$ that is not contained inside
a self-folded triangle, and let $\gamma=\tau(\gamma^\circ)$.
Then, for a lamination~$L$, we define $b_\gamma(T,L)$ by applying the
rule in Definition~\ref{def:shear} to the ordinary arc~$\gamma^\circ$
viewed inside the triangulation~$T^\circ$.
\end{enumerate}
Note that if $\gamma^\circ$ is contained inside a self-folded triangle
in~$T^\circ$ enveloping a puncture~$p$,
then we can first apply the tag-changing transformation~(i) to~$T$ at~$p$, and
then use the rule~(ii) to determine the shear coordinate in question.
\end{definition}

\begin{example}
\label{example:4-triang-of-punct-digon}
Let $T$ be the tagged triangulation of a punctured digon shown in
Figure~\ref{fig:tagged-triang-digon} (cf.\ also
Figure~\ref{fig:lamin-digon}),
and let $T_1$, $T_2$, and $T_{12}$ be the tagged triangulations
obtained from~$T$ as follows:
\begin{itemize}
\item
$T_1$ is obtained from $T$ by the tagged flip
replacing $\gamma_{ap}$
by the tagged arc~$\gamma_{b\bar p}$;
\item
$T_2$ is obtained from $T$ by the tagged flip
replacing $\gamma_{bp}$
by the tagged arc~$\gamma_{a\bar p}$;
\item
$T_{12}$ is obtained from $T$ by performing both of these (commuting)
tagged flips.
\end{itemize}
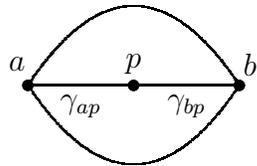
\begin{figure}[htbp]
\begin{center}
\setlength{\unitlength}{4pt}
\begin{picture}(25,12)(10,14)
\put(9,22){\makebox(0,0){$a$}}
\put(20,22){\makebox(0,0){$p$}}
\put(31,22){\makebox(0,0){$b$}}
\put(15,18){\makebox(0,0){$\gamma_{ap}$}}
\put(25,18){\makebox(0,0){$\gamma_{bp}$}}
\thinlines
  \put(10,20){\line(1,0){20}}
  \multiput(10,20)(10,0){3}{\circle*{1}}
\qbezier(10,20)(20,35)(30,20)
\qbezier(10,20)(20,5)(30,20)
\end{picture}
\end{center}
\caption{A triangulation $T$ of a once-punctured digon.
}
\label{fig:tagged-triang-digon}
\end{figure}
Let $L$ be a lamination in this once-punctured digon,
with shear coordinates $b_1(T,L)$ and  $b_2(T,L)$ corresponding to the
arcs $\gamma_{ap}$ and $\gamma_{bp}$ of~$T$, respectively.
We similarly define $b_1(T_s,L)$ and $b_2(T_s,L)$ for
each subscript $s\in\{1,2,12\}$.
For example, $b_1(T_{12},L)$ and $b_2(T_{12},L)$
refer to the shear coordinates associated with $\gamma_{b\bar p}$ and
$\gamma_{a\bar p}$, respectively,
since these tagged arcs replace $\gamma_{ap}$ and $\gamma_{bp}$,
respectively. 
Then the rules (i)--(ii) of Definition~\ref{def:shear-tagged} yield:
\[
b_j(T_s,L)=\begin{cases}
-b_j(T,L) & \text{if $j$ appears in~$s$;}\\
b_j(T,L) & \text{if $j$ does not appear in~$s$.}
\end{cases}
\]
For example, the shear coordinates of the laminations~$L_1$ and $L_3$
(in the notation of Figure~\ref{fig:lamin-digon}) with respect to the tagged
triangulations $T,T_1,T_2,T_{12}$ are:
\begin{align*}
b_1(T,L_1)&=-1 & b_2(T,L_1)&=0 &
  b_1(T,L_3)&=0  & b_2(T,L_3)&=1\\
b_1(T_1,L_1)&=1  & b_2(T_1,L_1)&=0  &
  b_1(T_1,L_3)&=0  & b_2(T_1,L_3)&=1\\
b_1(T_2,L_1)&=-1  & b_2(T_2,L_1)&=0  &
  b_1(T_2,L_3)&=0  & b_2(T_2,L_3)&=-1\\
b_1(T_{12},L_1)&=1  & b_2(T_{12},L_1)&=0  &
  b_1(T_{12},L_3)&=0  & b_2(T_{12},L_3)&=-1
\end{align*}


\end{example}

We now extend Definition~\ref{def:B-tilde-multi-lam} to the
tagged case.

\begin{definition}
\label{def:tilde-BL-tagged}
The \emph{extended exchange matrix} $\tilde B(T,\LL)$ of a
multi-lamination $\LL$ with respect to a tagged
triangulation~$T$ is defined in exactly the same way as in
Definition~\ref{def:B-tilde-multi-lam}, this time using the shear
coordinates from Definition~\ref{def:shear-tagged}.
\end{definition}

\begin{example}
\label{example:4-triang-of-punct-digon-continued}

Continuing with Example~\ref{example:4-triang-of-punct-digon},
let $\LL=(L_1)$. 
Then, e.g.,
\[
\tilde B(T,\LL)=\begin{bmatrix}
0 & 0\\
0 & 0\\
-1& 0
\end{bmatrix},\quad
\tilde B(T_1,\LL)=\begin{bmatrix}
0 & 0\\
0 & 0\\
1 & 0
\end{bmatrix}.
\]
\end{example}

We are now prepared to provide the promised generalizations of
Theorems~\ref{th:thurston-1} and~\ref{th:thurston-2}.

\begin{theorem}
\label{th::thurston-2-tagged}
Fix a multi-lamination~$\LL$.
If tagged triangulations $T$ and $T'$ are
related by a flip of a tagged arc~$k$,
then the corresponding matrices $\tilde B(T,\LL)$ and $\tilde
B(T',\LL)$ are related by a mutation in direction~$k$.
\end{theorem}

\begin{proof}
The proof is a straightforward albeit tedious case-by-case verification based on
Theorem~\ref{th:thurston-2}, Definition~\ref{def:shear-tagged}, and
Remark~\ref{rem:two-types-of-flips}.  
For a flip inside a punctured digon, the analysis involves the laminations $L_1,\dots,L_6$
from Example~\ref{example:elem-lamin-punctured-digon}.  (The more delicate
part concerns the transformation of the shear coordinates of the arcs on the boundary of the
digon.)
\end{proof}

In Chapter~\ref{sec:tropical-lambda} we will give an alternative 
(more conceptual and more detailed) 
proof of Theorem~\ref{th::thurston-2-tagged}.

\begin{theorem}
\label{th:thurston-param-tagged}
For a fixed tagged triangulation~$T$, the map
\[
L\mapsto (b_\gamma(T,L))_{\gamma\in T}
\]
is a bijection between integral unbounded measured
laminations and~$\ZZ^n$.
\end{theorem}

\begin{proof}
This statement follows from Theorem~\ref{th:thurston-1},
Theorem~\ref{th::thurston-2-tagged}, and
the invertibility of matrix mutations.
Proceed by induction
on the number of flips required to obtain $T$ from a triangulation
without self-folded triangles.
\end{proof}

\chapter{Tropical lambda lengths}
\label{sec:tropical-lambda}

A na\"ive definition of a tropical lambda length of an arc~$\gamma$ with
respect to a \hbox{(multi-)}lamination~$L$ is based on the notion of a
transverse
measure of~$\gamma$ with respect to~$L$, i.e., the corresponding
intersection number.
The latter notion is however ill defined, as a curve in $L$ spiraling into a
puncture~$p$ intersects the arcs incident to~$p$ infinitely many times.
(In the absence of punctures, this issue does not come up.)
We bypass this problem by passing to the opened surface where
$\gamma$ is replaced by a family of lifts~$\bargamma$,
as we did in Chapter~\ref{sec:lambda-opened} for lambda lengths.
This sets the stage for Chapter~\ref{sec:laminated-teich},
where the tropical lambda lengths of those lifts
are used as rescaling factors allowing us to construct the requisite
normalized patterns with arbitrary coefficients of geometric type.

\begin{remark}
One can alternatively define tropical lambda lengths via a
limiting procedure that degenerates a hyperbolic structure on $\SM$
into a discrete, or tropical, version thereof, in the spirit of
W.~Thurston's approach to compactifying Teichm\"uller spaces.
Further hints are provided in
Appendix~\ref{app:relative-lambda-lengths}.
In this paper, we do not systematically pursue this course of action,
as the limiting objects, non-integral measured laminations, are more
complicated than we need.
\end{remark}

\begin{definition}
  An (integral) \emph{lifted measured lamination} $\barL$ on
  $\barSM=(\Surf_{\barM}, \Mark_{\barM})$ (cf.~\ref{eq:barSM})
  consists of a choice of orientation on each opened puncture~$C_p$
  together with a
  finite collection of non-intersecting curves on $\barSM$, modulo
  isotopy relative to~$\Mark_{\barM}$, 
with the following restrictions.  First, each
  component is
  \begin{itemize}
  \item a closed curve, or else 
  \item a curve connecting two points on the boundary of $\Surf_{\barM}$ away
    from $\Mark_{\barM}$.
  \end{itemize}
  Second, the following types of curves are not allowed:
  \begin{itemize}
  \item
    curves that bound an unpunctured disk or a disk with a single
    (opened) puncture; and
  \item
    curves with two endpoints on the boundary of~$\Surf_{\barM}$ which are
    isotopic to a piece of boundary containing no marked
    points, or a single marked point.
  \end{itemize}
  There is a natural projection map taking lifted measured
  laminations~$\barL$ on $\barSM$ to unbounded measured
  laminations~$L$ on $\SM$: take the endpoints of~$\barL$
  that end at an opened puncture~$C_p$ and make them spiral around the
  corresponding puncture~$p$ in the direction \emph{opposite} to the
  orientation chosen (in~$\barL$) on~$C_p$.  In this case $\barL$ is
  said to be a \emph{lift} of~$L$.

  Since all curves in a lamination have a consistent direction
  of spiraling, every lamination has at least one lift.

This notion straightforwardly generalizes to multi-laminations:
a \emph{lifted multi\hyp lamination} $\barLL$ consists of an (uncoordinated)
collection of lifted laminations.  In particular, each of these 
lifted laminations has its own orientation on each~$C_p$.  The notion
of projection likewise carries over.
\end{definition}

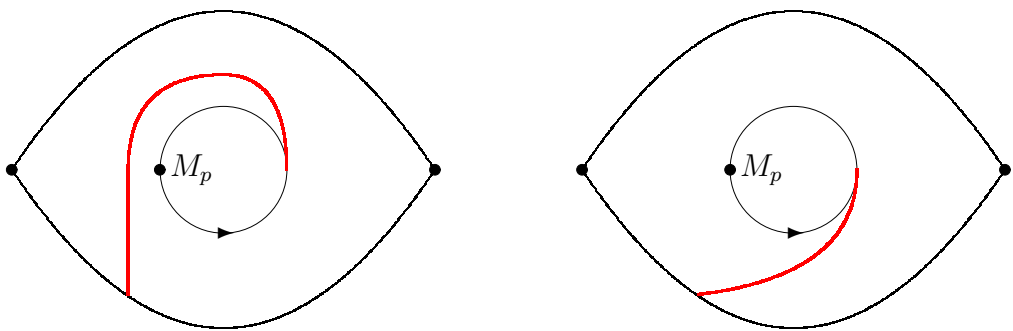
\begin{figure}[htbp]
\begin{center}

\setlength{\unitlength}{8pt}

\begin{picture}(25,15)(8,12.5)
\thinlines
\multiput(10,20)(20,0){2}{\circle*{0.5}}
\put(17,20){\circle*{0.5}}

\qbezier(10,20)(20,35)(30,20)
\qbezier(10,20)(20,5)(30,20)

\put(20,20){\circle{6}}

\thicklines

\darkred{\qbezier(15.5,14.1)(15.5,18)(15.5,20)}
\darkred{\qbezier(15.5,20)(15.5,24.5)(20,24.5)}
\darkred{\qbezier(23,20)(23,24.5)(20,24.5)}

 \put(18.5,20){\makebox(0,0){$M_p$}}

\put(20.5,17){\vector(1,0){0}}

\end{picture}
\quad
\begin{picture}(25,15)(8,12.5)
\thinlines
\multiput(10,20)(20,0){2}{\circle*{0.5}}
\put(17,20){\circle*{0.5}}

\qbezier(10,20)(20,35)(30,20)
\qbezier(10,20)(20,5)(30,20)

\put(20,20){\circle{6}}

\thicklines

\darkred{\qbezier(15.5,14.1)(23,15)(23,20)}

 \put(18.5,20){\makebox(0,0){$M_p$}}
\put(20.5,17){\vector(1,0){0}}

\end{picture}
\end{center}
\caption{Different lifts of the lamination~$L_1$ from
  Figure~\ref{fig:lamin-digon}.
}
\label{fig:lifts-lamin}
\end{figure}

\begin{remark}
  For a lamination~$L$ on $\SM$ with at least one curve spiraling
  into each puncture, the lifts~$\barL$ of~$L$ are parametrized by
  $\ZZ^{\text{number of punctures}}$.  These lifts differ by twisting
  $\barL$ around the opened circles~$C_p$, as illustrated in
  Figure~\ref{fig:lifts-lamin}.  This is analogous to
  Remark~\ref{rem:dTeich-lift}.  On the other hand, when $L$ has no
  spiraling into a given puncture~$p$, the lifts $\barL$ have no endpoints
  on~$C_p$---so locally, there are only two
  lifts, corresponding to the choices of orientation on~$C_p$.  To
  complete the analogy, we could forget the orientation here (in this
  case, it
  has no effect on the projection anyway), and admit
  curves in~$\barL$ that enclose a single marked point (either a
  closed curve enclosing~$C_p$, or a curve 
  cutting off a single marked point on the boundary); these extra curves are the
  analogues of the choice of a horocycle. 
  We do not pursue this further here  since for our
  main goal (Theorem~\ref{th:thurston-param-tagged}), we only need a
  single lift of any lamination, and since in order to get a full $\ZZ$'s worth
  of lifts and deal correctly with tagged arcs we would have to allow
  virtual (i.e., formally negative) curves enclosing a single marked point.
\end{remark}

\begin{figure}[htbp]
\begin{center}

\setlength{\unitlength}{8pt}

\begin{picture}(25,15)(8,12.5)
\thinlines
\qbezier(30,20)(17,10)(17,20)
\put(25,18){\makebox(0,0){$0$}}
\qbezier(30,20)(17,30)(17,20)
\put(25,22){\makebox(0,0){$1$}}

\multiput(10,20)(20,0){2}{\circle*{0.5}}
\put(17,20){\circle*{0.5}}

\qbezier(10,20)(20,35)(30,20)
\qbezier(10,20)(20,5)(30,20)

\put(20,20){\circle{6}}

\thicklines

\darkred{\qbezier(15.5,14.1)(15.5,18)(15.5,20)}
\darkred{\qbezier(15.5,20)(15.5,24.5)(20,24.5)}
\darkred{\qbezier(23,20)(23,24.5)(20,24.5)}

 \put(18.5,20){\makebox(0,0){$M_p$}}

\put(20.5,17){\vector(1,0){0}}

\end{picture}
\quad
\begin{picture}(25,15)(8,12.5)
\thinlines
\qbezier(30,20)(17,10)(17,20)
\put(25,18){\makebox(0,0){$-1$}}
\qbezier(30,20)(17,30)(17,20)
\put(25,22){\makebox(0,0){$0$}}

\multiput(10,20)(20,0){2}{\circle*{0.5}}
\put(17,20){\circle*{0.5}}

\qbezier(10,20)(20,35)(30,20)
\qbezier(10,20)(20,5)(30,20)

\put(20,20){\circle{6}}

\thicklines

\darkred{\qbezier(15.5,14.1)(23,15)(23,20)}

 \put(18.5,20){\makebox(0,0){$M_p$}}
\put(20.5,17){\vector(1,0){0}}

\end{picture}
\end{center}
\caption{Transverse measures of arcs with respect to a
  lifted lamination. Here $l_\barL(C_p)=1$, so $l_\barL(p)=-1$, and
  \eqref{eq:l(p)-twist-trop} becomes
$l_\barL(\psi_p\gamma) = l_\barL(\gamma)-1$.}
\label{fig:transverse-measures}
\end{figure}
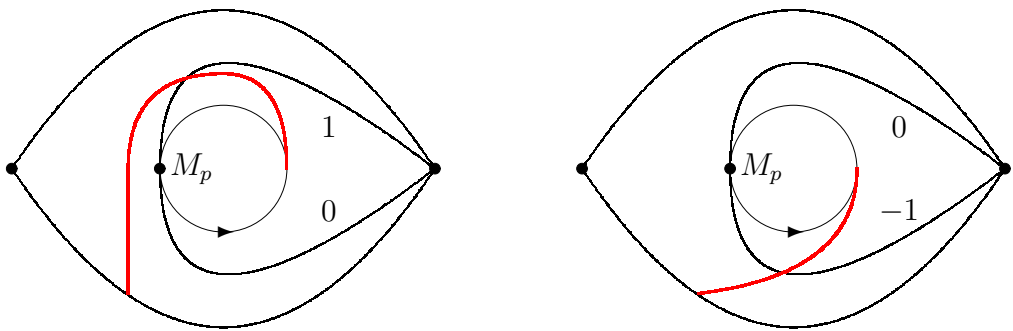

\begin{definition}[\emph{Transverse measures}]
\label{def:transverse-measure}
Let $\barL$ be a lifted lamination on the opened
surface $\barSM$. Let $\gamma$ be an (ordinary) arc in $\barAPSM$,
or a boundary segment in $\BSM$.
The \emph{transverse measure} of~$\gamma$ with respect to~$\barL$
is an integer denoted by $l_\barL(\gamma)$ and defined as follows.

If $\gamma$ does not have ends at~$M_p$
(for $p\in\barM$), then $l_\barL(\gamma)$ is simply the (non-negative)
minimal number of intersection points between the curves
in~$\barL$ and any arc homotopic
to~$\gamma$ (relative to endpoints).

Next suppose $\gamma$ has one or two ends that end at $M_p$ for
$p\in\barM$, and that $\gamma$ has no notches.
If $\gamma$ twists sufficiently far
around the corresponding opening(s) in the direction consistent with
the orientation of $\barL$ at each end, then, again, $l_\barL(\gamma)$
is equal to the minimal number
of intersections between $\gamma$ and~$\barL$.
(The notion of ``sufficiently far'' will depend on the
choice of a lift~$\barL$.)

We then extend the definition to all untagged arcs using the approach
utilized earlier to define $l(\gamma)$ in
Definition~\ref{def:lambda-opened}.  By analogy with
\eqref{eq:l-twist}, we require that
\begin{equation}
\label{eq:l(p)-twist-trop}
l_\barL(\psi_p\gamma)
= n_p(\gamma)l_\barL(p) + l_\barL(\gamma),
\end{equation}
where,
as before, we use the notation \eqref{eq:psi_p}--\eqref{eq:n_p}, and
\begin{equation}
\label{eq:l(p)-trop}
l_\barL(p) =
\begin{cases}
  -\text{$l_\barL(C_p)$}&p \in P, \text{if $C_p$ is oriented counterclockwise;}\\
  0& \text{if $p \not\in P$;}\\
  \text{$l_\barL(C_p)$}&p \in P, \text{if $C_p$ is oriented clockwise.}
\end{cases}
\end{equation}
As in the earlier case, property \eqref{eq:l(p)-twist-trop} is
consistent with the definition given above for the arcs that
twist sufficiently far.
Note that with this extended definition, the numbers $l_\barL(\gamma)$ may be
negative.
See Figure~\ref{fig:transverse-measures}.

Finally, for a tagged arc $\gamma\in\barATSM$ which is notched at~$p$ and
twists sufficiently far in the direction \emph{opposite} to the
orientation of $\barL$ on~$C_p$, define $l_\barL(\gamma)$ to be the
number of intersections of $\gamma$ with~$\barL$, plus
$\abs{l_\barL(p)}$.  (This extra term corresponds to the asymptotics
of $\nu(p)$ as described in Remark~\ref{rem:nu-asymptotics}, and will make
Lemma~\ref{lem:holed-monogon-trop} below come out uniformly.)
For notched arcs that do not twist sufficiently far, we extend the
definition using~\eqref{eq:l(p)-twist-trop}; as before, 
$l_\barL(\gamma)$ may be negative.  (Recall that $n_p(\gamma) < 0$ if
$\gamma$ is notched at~$p$.)
\end{definition}

We prefer to write these transverse measures multiplicatively, to make
the analogy
with geometric lambda lengths stronger and to match usual cluster
algebra notation.

\begin{definition}[\emph{Tropical semifield associated with a
    multi-lamination}]
\label{def:trop-lamin}
Let $\LL=(L_i)_{i\in I}$ be a multi-lamination in $\SM$; here $I$ is
a finite indexing set.
We introduce a formal variable~$q_i$ for each lamination~$L_i$, and
set (see~Definition~\ref{def:tropical})
\begin{equation}
\label{eq:PPLL}
\PP_\LL=\Trop(q_i:i\in I).
\end{equation}
\end{definition}

\begin{definition}[\emph{Tropical lambda lengths}]
\label{def:trop-lambda}
We continue in close analogy with Definition~\ref{def:lambda-opened}.
Let $\barLL=(\barL_i)_{i\in I}$ be a lift of a multi-lamination~$\LL$.
Instead of exponentiating the distances to get the
lambda lengths (cf.~\eqref{eq:lambda-opened}), we define the
\emph{tropical lambda length} of $\gamma\in\barATSM\cup\BSM$
with respect to~$\barLL$ by
\begin{align}
\label{eq:c(gamma)}
c_\barLL(\gamma)&=\prod_{i\in I} q_i^{-l_{\barL_i}(\gamma)/2}\in\PP_\LL.\\
\shortintertext{If we set}
\label{eq:c(p)}
c_\barLL(p)
&= \prod_{i\in I} q_i^{-l_{\barL_i}(p)/2}\in\PP_\LL\\
\intertext{(cf.~\eqref{eq:lambda(p)}),
then by~\eqref{eq:signed-lambda-twist} the tropical lambda lengths satisfy}
\label{eq:c-twist}
c_\barLL(\psi_p\gamma) &= c_\barLL(p)^{n_p(\gamma)}\,c_\barLL(\gamma).
\end{align}
\end{definition}

\begin{remark}
  The minus signs in the exponents of
  equations~\eqref{eq:c(gamma)}-\eqref{eq:c(p)} are there because
  transformations of transverse measures involve maxima (as in
  \eqref{eq:l-trop-ptolemy} below) whereas the tropical
  semifield is defined using minima.
\end{remark}

We can also define, in a similar way, the quantities
$l_\barLL(\gamma)$ and $c_\barLL(\gamma)$
when $\gamma$ is a tagged loop based at a marked point and
enclosing a sole puncture/opening.

\begin{remark}
\label{rem:depend-on-lift}
Our definition of a tropical lambda length of a (lifted) tagged arc
depends on the choice of a lift~$\barLL$ of the multi-lamination~$\LL$.
Different choices result in different notions of a
tropical lambda length. However, they all differ from each other by
\emph{gauge transformations} which simultaneously rescale the lambda
lengths of all arcs incident to a given puncture, and do not affect
the resulting cluster algebra structure.

On the other hand, the tropical (or ordinary) lambda lengths of
boundary segments or holes
do not depend on the choice of a lift~$\barLL$.
Consequently, we can use notation $c_\LL(\beta)=c_\barLL(\beta)$ for
$\beta\in\BSM$, or $c_\LL(p)=c_\barLL(p)$ for $p\in\barM$.
The (tropical) lambda lengths of arcs not incident to punctures are
similarly independent on the choice of~$\barLL$.
\end{remark}

It is perhaps worth emphasizing that all our lambda lengths, whether
tropical or ordinary, do not depend on a (tagged) triangulation
containing the (tagged) arc in question.

The main property of the tropical lambda lengths is that they satisfy
the tropical version of the exchange relations
\eqref{eq:ptolemy-opened} and~\eqref{eq:exch-opened-digon}.

\begin{lemma}
\label{lem:ptolemy-x-trop}
On the opened surface $\barSM$, consider a quadrilateral with sides
$\alpha,\beta,\gamma,\delta$ and diagonals $\eta$ and~$\theta$, cf.\
Figure~\ref{fig:ptolemy}.
Assume that the tagging of the arcs in
$\{\alpha,\beta,\gamma,\delta,\eta,\theta\}$ is consistent at each
marked point. Then
\begin{equation}
\label{eq:ptolemy-x-trop}
c_\barLL(\eta)c_\barLL(\theta)
=c_\barLL(\alpha)c_\barLL(\gamma)\oplus c_\barLL(\beta)c_\barLL(\delta),
\end{equation}
where $\oplus$ denotes the tropical addition in~$\PP_\LL$.
\end{lemma}

\begin{proof}
It is immediate from the definition of the tropical semifield (see
Definition~\ref{def:tropical}) that it suffices to
prove~\eqref{eq:ptolemy-x-trop} in the case when $\barLL$ consists of a
single lamination~$\barL$.
In that case, \eqref{eq:ptolemy-x-trop} becomes
\begin{equation}
\label{eq:l-trop-ptolemy}
l_\barL(\theta)+l_\barL(\eta)
=
\max(l_\barL(\alpha)+l_\barL(\gamma), l_\barL(\beta)+l_\barL(\delta)).
\end{equation}
If all the arcs twist sufficiently far around the openings (if any)
containing their ends, then $l_\barL$ is an
intersection number and this is
a well known (and easy to check) relation
(see, e.g., \cite[Section~3]{fg-dual-teich}).
The general case follows by noticing that \eqref{eq:l-trop-ptolemy} is
invariant under twists around openings. 
\end{proof}

Lemmas \ref{lem:holed-monogon-trop}--\ref{lem:exch-digon-trop} below
are tropical analogues of Lemmas~\ref{lem:holed-monogon}--\ref{lem:exch-opened-digon}, 
respectively. 

\begin{lemma}
  \label{lem:holed-monogon-trop}
  Consider an opened monogon as shown in
  Figure~\ref{fig:opened-monogon}, with boundary marked point~$q$ and
  opened puncture~$C_p$.  If $\delta$ and~$\varrho$ are
  compatible parallel tagged arcs connecting $q$ and~$M_p$, one of
  them plain and one notched at~$M_p$, then we have
  \begin{equation}
    c_\barLL(\delta)c_\barLL(\varrho) = c_\barLL(\eta).
  \end{equation}
\end{lemma}

\begin{proof}
  Straightforward to check, as there are only three elementary
  laminations in a punctured monogon.
\end{proof}

\begin{lemma}
\label{lem:exch-digon-trop}
Consider an opened digon as shown in Figure~\ref{fig:opened-digon}.
Under the assumptions of
Lemma~\ref{lem:exch-opened-digon}, we have
\begin{equation}
\label{eq:exch-opened-digon-1-trop}
c_\barLL(\varrho)\, c_\barLL(\theta)
= c_\barLL(\alpha)\oplus c_\barLL(p)^{-1}\,c_\barLL(\beta).
\end{equation}
\end{lemma}

\begin{proof}
Just as Lemma~\ref{lem:exch-opened-digon} follows from
Lemma~\ref{lem:holed-monogon}, this follows
straightforwardly from Lemma~\ref{lem:holed-monogon-trop}.
Alternatively, this can be checked by examining the six elementary
laminations on a punctured digon from Figure~\ref{fig:lamin-digon}.
\end{proof}
In the language of transverse measures
(cf.~\eqref{eq:l-trop-ptolemy}), the relation
\eqref{eq:exch-opened-digon-1-trop} corresponds to the identity
\begin{equation}
\label{eq:l-trop-digon}
l_\barL(\varrho)+l_\barL(\theta)
=
\max(l_\barL(\alpha),
l_\barL(\beta)-l_\barL(p)).
\qedhere
\end{equation}

We then obtain an analogue of Theorem~\ref{th:cluster-lambda-opened}.

\begin{corollary}
  \label{cor:cluster-tropical-opened}
  Fix a multi-lamination~$\LL$ and its lift~$\barLL$.  For each
  $\gamma\in\ATSM$, fix an arbitrary lift~$\bargamma$. Then the
  tropical lambda lengths $\{c_\barLL(\bargamma)\}$ satisfy 
(\ref{eq:exchange-rel-aux}), the exchange relations in the tropical
semifield~$\PP_\LL$, for the same exchange matrices and the same
choices of coefficients as in Theorem~\ref{th:cluster-lambda-opened}.
\end{corollary}

We now relate the tropical shear coordinates from
Chapters~\ref{sec:shear} and~\ref{sec:shear-tagged} to lifted lambda
lengths.
The statement below is a version of a well known formula (see, e.g.,
\cite[p.~44]{thurston-minimal-stretch},
\cite[Section~4.6]{chekhov-penner}) relating shear coordinates to
transverse measures.

\begin{lemma}\label{lem:b-quadrilateral-sum}
  Let $T$ be a triangulation of $\SM$, let $L$ be
  a lamination on $\SM$, and let $\eta$ be an arc in~$T$
  that is not contained inside a self-folded triangle.  Let
  $\alpha$, $\beta$, $\gamma$, $\delta$ be the arcs on the boundary of
  the quadrilateral containing $\eta$, arranged as in
  Figure~\ref{fig:ptolemy}, and let $\baralpha_0$, $\barbeta_0$,
  $\bargamma_0$, $\bardelta_0$ be compatible 
  lifts of them to $\barSM$, in the sense that they form a 
  quadrilateral on $\barSM$.
  Then for any lift $\barL$ of~$L$, we have
  \begin{equation}
    \label{eq:b-quadrilateral-sum}
    2b_\eta(T,L) = -l_\barL(\baralpha_0) + l_\barL(\barbeta_0) -
      l_\barL(\bargamma_0) + l_\barL(\bardelta_0).
  \end{equation}
\end{lemma}

Note that some of $\alpha$, $\beta$, $\gamma$, $\delta$ may not be
tagged arcs in the strict sense of Definition~\ref{def:tagged-arc}, as
they may enclose a
once-punctured monogon; still, $l_\barL$~is well defined.

\begin{proof}
  This is a consequence of Definition~\ref{def:shear}, as follows.
  Assume that $\barL$ is twisted sufficiently far in the sense of
  Definition~\ref{def:trop-lambda}.
  Each `S'-shaped intersection of $\barL$ with $\eta$ contributes $+2$ to the
  right hand side, each `Z'-shaped intersection contributes $-2$,
  while intersections of $\barL$ with the quadrilateral that cut off a
  corner does not effect the right hand side.  Changing
  the lift $\barL$ of~$L$ effectively changes the number of
  intersections that cut
  off corners. In
  particular, we can modify $\barL$ so that it twists sufficiently far
  without changing the right hand side.
\end{proof}

\begin{lemma}\label{lem:b-quadrilateral-sum-tagged}
  For $T$ a tagged triangulation of $\SM$ and $\eta\in T$ a tagged arc
  that is not parallel to any other arc, let $\alpha$, $\beta$,
  $\gamma$, $\delta$ be the tagged arcs on the boundary of the
  quadrilateral containing~$\eta$, arranged as in
  Figure~\ref{fig:ptolemy}.  (Some or all of $\alpha,
  \dots,\delta$ may be curves enclosing punctured monogons and not
  in~$T$ itself.)  Let $\baralpha_0$, $\barbeta_0$,
  $\bargamma_0$, $\bardelta_0$ be compatible 
  lifts to $\barSM$.
  Then for any lift $\barL$ of~$L$,
  \eqref{eq:b-quadrilateral-sum} holds.
\end{lemma}

\begin{proof}
  This is a consequence of Lemma~\ref{lem:b-quadrilateral-sum} and Part~(i) of
  Definition~\ref{def:shear-tagged}, as follows.  Let $L_1$ and $L_2$
  be two laminations on $\SM$ that differ only in the direction of
  spiraling
  at a puncture~$p$.  Then, if $\barL_1$ is a lift of
  $L_1$, we can find a lift $\barL_2$ of $L_2$ by simply changing the
  orientation of~$C_p$, without changing the actual curves
  corresponding to $\barL_1$ in $\barSM$.  Then, for any two lifted arcs
  $\baralpha_1$ and $\baralpha_2$ which are identical except that
  $\baralpha_1$ is plain at~$p$ and $\baralpha_2$ is notched at~$p$,
  by Definition~\ref{def:transverse-measure} we have
  \begin{equation}
    l_{\barL_1}(\baralpha_2) = l_{\barL_2}(\baralpha_1) +
      \abs{l_{\barL_1}(p)},
  \end{equation}
  as when $\baralpha_1$ is twisted sufficiently far with respect to
  $\barL_1$, $\baralpha_2$ is twisted sufficiently far with
  respect to $\barL_2$, as both are trying to twist in the same
  direction, $\baralpha_1$ with the orientation of $C_p$ and
  $\baralpha_2$ against the reversed orientation of~$C_p$.

  Thus, for each vertex of the quadrilateral at which there are notches,
  we can simultaneously remove them and reverse the direction
  of spiraling of~$L$ without changing the sum on the right
  of~\eqref{eq:b-quadrilateral-sum}, in accordance with Part~(i) of
  Definition~\ref{def:shear-tagged}.
\end{proof}

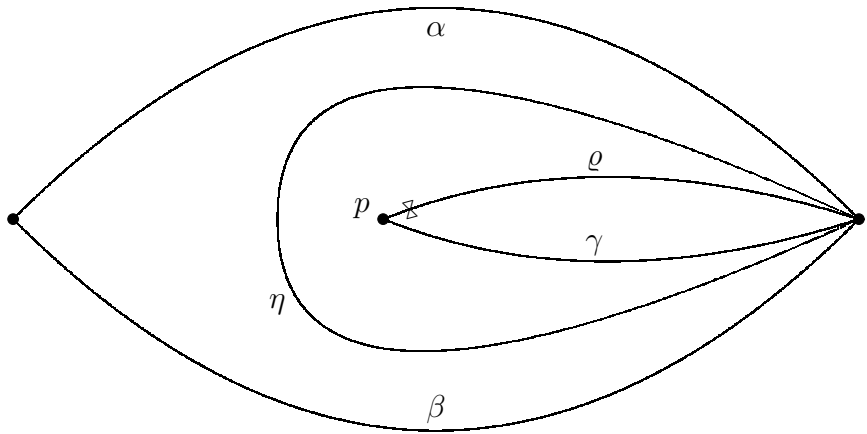
\begin{figure}[htbp]
\begin{center}
\setlength{\unitlength}{4pt}
\begin{picture}(80,40)(0,0)
\thinlines
\put(35,20){\circle*{1}}
\put(0,20){\circle*{1}}
\put(80,20){\circle*{1}}

\qbezier(35,20)(55,12)(80,20)
\qbezier(35,20)(55,28)(80,20)

\qbezier(25,20)(25,-5)(80,20)
\qbezier(25,20)(25,45)(80,20)

\qbezier(0,20)(20,40)(40,40)
\qbezier(0,20)(20,0)(40,0)
\qbezier(40,0)(60,0)(80,20)
\qbezier(40,40)(60,40)(80,20)

\put(37.5,20.9){\rotatebox{-75}{\makebox(0,0){$\notch$}}}

\put(33,21){\makebox(0,0){$p$}}

\put(40,38){\makebox(0,0){$\alpha$}}
\put(40,2){\makebox(0,0){$\beta$}}
\put(25,12){\makebox(0,0){$\eta$}}
\put(55,17.5){\makebox(0,0){$\gamma$}}
\put(55,25.5){\makebox(0,0){$\varrho$}}


\end{picture}
\end{center}
\caption{Proof of Lemma~\ref{lem:b-digon-sum}.}
\label{fig:b-digon-sum}
\end{figure}

\begin{lemma}\label{lem:b-digon-sum}
  Let $T$ be a tagged triangulation of $\SM$, $L$ be
  a lamination on $\SM$, and let $\varrho$ and~$\gamma$ be parallel arcs
  in~$T$ that differ in tagging at one
  end, call it~$p$.  The arcs $\varrho$ and~$\gamma$ are contained inside a digon in~$T$
  where the other sides are labeled $\alpha$ and~$\beta$, as in
  Figure~\ref{fig:b-digon-sum}.  Let $\baralpha_0$ and
  $\barbeta_0$ be compatible lifts of $\alpha$ and $\beta$ to
  $\barSM$.
  Then for any lift $\barL$ of~$L$,
  \begin{equation}
    \label{eq:b-digon-sum}
    2b_\varrho(T, L) = -l_\barL(\baralpha_0) + l_\barL(\barbeta_0) +
    n_p(\varrho) l_\barL(p).
  \end{equation}
\end{lemma}

As before, $\alpha$ and $\beta$ may enclose a once-punctured monogon
and so not be tagged arcs.

\begin{proof}
  Suppose that~$\varrho$ is notched at~$p$,
and suppose there is an ordinary
  triangulation~$T^\circ$ so that $\tau(T^\circ) = T$.  Let $\eta$ be
  the arc that projects
  to~$\varrho$. Then by Part~(ii) of
  Definition~\ref{def:shear-tagged}, $b_\varrho(T,L)=b_\eta(T^\circ,L)$.
  Now by Lemma~\ref{lem:b-quadrilateral-sum},
  \begin{equation}
    b_\eta(T^\circ, L) = -l_\barL(\baralpha_0) + l_\barL(\barbeta_0)
      - l_\barL(\bargamma_0) + l_\barL(\bardelta_0),
  \end{equation}
  where $\baralpha_0,\barbeta_0,\bargamma_0$, and $\bardelta_0=\psi_p(\bargamma_0)$ 
are compatible lifts of
  $\alpha,\beta,\gamma$ and again~$\gamma$.
The result then follows by
  equation~\eqref{eq:l(p)-twist-trop}.
  
  If the tags differ from the case above, then the result follows
  by applying tag-changing transformations and changing the spiraling
  of~$L$, as in Lemma~\ref{lem:b-quadrilateral-sum-tagged}.
\end{proof}

We can now give a conceptual proof that tropical shear coordinates behave as
expected under mutation.  For any triangulation~$T$ (ordinary or
tagged), multi-lamination~$\LL$, and arc $\eta\in T$, define
\begin{equation}
  r_\eta(T,\LL) = \prod_{i\in I} q_i^{-b_\eta(T,L_i)}. 
\end{equation}

\begin{proof}[Proof of Theorem~\ref{th:thurston-2}]
  Let us fix an ordinary triangulation~$T$ without self-folded triangles and
  a multi-lamination~$\LL$ with lift~$\barLL$.  For each arc~$\gamma\in
  T$, fix a lift~$\bargamma$.  In the $r$ variables,
  \eqref{eq:b-quadrilateral-sum} becomes
  \begin{equation}
    \label{eq:r-trans}
    r_\eta(T,\LL) = c_\barLL(\baralpha_0)^{-1}\cdot c_\barLL(\barbeta_0)
      \cdot c_\barLL(\bargamma_0)^{-1}\cdot c_\barLL(\bardelta_0).
  \end{equation}
Rewriting equation~\eqref{eq:r-trans} in terms of our
  chosen lifts
  $\baralpha,\dots,\bardelta$ (which are not necessarily compatible)
  and the exchange matrix $B(T)$, we
  obtain 
  \begin{equation}\label{eq:r-is-y}
    r_\eta(T,\LL) = \frac{p^+_\eta}{p^-_\eta}
      \prod_{\theta\in T} c_\barLL(\bartheta)^{B(T)_{\theta\eta}},
  \end{equation}
  where the product runs over all edges $\theta\in T$, and $p^\pm_\eta$ are
  the coefficients from Theorem~\ref{th:cluster-lambda-opened}.
  (An edge~$\theta$ in~$T$ that is not adjacent to~$\eta$ will not
  contribute to the product, as then $B(T)_{\theta\eta} = 0$.  The
  factor $p^+_\eta/p^-_\eta$ makes up for taking the lifts
  $\baralpha,\dots,\bardelta$ rather than
  $\baralpha_0,\dots,\bardelta_0$.)

  Thus the variables $r_\eta(T,\LL)$ are defined just like the $\hat
  y$ variables in~\eqref{eq:yhat}, but with the variables
  $c_\barLL$.  By Corollary~\ref{cor:cluster-tropical-opened}, the
  $c_\barLL$ transform as cluster variables in the tropical
  semifield $\PP_\LL$.  Thus Proposition~\ref{pr:yhat} applies, and
  the $r_\eta(T,\LL)$ transform according to
  equation~\eqref{eq:yhat-mutation} (interpreted tropically).  As
  noted in Remark~\ref{rem:mutation-yhat}, this is how the
  coefficients in the $B$ matrix transform, as claimed in the statement of 
  Theorem~\ref{th:thurston-2}.
\end{proof}

\begin{proof}[Proof of Theorem~\ref{th::thurston-2-tagged}]
  For the more general setting of Theorem~\ref{th::thurston-2-tagged},
  we must extend the above arguments to the case of tagged triangulations.
  As in
  Definition~\ref{def:shear-tagged}, there is an untagged
  triangulation $T^\circ$ so that $T$ is obtained from $T^\circ$ by
  applying tag-changing transformations to
  $\tau(T^\circ)$.  An arc~$\eta \in T$ can have a parallel copy with
  different tagging (in which case the corresponding arc $\eta^\circ
  \in T^\circ$ is part of a self-folded triangle), or not.  In the
  second case, by Lemma~\ref{lem:b-quadrilateral-sum-tagged},
  \eqref{eq:r-trans} holds (where, as before,
  $\baralpha_0,\dots,\bardelta_0$ are compatible lifts of the
  quadrilateral containing~$\eta$).  If one of $\alpha,\dots,\delta$
  is a tagged curve enclosing a once-punctured monogon,
  Lemma~\ref{lem:holed-monogon-trop} applies.  This combines with the
  behaviour of $B(T)$ in the case of self-folded triangles (see
  \cite[Definitions~4.1 and~9.7]{cats1}) to show that \eqref{eq:r-is-y}
  is true in this case as well.  (As before, the factor
  $p_\eta^+/p_\eta^-$ makes up for taking the original lifts rather
  than compatible lifts.)

  On the other hand, let $\varrho\in T$ be a tagged arc with a parallel
  copy~$\gamma$.  Then we can apply Lemma~\ref{lem:b-digon-sum} to
  conclude
  \begin{equation}\label{eq:r-trans-bigon}
    r_\varrho(T, \LL) = c_\barLL(\baralpha_0)^{-1} \cdot
      c_\barLL(\barbeta_0) \cdot c_\barLL(p)^{n_p(\varrho)}.
  \end{equation}
  But this is yet another form of~\eqref{eq:r-is-y} for the
  arc~$\varrho$.  (The last factor
  in~\eqref{eq:r-trans-bigon} becomes part of the coefficient factor
  $p_\varrho^+/p_\varrho^-$.)

  Thus for all arcs in~$T$, \eqref{eq:r-is-y} holds and the $r_\eta$
  are defined like the $\hat y$ variables with respect to the tropical
  lambda lengths $c_\barLL$.  The result follows by
  Remark~\ref{rem:mutation-yhat}.
\end{proof}

\chapter{Laminated Teichm\"uller spaces}
\label{sec:laminated-teich}

In this chapter, we use the notions developed in previous
chapters---specifically, ordinary and tropical lambda lengths of
tagged arcs on an opened surface---to present our main construction, a
geometric realization of cluster algebras associated with surfaces.
The main idea is rather natural.
As the lifts~$\bargamma$ of an arc~$\gamma$ vary,
the corresponding tropical lambda lengths, just like the ordinary
ones, form a geometric progression (cf.~\eqref{eq:signed-lambda-twist}
and~\eqref{eq:c-twist}).
After making sure that the ratios of the two progressions (ordinary and
tropical) coincide,
we proceed by dividing an ordinary (rescaled) lambda length
of~$\bargamma$ by a tropical one, thus obtaining an invariant
of~$\gamma$ that plays the role of a cluster variable.

\begin{definition}[\emph{Laminated Teichm\"uller space}]
\label{def:teich-lamin}
Let $\LL=(L_i)_{i\in I}$ be a multi\hyp lamination in $\SM$.
The \emph{laminated Teichm\"uller space} $\TTSML$ is defined as
follows. A point $(\sigma,q)\in\TTSML$
is a decorated hyperbolic structure $\sigma\in\TTSM$ together with a
collection of positive real weights~$q=(q_i)_{i\in I}$
which are chosen so that the following boundary conditions hold:
\begin{itemize}
\item
for each boundary segment~$\beta\in\BSM$,
we have $\lambda_\sigma(\beta)=c_\LL(\beta)$;
\item
for each hole~$C_p$, with $p\in\barM$,
we have $\lambda_\sigma(p)=c_\LL(p)$.
\end{itemize}
In these equations, the quantities $c_\LL(\beta)$ and
$c_\LL(p)$ (cf.\ Remark~\ref{rem:depend-on-lift})
are given by the formulas \eqref{eq:c(gamma)} and
\eqref{eq:c(p)}, with each $q_i$ specialized to the given positive real
value.
In more concrete terms, we require that for each boundary segment~$\beta$,
\begin{equation}
\label{eq:l=ltrop}
l_\sigma(\beta)=-\sum_{i\in I} l_{\barL_i}(\beta) \ln(q_i),
\end{equation}
and similarly with $\beta$ replaced by~$p$.
Informally, our boundary conditions \eqref{eq:l=ltrop} ask
that for each hole (resp., boundary segment),
the total weighted sum of its transverse measures
with respect to the laminations in~$\LL$ is the negative of the ordinary
hyperbolic length, with respect to~$\sigma$,
of the hole (resp., segment between horocycles).
\end{definition}

We next coordinatize the laminated Teichm\"uller space $\TTSML$.
A system of coordinates will include the weights of laminations plus
the lambda lengths of the arcs in a (tagged) triangulation.

\pagebreak[3]

\begin{proposition}
\label{pr:param-teich-lamin}
Let $\LL=(L_i)_{i\in I}$ be a multi\hyp lamination in~$\SM$.
The laminated Teichm\"uller space $\TTSML$ can be coordinatized
as follows.
Fix an ideal (or tagged) triangulation~$T$ of~$\SM$.
Choose a lift of each of the $n$ arcs $\gamma \in T$ to an arc
  $\bargamma \in \barAPSM$.
Then the map
  \[
\Psi: \TTSML \to \RR_{>0}^{n+\abs{I}}
  \]
defined by
\[
\Psi(\sigma,q)
=
(\lambda_\sigma(\bargamma))_{\gamma \in T} \times q
\]
is a homeomorphism.
\end{proposition}

\begin{proof}
This is a direct consequence of Proposition~\ref{pr:param-teich-opened} and
Definition~\ref{def:teich-lamin}.
\end{proof}

Proposition~\ref{pr:param-teich-lamin} enables us to view the
coordinate functions $q_i$ and $\lambda(\bargamma)$ as ``variables,''
similarly to our treatment of the earlier versions of the
Teichm\"uller space.

We next define the quantities that will play the role of cluster
variables in our main construction.

\begin{definition}[\emph{Laminated lambda lengths}]
\label{def:lamin-lambda}
Let us fix a lift~$\barLL$ of a multi-lamination~$\LL$.
For a tagged arc $\gamma\in\ATSM$, the \emph{laminated lambda length}
$x_\barLL(\gamma)$ is a function on the laminated
Teichm\"uller space $\TTSML$ defined by
\begin{equation}
\label{eq:x-lamin}
x_\barLL(\gamma)=\lambda(\bargamma)/c_\barLL(\bargamma),
\end{equation}
where $\bargamma$ is an arbitrary lift of~$\gamma$
and $c_\barLL(\bargamma)$ is the
tropical lambda length of Definition~\ref{def:trop-lambda}.
The value of $x_\barLL(\gamma)$ does not depend on the choice of
the lift~$\bargamma$ since $\lambda(\bargamma)$ and $c_\barLL(\bargamma)$ rescale
by the same factor $c_\LL(p)^{\pm 1}=\lambda_\sigma(p)^{\pm 1}$
(see Definition~\ref{def:teich-lamin}) as $\bargamma$
twists around the opening~$C_p$ (cf.\ 
\eqref{eq:signed-lambda-twist} and~\eqref{eq:c-twist}).

On the other hand, $x_\barLL(\gamma)$ does depend on the choice of the
lifted multi\hyp lamination~$\barLL$.
This will not create problems as the resulting cluster structure will be
unique up to gauge transformations (see Remark~\ref{rem:depend-on-lift}),
hence up to a canonical isomorphism.
\end{definition}

\begin{remark}
The laminated lambda lengths $x_\barLL(\gamma)$ can be intuitively
interpreted as follows.
Suppose that the lift $\bargamma$ is such that it twists sufficiently
far around each of its ends lying on opened components.
(Since
some of the $\barL_i$ may have opposite orientation on a given $C_p$, it may
be impossible to satisfy this condition for all~$i$
simultaneously, as
$\bargamma$ may be required to spiral in opposite
directions; still, let us assume that we can.)
Then we combine the definition \eqref{eq:x-lamin} with
\eqref{eq:lambda-opened} and
\eqref{eq:c(gamma)} to get
\begin{equation}
\label{eq:x-lamin-unravel}
x_\barLL(\gamma)=e^{l(\bargamma)/2}\,
\prod_{i\in I} q_i^{l_{\barL_i}(\bargamma)/2}.
\end{equation}
Thus $x_\barLL(\gamma)$ is obtained by exponentiating
a sum of two kinds  of contributions:
\begin{itemize}
\item
the hyperbolic length of the lifted arc~$\bargamma$ between
appropriate horocycles (``the cost of fuel while traveling
along~$\bargamma$''), and
\item
a fixed contribution, depending on the~$q_i$,
associated with each crossing of~$\barL_i$ by~$\bargamma$
(``the tolls'').
\end{itemize}
\end{remark}

\begin{corollary}
\label{cor:param-lamin-cluster}
Let $\LL=(L_i)_{i\in I}$ be a multi\hyp lamination in~$\SM$.
For a tagged triangulation~$T$ and any choice of lift $\barLL$ of
$\LL$, the map
\begin{align*}
\TTSML &\to \RR_{>0}^{n+\abs{I}} \\
(\sigma,q)&\mapsto
(x_\barLL(\gamma))_{\gamma \in T} \times q
\end{align*}
is a homeomorphism.
\end{corollary}

\begin{proof}
This is a corollary of Proposition~\ref{pr:param-teich-lamin} and
Definition~\ref{def:lamin-lambda}.
\end{proof}

We are finally prepared to present our main construction.

\pagebreak[3]

\begin{theorem}
\label{th:cluster-lamin}
For a given multi-lamination $\LL=(L_i)_{i\in I}$,
there exists a unique normalized exchange pattern $(\Sigma_T)$
of geometric type with the following properties:
\begin{itemize}
\item
the coefficient semifield is the tropical semifield
$\,\PP_\LL=\Trop(q_i:i\in I)$;
\item
the cluster variables $x_\LL(\gamma)$ are labeled by the tagged arcs
$\gamma\in\ATSM$;
\item
the seeds $\Sigma_T=(\xx(T),\tilde B(T,\LL))$
are labeled by the tagged triangulations~$T$;
\item
the exchange graph is the graph $\ESM$ of
tagged flips, see Definition~\ref{def:ESM};
\item
each cluster $\xx(T)$ consists of cluster variables
$x_\LL(\gamma)$, for $\gamma\in T$;
\item
the ambient field $\Fcal$ is generated over $\PP_\LL$
by any given cluster~$\xx(T)$;
\item
the extended exchange matrix $\tilde B(T,\LL)$
is described in Definition~\ref{def:tilde-BL-tagged}.
\end{itemize}
This exchange pattern has a positive realization (see
Definition~\ref{def:positive-realization})
by functions on the laminated Teichm\"uller space $\TTSML$.
To obtain this realization,
choose a lift $\barLL$ of the multi-lamination~$\LL$;
then represent each cluster variable $x_\LL(\gamma)$ by the
corresponding laminated lambda length $x_\barLL(\gamma)$,
and each coefficient variable~$q_i$
by the corresponding function on $\TTSML$.
\end{theorem}

\begin{proof}
We already know (see Theorem~\ref{th::thurston-2-tagged}) that the
matrices $\tilde B(T,\LL)$ are related to each other by mutations
associated with the tagged flips.
In view of Propositions~\ref{prop:positive-realization}
and~\ref{pr:param-teich-lamin},
all we need to prove is that the laminated lambda lengths
$x_\barLL(\gamma)$ satisfy the exchange relations encoded by the
matrices $\tilde B(T,\LL)$.

Fix arbitrary lifts $\bargamma$ of all tagged arcs $\gamma\in\ATSM$.
By Theorem~\ref{th:cluster-lambda-opened}, the lambda lengths of
lifted arcs $x(\bargamma)$ form a non-normalized exchange pattern on
$\ESM$. Moreover this statement remains true
if we view $x(\bargamma)$ as a function on $\TTSML$ rather than $\TTSM$,
replacing the coefficient group $\PP\SM$
(cf.~\eqref{eq:coeff-group-non-normalized}) 
by $\Trop(q_i)$.
Indeed, the monomial mutation rules
\eqref{eq:p-mutation1}--\eqref{eq:p-mutation2} are satisfied by the
coefficients in exchange relations for $x(\bargamma)$'s viewed as
functions on $\TTSM$, and therefore they are satisfied after the
lambda lengths of boundary segments
and holes are replaced by monomials in the lamination weights.

By Lemmas \ref{lem:ptolemy-x-trop}
and~\ref{lem:exch-digon-trop}, their tropical counterparts
$c_\barLL(\bargamma)$ satisfy the tropical versions of the same
exchange relations.
By Proposition~\ref{pr:rescale-normalize}, this is exactly what is
needed in order for the rescaled lambda lengths
$x(\bargamma)/c_\barLL(\bargamma)$ to form a normalized exchange
pattern.

It remains to verify that the extended exchange matrices of
this exchange pattern are the matrices $\tilde B(T,\LL)$.
(This in itself implies the uniqueness statement in the theorem.)
In fact, by Theorem~\ref{th::thurston-2-tagged}, it suffices to do
this for \emph{some} triangulation, say an ordinary triangulation~$T$
without self-folded triangles.
In this setting, we need to verify that the coefficients of
the exchange relations associated with the flips from~$T$
(recall that these coefficients are given by~\eqref{eq:rescale-p})
coincide with the coefficients encoded by~$\tilde B(T,\LL)$.

These exchange relations are obtained by dividing the Ptolemy
relations~\eqref{eq:ptolemy-opened} by their tropical
counterparts~\eqref{eq:ptolemy-x-trop}.
If arcs $\alpha,\beta,\gamma,\delta\in T$ form a quadrilateral
with diagonals $\eta\in T$ and~$\theta$, as in Figure~\ref{fig:ptolemy},
then the corresponding exchange relation is
\[
x_\barLL(\eta)\,x_\barLL(\theta)
= p_\eta^+ \,x_\barLL(\alpha)\,x_\barLL(\gamma)
+ p_\eta^- \,x_\barLL(\beta) \,x_\barLL(\delta),
\]
where
\begin{equation}
\label{eq:p-eta}
p_\eta^+=
\frac{c_\barLL(\baralpha)\,c_\barLL(\bargamma)}{c_\barLL(\bareta)\,c_\barLL(\bartheta)},\qquad
p_\eta^-=
\frac{c_\barLL(\barbeta)\,c_\barLL(\bardelta)}{c_\barLL(\bareta)\,c_\barLL(\bartheta)},
\end{equation}
and $\baralpha,\barbeta,\bargamma,\bardelta,\bareta,\bartheta$ are
lifts of the corresponding arcs in~$T$, coordinated so that
$\baralpha,\barbeta,\bargamma,\bardelta$ form a quadrilateral with
diagonals $\bareta$ and~$\bartheta$.
If, for example, $\alpha$~is a boundary segment rather than an arc,
then $x(\alpha)=\lambda(\alpha)=c_\LL(\alpha)=c_\barLL(\alpha)$
(see Definition~\ref{def:teich-lamin}), and the exchange relation
becomes
\[
x_\barLL(\eta)\,x_\barLL(\theta)
= p_\eta^+ \,x_\barLL(\gamma)
+ p_\eta^- \,x_\barLL(\beta) \,x_\barLL(\delta),
\]
with the coefficients $p_\eta^\pm$ given by \eqref{eq:p-eta} as
before.
Since these coefficients satisfy the normalization condition
$p_\eta^+\oplus p_\eta^-=1$ (by
Proposition~\ref{pr:rescale-normalize}, or by
\eqref{eq:ptolemy-x-trop}), the claim reduces to showing that
their ratio coincides with the Laurent monomial encoded by the
appropriate column of the matrix~$\tilde B(T,\LL)$. That is, we need
to check that
\[
\frac{p_\eta^+}{p_\eta^-}=\frac{c_\barLL(\baralpha)\,c_\barLL(\bargamma)}
{c_\barLL(\barbeta)\,c_\barLL(\bardelta)}
=\prod_{i\in I} q_i^{b_\eta(T,L_i)},
\]
where, as in~\eqref{eq:bijTL}, $b_\eta(T,L_i)$ denotes the shear
coordinate of the arc~$\eta$ with respect to the lamination~$L_i$.
In view of~\eqref{eq:c(gamma)}, this is equivalent
to~\eqref{eq:b-quadrilateral-sum}.

The latter is a version of a well known formula (see, e.g.,
\cite[p.~44]{thurston-minimal-stretch},
\cite[Section~4.6]{chekhov-penner}) relating shear coordinates to
transverse measures.
\end{proof}

\begin{definition}
\label{def:ca-multi-lamin}
Let $\LL$ be a multi-lamination on~$\SM$,
a bordered surface with marked points.
Then the cluster algebra of geometric type described in
Theorem~\ref{th:cluster-lamin} is denoted
by~$\Acal(\Surf,\Mark,\LL)$.
\end{definition}

\begin{remark}
\label{rem:border-coeff-as-lamin}
The exchange pattern in Theorem~\ref{th:cluster-lambda}---with the
coefficient variables corresponding to the boundary segments in $\BSM$---can be
obtained as a particular case of the main construction of this
chapter.
The corresponding multi-lamination
$\LL=\{L_\beta\}_{\beta\in\BSM}$
contains one lamination $L_\beta$ for each boundary segment
$\beta\in\BSM$.
The lamination $L_\beta$ consists of a single curve (also denoted
by~$L_\beta$) defined as follows.
Let $p,q\in\Mark$ be the endpoints of~$\beta$.
If $p=q$ (so that $\beta$ lies on a boundary component with a single
marked point), then $L_\beta$ is a closed curve in~$\Surf$
encircling~$\beta$.
Otherwise, let $p'$ and $q'$ be two points
on~$\partial\Surf\setminus\beta$ located
very close to $p$ and~$q$, respectively.
Then $L_\beta$ connects $p'$ and~$q'$ within a small neighborhood
of~$\beta$.

See Example~\ref{example:G2n+3} for a particular instance of this
construction.
\end{remark}

\begin{proof}[Proof of Theorem~\ref{th:patterns-on-SM}]

Uniqueness of an exchange pattern described in the theorem is clear:
the mutation rules determine all seeds uniquely from the initial one.
To prove existence, we need to show that a sequence of mutations that
returns us back to the same tagged triangulation recovers the
original seed.
To use the terminology of \cite[Definition~4.5]{ca4}, we want to show
that the exchange graph of our pattern is \emph{covered} by (and in fact is
identical to) the graph $\ESM$ of tagged flips.

We first reduce this claim to the case of patterns of geometric type.
Indeed, \cite[Theorem~4.6]{ca4} asserts that 
the exchange graph of any cluster algebra~$\Acal$ is
covered by the exchange graph of a particular cluster algebra 
of geometric type, namely one that has the same
exchange matrices as~$\Acal$, and whose coefficients are 
\emph{principal} at an arbitrarily chosen seed. 
(See \cite[Definition~3.1]{ca4} or Definition~\ref{def:principal} below.) 

For an exchange pattern of geometric type,
we produce a positive realization
(cf.\ Definition~\ref{def:positive-realization})
using the construction of Theorem~\ref{th:cluster-lamin}.
The key role in the argument is played by our generalization of
Thurston's coordinatization theorem (Theorems~\ref{th:thurston-1}
and~\ref{th:thurston-param-tagged}), which guarantees that
we can get any coefficients (of geometric type) by making an
appropriate (unique) choice of a multi-lamination~$\LL$.

It remains to show that all laminated lambda lengths
$x_\barLL(\gamma)$, for $\bargamma\in\APSM$, are distinct.
There are several ways to prove this; here we give a sketch of one such
proof.
It is easy to see that if in some normalized cluster algebra,
a particular sequence of mutations yields a
seed containing a cluster variable equal to one of the cluster
variables of the initial cluster, then the same phenomenon must hold
in the cluster algebra defined by the same initial exchange matrix over the
one-element semifield~$\{1\}$.
It is therefore enough to check the claim in the case of
trivial coefficients. In our setting, this case can be viewed as a
special instance of the construction presented in
Theorem~\ref{th:cluster-lambda}, with the length of each boundary
segment equal to~$1$.
In that instance, the cluster variables $x(\gamma)$ are the
lambda lengths of certain (pairwise non-isotopic) geodesics between
appropriate horocycles.

We can use the torus action associated with changing the horocycles to
show that we only need to consider pairs of tagged arcs connecting the
same horocycles (around the same marked points).

So suppose we have two arcs, $\gamma_1$ and~$\gamma_2$, connecting the
same pair of marked points $p$ and~$q$.  If $p \ne q$, consider a loop
$\ell$ that surrounds $\gamma_1$, and consider a family of hyperbolic
metrics that pinch at $\ell$ (i.e., the length of $\ell$ goes to~$0$). Then
the length of $\gamma_2$ goes to~$\infty$ whereas the decorations can be
chosen to keep the length of
$\gamma_1$ bounded.

If $p = q$, consider instead the two loops $\ell_1$ and $\ell_2$
obtained by pushing $\gamma_1$ around the puncture $p$ on the two
different sides.  Again, consider a family of hyperbolic metrics that
pinch at $\ell_1$ and~$\ell_2$.  In this family, the length of
$\gamma_2$ goes to~$\infty$ while the length of $\gamma_1$ can remain
bounded.

(We thank Y.~Eliashberg for suggesting this last argument.)
\end{proof}

\begin{proof}[Proof of Corollary~\ref{cor:structural-clust}]
By Theorem~\ref{th:patterns-on-SM},
the cluster variables in $\Acal$ are in bijection with tagged arcs, and
the seeds with tagged triangulations.
The claims in the Corollary then follow from the basic properties of
the tagged arc complex (see \cite[Section~7]{cats1}).
\end{proof}

\chapter{Topological realizations of some coordinate rings}
\label{sec:topol-coord-ring}

The main construction of Chapter~\ref{sec:laminated-teich} can be used
to produce topological realizations of (well known) cluster structures
in coordinate rings of various algebraic varieties. Several such
examples are presented in this chapter.

In each of these examples, the corresponding generalized decorated
Teichm\"uller space $\TTSML$ is naturally interpreted as the
\emph{totally positive} part of the corresponding \emph{cluster
  variety}~$X$, the spectrum of the associated cluster algebra. 
In other words, $\TTSML$ can be identified with the set of those points
in~$X$ at which all cluster and coefficient variables (equivalently,
those belonging to a given extended cluster) take positive values. 
In each case, one recovers the usual notion of total positivity of
matrices or its well known extensions to Grassmannians and
other~$G/P$'s (see, e.g., \cite{dbc, tptp, lusztig-intro} and
references therein). 

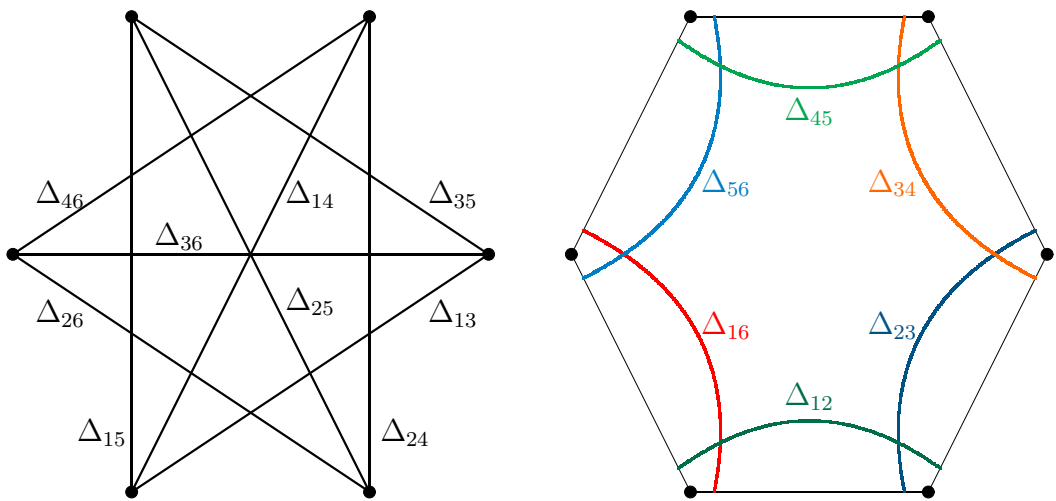
\begin{figure}[htbp]
\begin{center}
\setlength{\unitlength}{4.5pt}
\begin{picture}(40,42)(0,0)
\thinlines

  \put(0,20){\circle*{1}}
  \put(10,0){\circle*{1}}
  \put(10,40){\circle*{1}}
  \put(30,0){\circle*{1}}
  \put(30,40){\circle*{1}}
  \put(40,20){\circle*{1}}

\thicklines
  \put(0,20){\line(3,2){30}}
\put(4,25){\makebox(0,0){$\Delta_{46}$}}
  \put(0,20){\line(3,-2){30}}
\put(4,15){\makebox(0,0){$\Delta_{26}$}}
  \put(30,0){\line(0,1){40}}
\put(33,5){\makebox(0,0){$\Delta_{24}$}}

  \put(40,20){\line(-3,2){30}}
\put(37,25){\makebox(0,0){$\Delta_{35}$}}
  \put(40,20){\line(-3,-2){30}}
\put(37,15){\makebox(0,0){$\Delta_{13}$}}
  \put(10,0){\line(0,1){40}}
\put(7.5,5){\makebox(0,0){$\Delta_{15}$}}

  \put(10,0){\line(1,2){20}}
\put(25,25){\makebox(0,0){$\Delta_{14}$}}
  \put(10,40){\line(1,-2){20}}
\put(25,16){\makebox(0,0){$\Delta_{25}$}}
  \put(0,20){\line(1,0){40}}
\put(14,21.5){\makebox(0,0){$\Delta_{36}$}}

\end{picture}
\qquad
\setlength{\unitlength}{4.5pt}
\begin{picture}(40,42)(0,0)
\thinlines
  \put(10,0){\line(-1,2){10}}
  \put(40,20){\line(-1,2){10}}
  \put(10,0){\line(1,0){20}}
  \put(10,40){\line(1,0){20}}
  \put(0,20){\line(1,2){10}}
  \put(30,0){\line(1,2){10}}

  \put(0,20){\circle*{1}}
  \put(10,0){\circle*{1}}
  \put(10,40){\circle*{1}}
  \put(30,0){\circle*{1}}
  \put(30,40){\circle*{1}}
  \put(40,20){\circle*{1}}

\thicklines
\darkred{\qbezier(12,0)(15,15)(1,22)
\put(13,14){\makebox(0,0){$\Delta_{16}$}}
}
\lightblue{\qbezier(12,40)(15,25)(1,18)
\put(13,26){\makebox(0,0){$\Delta_{56}$}}
}
\darkblue{\qbezier(28,0)(25,15)(39,22)
\put(27,14){\makebox(0,0){$\Delta_{23}$}}
}
\lightred{\qbezier(28,40)(25,25)(39,18)
\put(27,26){\makebox(0,0){$\Delta_{34}$}}
}
\darkgreen{\qbezier(9,2)(20,10)(31,2)
\put(20,8){\makebox(0,0){$\Delta_{12}$}}
}
\lightgreen{\qbezier(9,38)(20,30)(31,38)
\put(20,32){\makebox(0,0){$\Delta_{45}$}}
}

\end{picture}
\end{center}
\caption{Representing the cluster structure on the Grassmannian
  $\operatorname{Gr}_{2,6}$ by a multi-lamination of a hexagon.
}
\label{fig:G26-lamin}
\end{figure}

\begin{example}[\emph{Grassmannians
    $\operatorname{Gr}_{2,n+3}(\CC)$}]
\label{example:G2n+3}
We already considered this cluster algebra in
Example~\ref{ex:unpunctured-disk}.
The homogeneous coordinate ring of the Grassmannian
$\operatorname{Gr}_{2,n+3}(\CC)$ (with respect to its Pl\"ucker
embedding) is generated by the Pl\"ucker coordinates $\Delta_{ij}$,
for $1\le i<j\le n+3$, subject to the Grassmann-Pl\"ucker relations
\[
\Delta_{ik}\,\Delta_{jl}
=\Delta_{ij}\,\Delta_{kl}+\Delta_{il}\,\Delta_{jk}
\qquad (i<j<k<l).
\]
These relations can be viewed as exchange relations in the cluster
algebra $\Acal_n=\CC[\operatorname{Gr}_{2,n+3}]$ of cluster type~$A_n$.
This cluster algebra has $n+3$ coefficient variables
\[
\Delta_{12},\,\Delta_{23},\dots,\Delta_{n+2,n+3},\,\Delta_{1,n+3},
\]
naturally corresponding to the sides of a convex $(n+3)$-gon.
The remaining $\frac{n(n+3)}{2}$ Pl\"ucker coordinates~$\Delta_{ij}$
are the cluster variables; they are naturally labeled by the diagonals
of the $(n+3)$-gon.
Since the exchange relations in $\Acal_n$ can be viewed as Ptolemy
relations between the lambda lengths of sides and diagonals of an
$(n+3)$-gon (=unpunctured disk with $n+3$ marked points on the
boundary), we find ourselves in a situation described in
Remark~\ref{rem:border-coeff-as-lamin}, and can apply the construction
discussed therein.
The case $n=3$ is illustrated in Figure~\ref{fig:G26-lamin}.
\end{example}

\begin{remark}\label{rem:geom-grassmannian}
  We can also give a geometric correspondence between the Grassmannian
  $\Gr_{2,n+3}(\RR)$ and the
  moduli space of decorated ideal $(n+3)$-gons as in
  Example~\ref{ex:unpunctured-disk}, as follows.
  The vector space of $2 \times 2$ symmetric matrices $M = \left(\begin{smallmatrix}
      a&b\\b&c \end{smallmatrix}\right)$ carries a natural quadratic form
\[
-\det M = b^2-ac.  
\]
This quadratic form has
  signature~$(2,1)$ and is invariant under the action of $\operatorname{SL}_2(\RR)$
  by $A \cdot M = A^T M A$.
  Consider the subvariety
\[
H=\left\{\left(\begin{matrix}
      a&b\\b&c \end{matrix}\right)\,:\, ac-b^2 = 1, a>0
\right\}
\]
of positive definite
  symmetric matrices with determinant~$1$.
If we restrict $-\det(\cdot)$ to the
  tangent space to~$H$ at some point, we get a positive definite form.
This gives $H$
  the structure of a Riemannian manifold, which is none
  other than the standard hyperboloid model of the hyperbolic
  plane~$\HH^2$. 
  The ideal boundary of
  $\HH^2$ is the set of positive \emph{null rays}, i.e., the
  set of rank~$1$ positive semidefinite symmetric
  matrices, considered up to scale.  A decorated $(n+3)$-gon is a
  choice of $n+3$ ideal points on the boundary of $\HH^2$ together
  with a choice of horocycle around each ideal point, modulo
  symmetries of~$\HH^2$.  The choice of
  horocycle is the same as a choice of positive vector on the ideal ray, so a
  decorated $(n+3)$-gon is a choice of $n+3$ rank~$1$ positive
  semi-definite matrices, up to simultaneous action by
  $\operatorname{SL}_2(\RR)$ (but not scaling).

A generic point in $\Gr_{2,n+3}(\RR)$ 
can be thought of as a sequence of $n+3$ vectors 
$\biggl(\begin{matrix}
      p_i\\q_i \end{matrix}\biggr)\in\RR^2$ modulo the action 
of~$\operatorname{SL}_2(\RR)$.  
The map
\[
\left(\begin{matrix}
      p_i\\q_i \end{matrix}\right)\,\longmapsto\,
w_i=\begin{pmatrix}
    p_i^2 & p_i q_i \\ p_i q_i & q_i^2
  \end{pmatrix}
\]
identifies such a sequence with a collection $(w_i)$ of rank~1 positive 
semi-definite matrices.
%
  To establish the dictionary between respective coordinates, recall
  \cite{penner-decorated} that the $\lambda$-length between two null
  vectors $w_i, w_j$ in the hyperboloid model can be defined as
  $\lambda(w_i, w_j) = \sqrt{-2\langle w_i, w_j \rangle}$, where
  $\langle\cdot,\cdot\rangle$ denotes the polarization of the quadratic
  form 
(the negated determinant). 
  Elementary computation shows that 
  \[
  \langle w_i , w_j \rangle = -\frac{1}{2}\left(\det
    \begin{pmatrix}
      p_i & p_j \\ q_i & q_j
    \end{pmatrix}\right)^2,
  \]
implying that the lambda lengths $\lambda(w_i, w_j)$ and the 
Pl\"ucker coordinates $\Delta_{ij}$
agree up to sign.  If the vectors in $\RR^2$
  are in one half-space and
  cyclically ordered, then all the $\Delta_{ij}$ are positive and the
  coordinates agree on the nose.  The correct way to get the signs to
  match in general is to consider twisted systems of vectors;
  see \cite[Section~11.1]{fock-goncharov1}.
\end{remark}

\begin{example}[\emph{Grassmannian $\operatorname{Gr}_{3,6}$;
    cf.~\cite{scott}}]
\label{example:G36}

The homogeneous coordinate ring $\Acal=\CC[\operatorname{Gr}_{3,6}]$
has a natural structure of a cluster algebra of (cluster) type~$D_4$,
described in detail by J.~Scott~\cite{scott}.
As a ring, $\Acal$~is generated by the Pl\"ucker coordinates
$\Delta_{ijk}$, for $1\le i<j<k\le 6$.
As a cluster algebra, $\Acal$~has 16 cluster variables,
which include the cluster
\begin{equation}
\label{eq:G36-initial-cluster}
\xx=(x_1,x_2,x_3,x_4)=(\Delta_{245}, \Delta_{256}, \Delta_{125},
\Delta_{235}),
\end{equation}
which we will use as an initial cluster.
The cluster algebra $\Acal$ has 6 coefficient variables
\[
(x_5,\dots,x_{10})=
(\Delta_{123}, \Delta_{234}, \Delta_{345}, \Delta_{456}, \Delta_{156},
\Delta_{126}),
\]
which generate its ground ring.
The exchange relations from the initial~seed,
and the corresponding extended exchange matrix, with rows labeled by the
cluster variables (rows~1--4) and coefficient variables (rows~5--10), are:
\[
\begin{array}{l}
\Delta_{245}\,\Delta_{356} =
  \Delta_{345}\,\Delta_{256} + \Delta_{456}\,\Delta_{235}\\[.1in]
\Delta_{256}\,\Delta_{145} =
  \Delta_{456}\,\Delta_{125} + \Delta_{156}\,\Delta_{245}\\[.1in]
\Delta_{125}\,\Delta_{236} =
  \Delta_{126}\,\Delta_{235} + \Delta_{123}\,\Delta_{256}\\[.1in]
\Delta_{235}\,\Delta_{124} =
  \Delta_{123}\,\Delta_{245} + \Delta_{234}\,\Delta_{125}
\end{array}
\qquad\qquad
\begin{matrix}
\\[-.15in]
\Delta_{245}\\[.05in]
\Delta_{256}\\[.05in]
\Delta_{125}\\[.05in]
\Delta_{235}\\[.08in]
\Delta_{123}\\[.05in]
\Delta_{234}\\[.05in]
\Delta_{345}\\[.05in]
\Delta_{456}\\[.05in]
\Delta_{156}\\[.05in]
\Delta_{126}
\end{matrix}
\begin{bmatrix}
 0 &-1 & 0 & 1 \\[.05in]
 1 & 0 &-1 & 0 \\[.05in]
 0 & 1 & 0 &-1 \\[.05in]
-1 & 0 & 1 & 0 \\[.05in]
\hline\\[-.15in]
 0 & 0 &-1 & 1 \\[.05in]
 0 & 0 & 0 &-1 \\[.05in]
 1 & 0 & 0 & 0 \\[.05in]
-1 & 1 & 0 & 0 \\[.05in]
 0 &-1 & 0 & 0 \\[.05in]
 0 & 0 & 1 & 0
\end{bmatrix}
\]
The topological realization of this cluster algebra,
in which the marked surface $\SM$ is a
once-punctured quadrilateral, is shown in
Figure~\ref {fig:lamin-G36}.
The initial cluster that we chose in~\eqref{eq:G36-initial-cluster}
corresponds to the triangulation shown in Figure~\ref{fig:arcs-G36}.

\begin{figure}[htbp]
\begin{center}
\begin{tabular}{ccc}
\setlength{\unitlength}{3.5pt}
\begin{picture}(33,20)(2,10)
\darkred{\put(3,20){$\Delta_{123}$}}
\thinlines
\multiput(10,10)(20,0){2}{\circle*{1}}
\multiput(10,30)(20,0){2}{\circle*{1}}
\put(20,20){\circle*{1}}
\multiput(10,10)(0,20){2}{\line(1,0){20}}
\multiput(10,10)(20,0){2}{\line(0,1){20}}

\thicklines

\darkred{\put(10,27){\line(1,0){20}}}
\end{picture}
&
\hspace{-.2in}
\setlength{\unitlength}{3.5pt}
\begin{picture}(33,20)(2,10)
\darkred{\put(3,20){$\Delta_{234}$}}
\thinlines

\multiput(10,10)(20,0){2}{\circle*{1}}
\multiput(10,30)(20,0){2}{\circle*{1}}
\put(20,20){\circle*{1}}
\multiput(10,10)(0,20){2}{\line(1,0){20}}
\multiput(10,10)(20,0){2}{\line(0,1){20}}

\thicklines

\darkred{\qbezier(23.5,30)(23.5,22)(23.5,20)}
\darkred{\qbezier(23.5,20)(23.5,16)(20,16)}
\darkred{\qbezier(17,20)(17,16)(20,16)}
\darkred{\qbezier(17,20)(17,23)(20,23)}
\darkred{\qbezier(20,23)(22,23)(22,20)}
\darkred{\qbezier(20,18)(22,18)(22,20)}
\darkred{\qbezier(20,18)(18.7,18)(18.7,20)}
\darkred{\qbezier(20,21.5)(18.7,21.5)(18.7,20)}
\darkred{\qbezier(20,21.5)(21,21.5)(21,20)}
\darkred{\qbezier(20,19)(21,19)(21,20)}

\end{picture}
&
\hspace{-.2in}
\setlength{\unitlength}{3.5pt}
\begin{picture}(33,20)(2,10)
\darkred{\put(3,20){$\Delta_{345}$}}
\thinlines

\multiput(10,10)(20,0){2}{\circle*{1}}
\multiput(10,30)(20,0){2}{\circle*{1}}
\put(20,20){\circle*{1}}
\multiput(10,10)(0,20){2}{\line(1,0){20}}
\multiput(10,10)(20,0){2}{\line(0,1){20}}

\thicklines

\darkred{\qbezier(23.5,10)(23.5,18)(23.5,20)}
\darkred{\qbezier(23.5,20)(23.5,24)(20,24)}
\darkred{\qbezier(17,20)(17,24)(20,24)}
\darkred{\qbezier(17,20)(17,17)(20,17)}
\darkred{\qbezier(20,17)(22,17)(22,20)}
\darkred{\qbezier(20,22)(22,22)(22,20)}
\darkred{\qbezier(20,22)(18.7,22)(18.7,20)}
\darkred{\qbezier(20,18.5)(18.7,18.5)(18.7,20)}
\darkred{\qbezier(20,18.5)(21,18.5)(21,20)}
\darkred{\qbezier(20,21)(21,21)(21,20)}

\end{picture}
\\[.3in] 

\setlength{\unitlength}{3.5pt}
\begin{picture}(33,20)(2,10)
\darkred{\put(3,20){$\Delta_{456}$}}
\thinlines

\multiput(10,10)(20,0){2}{\circle*{1}}
\multiput(10,30)(20,0){2}{\circle*{1}}
\put(20,20){\circle*{1}}
\multiput(10,10)(0,20){2}{\line(1,0){20}}
\multiput(10,10)(20,0){2}{\line(0,1){20}}

\thicklines

\darkred{\put(10,13){\line(1,0){20}}}

\end{picture}
&
\hspace{-.2in}
\setlength{\unitlength}{3.5pt}
\begin{picture}(33,20)(2,10)
\darkred{\put(3,20){$\Delta_{156}$}}
\thinlines

\multiput(10,10)(20,0){2}{\circle*{1}}
\multiput(10,30)(20,0){2}{\circle*{1}}
\put(20,20){\circle*{1}}
\multiput(10,10)(0,20){2}{\line(1,0){20}}
\multiput(10,10)(20,0){2}{\line(0,1){20}}

\thicklines

\darkred{\qbezier(16.5,10)(16.5,18)(16.5,20)}
\darkred{\qbezier(16.5,20)(16.5,24)(20,24)}
\darkred{\qbezier(23,20)(23,24)(20,24)}
\darkred{\qbezier(23,20)(23,17)(20,17)}
\darkred{\qbezier(20,17)(18,17)(18,20)}
\darkred{\qbezier(20,22)(18,22)(18,20)}
\darkred{\qbezier(20,22)(21.3,22)(21.3,20)}
\darkred{\qbezier(20,18.5)(21.3,18.5)(21.3,20)}
\darkred{\qbezier(20,18.5)(19,18.5)(19,20)}
\darkred{\qbezier(20,21)(19,21)(19,20)}

\end{picture}
&
\hspace{-.2in}
\setlength{\unitlength}{3.5pt}
\begin{picture}(33,20)(2,10)
\darkred{\put(3,20){$\Delta_{126}$}}
\thinlines

\multiput(10,10)(20,0){2}{\circle*{1}}
\multiput(10,30)(20,0){2}{\circle*{1}}
\put(20,20){\circle*{1}}
\multiput(10,10)(0,20){2}{\line(1,0){20}}
\multiput(10,10)(20,0){2}{\line(0,1){20}}

\thicklines

\darkred{\qbezier(16.5,30)(16.5,22)(16.5,20)}
\darkred{\qbezier(16.5,20)(16.5,16)(20,16)}
\darkred{\qbezier(23,20)(23,16)(20,16)}
\darkred{\qbezier(23,20)(23,23)(20,23)}
\darkred{\qbezier(20,23)(18,23)(18,20)}
\darkred{\qbezier(20,18)(18,18)(18,20)}
\darkred{\qbezier(20,18)(21.3,18)(21.3,20)}
\darkred{\qbezier(20,21.5)(21.3,21.5)(21.3,20)}
\darkred{\qbezier(20,21.5)(19,21.5)(19,20)}
\darkred{\qbezier(20,19)(19,19)(19,20)}


\end{picture}
\end{tabular}
\end{center}
\caption{Representing the cluster structure on the Grassmannian
  $\operatorname{Gr}_{3,6}$ by a multi-lamination $\LL$
on a once-punctured quadrilateral. Each of the 6~laminations
in~$\LL$ consists of a single curve, and corresponds to a
particular coefficient variable (a Pl\"ucker coordinate).
}
\label{fig:lamin-G36}
\end{figure}
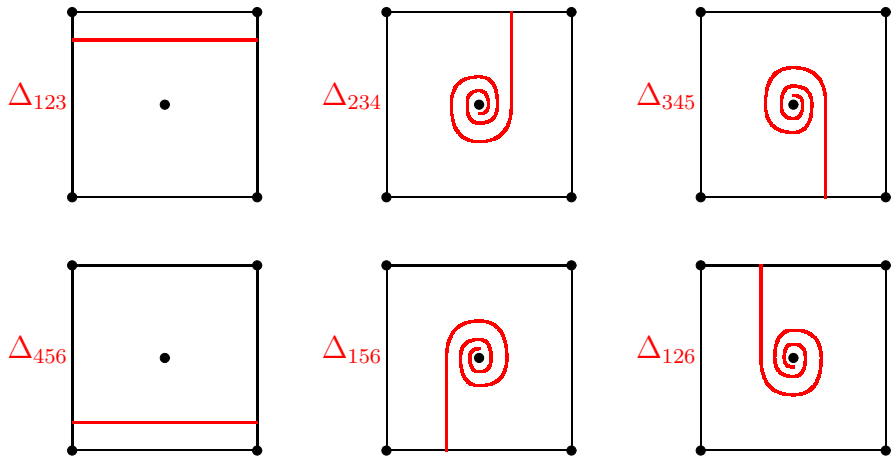

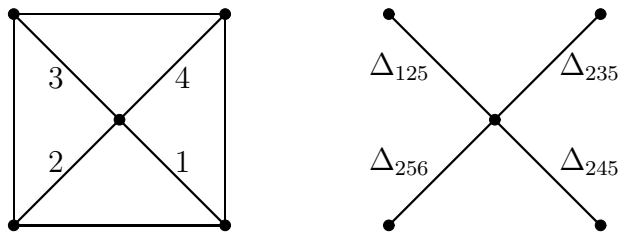
\begin{figure}[htbp]
\begin{center}
\begin{tabular}{ccc}
\setlength{\unitlength}{4pt}
\begin{picture}(33,22)(2,10)
\thinlines

\multiput(10,10)(20,0){2}{\circle*{1}}
\multiput(10,30)(20,0){2}{\circle*{1}}
\put(20,20){\circle*{1}}
\multiput(10,10)(0,20){2}{\line(1,0){20}}
\multiput(10,10)(20,0){2}{\line(0,1){20}}

\put(26,16){\makebox(0,0){$1$}}
\put(14,16){\makebox(0,0){$2$}}
\put(14,24){\makebox(0,0){$3$}}
\put(26,24){\makebox(0,0){$4$}}
\thicklines
\put(10,10){\line(1,1){20}}
\put(10,30){\line(1,-1){20}}

\end{picture}
&
\begin{picture}(33,22)(2,10)
\thinlines
\multiput(10,10)(20,0){2}{\circle*{1}}
\multiput(10,30)(20,0){2}{\circle*{1}}
\put(20,20){\circle*{1}}

\put(29,16){\makebox(0,0){$\Delta_{245}$}}
\put(11,16){\makebox(0,0){$\Delta_{256}$}}
\put(11,25){\makebox(0,0){$\Delta_{125}$}}
\put(29,25){\makebox(0,0){$\Delta_{235}$}}
\thicklines
\put(10,10){\line(1,1){20}}
\put(10,30){\line(1,-1){20}}

\end{picture}
\end{tabular}
\end{center}
\caption{The triangulation representing the initial
cluster \eqref{eq:G36-initial-cluster}
for the Grassmannian  $\operatorname{Gr}_{3,6}$.
The remaining 12 tagged arcs (not shown in the picture) correspond to
the rest of the cluster variables.
}
\label{fig:arcs-G36}
\end{figure}
\end{example}

\begin{example}[\emph{The ring $\CC[\operatorname{Mat}_{3,3}]$}]
\label{example:mat33}
The ring of polynomials in 9 variables $z_{ij}$ ($i,j\in\{1,2,3\}$)
viewed as matrix entries of a $3\times 3$ matrix
\[
z=\begin{bmatrix}
z_{11} & z_{12} & z_{13} \\[.05in]
z_{21} & z_{22} & z_{23} \\[.05in]
z_{31} & z_{32} & z_{33}
\end{bmatrix}
\in\operatorname{Mat}_{3,3}\cong\CC^9
\]
carries a natural cluster algebra structure of cluster type~$D_4$;
see, e.g.,
\cite{dbc, cdm, skandera}.
This cluster structure is very similar to the one discussed in
Example~\ref{example:G36}, and is obtained as follows.
The map $\operatorname{Mat}_{3,3}\to \operatorname{Gr}_{3,6}$ defined
by 
\[
z\mapsto
\operatorname{rowspan}\begin{pmatrix}
z_{11} & z_{12} & z_{13} &0&0&1\\[.05in]
z_{21} & z_{22} & z_{23} &0&-1&0\\[.05in]
z_{31} & z_{32} & z_{33} &1&0&0
\end{pmatrix}
\]
induces a ring homomorphism
\[
\varphi: \CC[\operatorname{Gr}_{3,6}]\to\CC[z_{11},\dots,z_{33}].
\]
For example, $\varphi(\Delta_{145})=z_{11}$.
More generally, $\varphi$ maps the 19 Pl\"ucker coordinates
$\Delta_{ijk}$---all but~$\Delta_{456}$---into the 19 minors
of~$z$. 
In fact, $\varphi$ sends each of the 16 cluster variables in
$\CC[\operatorname{Gr}_{3,6}]$ into a cluster variable in
$\CC[z_{11},\dots,z_{33}]$, and sends all but one coefficient
variables in $\CC[\operatorname{Gr}_{3,6}]$ into coefficient
variables in $\CC[z_{11},\dots,z_{33}]$,
the sole exception being $\varphi(\Delta_{456})=1$.
The exchange relations in the cluster algebra
$\CC[z_{11},\dots,z_{33}]$ are obtained from those in
$\CC[\operatorname{Gr}_{3,6}]$ by applying~$\varphi$.
As a result, we get a cluster structure that can be realized by a
multi-lamination on a once-punctured quadrilateral that consists of
5~laminations: the ones shown in Figure~\ref{fig:lamin-G36} with the
exception of the leftmost lamination in the bottom row.
\end{example}

\pagebreak[3]

\begin{example}[\emph{The special linear group $\operatorname{SL}_3$;
cf.~\cite{dbc}}]
\label{example:sl3}
This cluster algebra (also of cluster type~$D_4$) is obtained from the
one in Example~\ref{example:mat33} by further specializing the
coefficients, this time setting $\det(z)=1$.
Equivalently, we send both coefficient variables
$\Delta_{123}$ and $\Delta_{456}$ in Example~\ref{example:G36} to~1.
This translates into removing the leftmost laminations in each
of the two rows in Figure~\ref{fig:lamin-G36}.
\end{example}

\begin{example}[\emph{Affine base space $\operatorname{SL}_4/N$}]
\label{example:SL4/N} 

This example---with $\operatorname{SL}_4$ replaced by an arbitrary
complex semisimple Lie group---was one of the main examples that
motivated the introduction of cluster algebras~\cite{ca1}.
Let $N$ denote the subgroup of uni\-potent
  upper-triangular matrices in $\operatorname{SL}_4(\mathbb{C})$.
The group~$N$ acts on $\operatorname{SL}_4$ by left multiplication;
let
\[
\Acal=\CC[\operatorname{SL}_4/N]
\subset \mathbb{C}[z_{11},\dots,z_{44}]/\langle\det(z)-1\rangle
\]
be the ring of $N$-invariant polynomials in the matrix entries
$z_{ij}$.
The ring $\Acal$~is
generated by the \emph{flag minors} 
\[
\Delta_I: z=(z_{ij})\mapsto \det(z_{ij}|i\in I, j\le|I|),
\]
for $I\subsetneq\{1,2,3,4\}$, $I\neq\emptyset$.
These flag minors satisfy well-known Pl\"ucker-type relations.

The ring $\Acal$ carries a natural cluster structure of type~$A_3$
(see, e.g., \cite[Section~2.6]{ca3}),
with 6 coefficient variables
\begin{equation}
\label{eq:SL4/N-coeff-var}
\Delta_1, \Delta_{12}, \Delta_{123}, \Delta_4, \Delta_{34},
\Delta_{234},
\end{equation}
and 9 cluster variables
\begin{equation}
\label{eq:SL4/N-cluster-var}
\Delta_2, \Delta_3, \Delta_{13}, \Delta_{14}, \Delta_{23},
\Delta_{24}, \Delta_{124}, \Delta_{134},
\Omega=-\Delta_1\Delta_{234}+\Delta_2\Delta_{134}\,.
\end{equation}
Let the initial cluster be $\xx=(\Delta_2, \Delta_3, \Delta_{23})$.
Then the exchange relations from the initial~seed are:
\begin{align*}
\Delta_2\Delta_{13}&=\Delta_{12}\Delta_3+\Delta_1\Delta_{23}\\
\Delta_3\Delta_{24}&=\Delta_{4}\Delta_{23}+\Delta_{34}\Delta_{2}\\
\Delta_{23}\,\Omega\,\,&=\Delta_{123}\Delta_{34}\Delta_2+\Delta_{12}\Delta_{234}\Delta_3
\end{align*}

This algebra
can be described using the multi-lamination $\LL$ of a hexagon shown in
Figure~\ref{fig:sl-4-lamin-dictionary}.

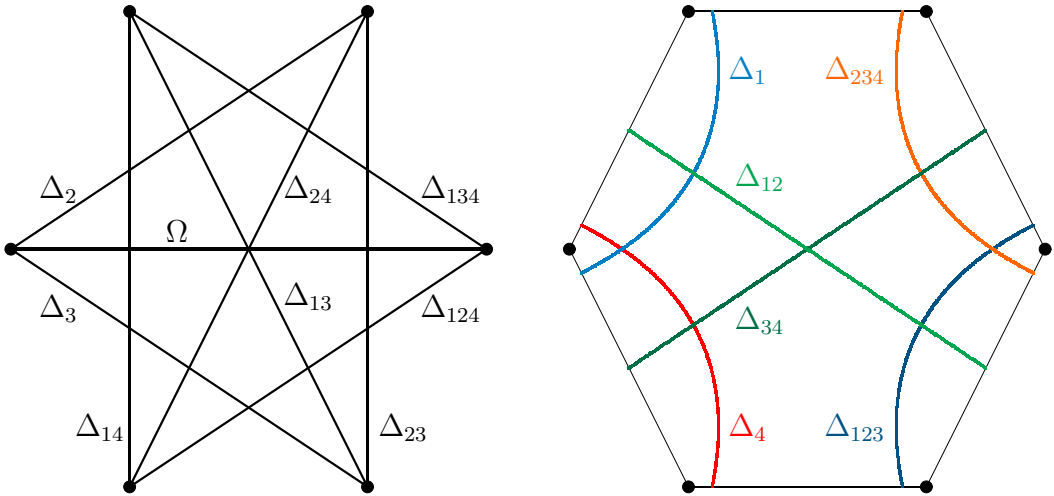
\begin{figure}[htbp]
\begin{center}
\setlength{\unitlength}{4.5pt}
\begin{picture}(40,42)(0,0)
\thinlines

  \put(0,20){\circle*{1}}
  \put(10,0){\circle*{1}}
  \put(10,40){\circle*{1}}
  \put(30,0){\circle*{1}}
  \put(30,40){\circle*{1}}
  \put(40,20){\circle*{1}}

\thicklines
  \put(0,20){\line(3,2){30}}
\put(4,25){\makebox(0,0){$\Delta_2$}}
  \put(0,20){\line(3,-2){30}}
\put(4,15){\makebox(0,0){$\Delta_3$}}
  \put(30,0){\line(0,1){40}}
\put(33,5){\makebox(0,0){$\Delta_{23}$}}

  \put(40,20){\line(-3,2){30}}
\put(37,25){\makebox(0,0){$\Delta_{134}$}}
  \put(40,20){\line(-3,-2){30}}
\put(37,15){\makebox(0,0){$\Delta_{124}$}}
  \put(10,0){\line(0,1){40}}
\put(7.5,5){\makebox(0,0){$\Delta_{14}$}}

  \put(10,0){\line(1,2){20}}
\put(25,25){\makebox(0,0){$\Delta_{24}$}}
  \put(10,40){\line(1,-2){20}}
\put(25,16){\makebox(0,0){$\Delta_{13}$}}
  \put(0,20){\line(1,0){40}}
\put(14,21.5){\makebox(0,0){$\Omega$}}

\end{picture}
\qquad
\setlength{\unitlength}{4.5pt}
\begin{picture}(40,42)(0,0)
\thinlines
  \put(10,0){\line(-1,2){10}}
  \put(40,20){\line(-1,2){10}}
  \put(10,0){\line(1,0){20}}
  \put(10,40){\line(1,0){20}}
  \put(0,20){\line(1,2){10}}
  \put(30,0){\line(1,2){10}}

  \put(0,20){\circle*{1}}
  \put(10,0){\circle*{1}}
  \put(10,40){\circle*{1}}
  \put(30,0){\circle*{1}}
  \put(30,40){\circle*{1}}
  \put(40,20){\circle*{1}}

\thicklines
\darkred{\qbezier(12,0)(15,15)(1,22)
\put(15,5){\makebox(0,0){$\Delta_4$}}
}
\lightblue{\qbezier(12,40)(15,25)(1,18)
\put(15,35){\makebox(0,0){$\Delta_1$}}
}
\darkblue{\qbezier(28,0)(25,15)(39,22)
\put(24,5){\makebox(0,0){$\Delta_{123}$}}
}
\lightred{\qbezier(28,40)(25,25)(39,18)
\put(24,35){\makebox(0,0){$\Delta_{234}$}}
}
\darkgreen{\qbezier(5,10)(20,20)(35,30)
\put(16,14){\makebox(0,0){$\Delta_{34}$}}
}
\lightgreen{\qbezier(5,30)(20,20)(35,10)
\put(16,26){\makebox(0,0){$\Delta_{12}$}}
}

\end{picture}
\end{center}
\caption{Representing the cluster structure on
the affine base space
$\operatorname{SL}_4/N$
by a multi-lamination of a hexagon.
The 9~cluster variables (cf.~\eqref{eq:SL4/N-cluster-var})
correspond to the 9 arcs (diagonals) as shown
on the left. The 6 coefficient variables
(cf.~\eqref{eq:SL4/N-coeff-var}) correspond to the
6~single-curve laminations as shown on the right.
}
\label{fig:sl-4-lamin-dictionary}
\end{figure}

\end{example}

\begin{remark}
The cluster algebra of Example~\ref{example:SL4/N} illustrates a
relatively rare phenomenon where the same algebra 
can be given two \emph{different} topological realizations---because
its exchange matrices can be obtained from two topologically different \linebreak[3]
surfaces, in this case an unpunctured hexagon and a once-punctured
triangle. (Informally speaking, types $A_3$ and~$D_3$ coincide.) 
Consequently, one can choose to represent the cluster structure in
$\CC[\operatorname{SL}_4/N]$ by a multi-lamination in a once-punctured
triangle.
Cf.\ \cite[Examples~4.3 and~4.5]{cats1} and the discussion at the
beginning of \cite[Section~14]{cats1}. 
\end{remark}

\begin{example}[\emph{Unipotent subgroup $N_-\subset\operatorname{SL}_4$}]
The coordinate ring of the subgroup $N_-$ of unipotent
lower-triangular $4\times 4$ matrices carries a cluster algebra
structure that can 
be obtained by specialization of coefficients from the cluster
structure in $\CC[\operatorname{SL}_4/N]$ discussed in
Example~\ref{example:SL4/N}: we set
$\Delta_1=\Delta_{12}=\Delta_{123}=1$. 
(Cf.\ \cite{ca3} and \cite[Sections~4.2.6--4.3]{gls}.) 
This corresponds to removing the three laminations with these labels
from Figure~\ref{fig:sl-4-lamin-dictionary}. 
\end{example}

\begin{example}[\emph{Affine base space $\operatorname{SL}_5/N$
and unipotent subgroup $N_-\subset\operatorname{SL}_5$}]
In these two cases, one gets cluster algebras of type~$D_6$, 
which can be represented by a collection of~$8$ 
(respectively,~$4$) single-curve laminations in a once-punctured hexagon. 
We omit the details for the sake of brevity. 
\end{example}

Further examples of similar nature include various coordinate rings
discussed in \cite[Examples~6.7, 6.9, 6.10]{cats1}, among others.

\chapter{Principal and universal coefficients}
\label{sec:prin-univ-coeff}

In this chapter, we work out two particular cases of our main
construction, yielding two distinguished choices of coefficients
introduced in~\cite{ca4}.

\begin{definition}[\emph{Principal coefficients \cite[Definition~3.1,
    Remark~3.2]{ca4}}]
\label{def:principal}
Let $(\Sigma_t)=(\xx(t),\tilde B(t))$ be an exchange pattern of
geometric type.
We say that this pattern (or the corresponding cluster
algebra~$\Acal$) has \emph{principal coefficients} with respect to an initial
vertex $t_0$ if the extended exchange matrix $\tilde B(t_0)$
has order $2n\times n$ (as always, $n$~is the rank of the pattern),
and the bottom $n\times n$ part of $\tilde B(t_0)$ is the identity
matrix. We denote $\Acal=\Acal_\bullet(B(t_0))$, where $B(t_0)$ is the
initial $n\times n$ exchange matrix.

To rephrase, let $\xx(t_0)=(x_1,\dots,x_n)$ be the initial cluster, and
let $B(t_0)=(b_{ij})$.
Then the
algebra $\Acal=\Acal_\bullet(B(t_0))$ is defined as follows.
The ground semifield of $\Acal$ is $\Trop(q_1,\dots,q_n)$,
and the exchange relations \eqref{eq:exchange-rel-E} out of the
initial seed $\Sigma_{t_0}$ are
\[
x_k x_k' =
q_k \prod_{\substack{1\le i\le n\\ b_{ik}>0}} x_i^{b_{ik}}
  + \prod_{\substack{1\le i\le n\\ b_{ik}<0}} x_i^{-b_{ik}}.
\]
\end{definition}

To construct a topological model for a cluster algebra with principal
coefficients, we will need the following notion.

\begin{figure}[htbp]
\begin{center}
\setlength{\unitlength}{2pt}
\begin{picture}(40,35)(20,-3)
\thinlines
\put(0,25){\circle{20}}
\put(0,15){\circle*{2}}
\put(0,35){\circle*{2}}
\put(-1,19){\makebox(0,0){$a$}}
\qbezier(0,15)(0,0)(20,0)
\qbezier(40,15)(40,0)(20,0)

%
\put(4,0){\makebox(0,0){$\gamma$}}

\put(40,25){\circle{20}}
\put(30,25){\circle*{2}}
\put(50,25){\circle*{2}}
\put(40,15){\circle*{2}}
\put(40,19){\makebox(0,0){$b$}}

\thicklines
\darkred{\qbezier(2,15.2)(2,3)(16,2)}
\darkred{\qbezier(16,2)(17.4,1.9)(20,0)}
\darkred{\qbezier(20,0)(22.2,-1.9)(24,-2)}
\darkred{\qbezier(42,15.2)(42,-1)(24,-2)}
\darkred{\put(2,15.2){\circle*{2}}}
\darkred{\put(6,14.5){\makebox(0,0){$a'$}}}
\darkred{\put(42,15.2){\circle*{2}}}
\darkred{\put(46,14){\makebox(0,0){$b'$}}}
\darkred{\put(42,0){\makebox(0,0){$L_\gamma$}}}
\end{picture}
\qquad\qquad\qquad
\setlength{\unitlength}{4pt}
\begin{picture}(20,17.5)(30,6)
\thinlines
  \put(20,20){\circle*{1}}
  \put(60,20){\circle*{1}}
\qbezier(20,20)(20,5)(40,5)
\qbezier(60,20)(60,5)(40,5)
\put(59.5,15.5){\makebox(0,0){$\notch$}}

  \put(19,19){\makebox(0,0){$a$}}
  \put(22,8){\makebox(0,0){$\gamma$}}
\thicklines

\darkred{\qbezier(40,9)(51,9)(55,20)}
\darkred{\qbezier(40,9)(29,9)(25,20)}
\darkred{\qbezier(55,20)(56.5,24)(60,24)}
\darkred{\qbezier(25,20)(23.5,24)(20,24)}
\darkred{\qbezier(63,20)(63,24)(60,24)}
\darkred{\qbezier(17,20)(17,24)(20,24)}
\darkred{\qbezier(63,20)(63,17)(60,17)}
\darkred{\qbezier(17,20)(17,17)(20,17)}
\darkred{\qbezier(20,17)(22,17)(22,20)}
\darkred{\qbezier(60,17)(58,17)(58,20)}
\darkred{\qbezier(20,21.5)(22,21.5)(22,20)}
\darkred{\qbezier(60,21.5)(58,21.5)(58,20)}
\darkred{\qbezier(20,21.5)(18.7,21.5)(18.7,20)}
\darkred{\qbezier(60,21.5)(61.3,21.5)(61.3,20)}

\darkred{\put(40,11){\makebox(0,0){$L_\gamma$}}}

\end{picture}

\end{center}
\caption{Elementary laminations associated with tagged arcs}
\label{fig:elem-lamin-arc}
\end{figure}
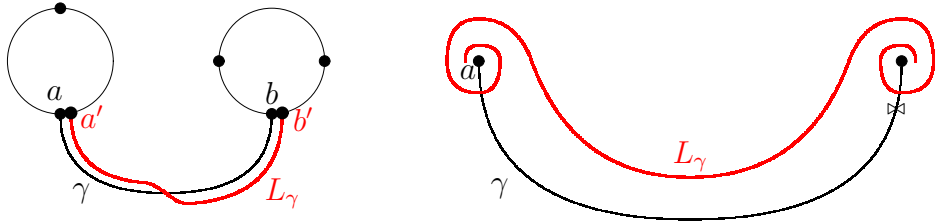

\begin{definition}[\emph{Elementary lamination associated with a
    tagged arc}]
\label{def:elem-lamin-arc}
Let $\gamma\in\APSM$ be a tagged arc in~$\SM$.
We denote by $L_\gamma$ a lamination consisting of a single curve
(also denoted by $L_\gamma$, by an abuse of notation) defined as
follows (up to isotopy).
The curve $L_\gamma$ runs along $\gamma$ in a small neighborhood
thereof; to complete its description, we only need to specify what
happens near the ends.
Assume that $\gamma$ has an endpoint~$a$ on the boundary of~$\Surf$,
more specifically on a circular component~$C$.
Then  $L_\gamma$ begins at a point $a'\in C$ located near~$a$ in
the counterclockwise direction, and proceeds along~$\gamma$ as shown in
Figure~\ref{fig:elem-lamin-arc} on the left.

If  $\gamma$ has an endpoint~$a\in\barM$ (a puncture), then $L_\gamma$
spirals into~$a$: counterclockwise if $\gamma$ is tagged plain at~$a$,
and clockwise if it is notched. See Figure~\ref{fig:elem-lamin-arc} on
the right.
\end{definition}

The following statement is straightforward to check.

\begin{proposition}
Let $T$ be a tagged triangulation with a signed adjacency
matrix~$B(T)$.
Then $\Acal_\bullet(B(T))\cong\Acal(\Surf,\Mark,\LL_T)$
(cf.\ Definition~\ref{def:ca-multi-lamin}), where
$\LL_T=(L_\gamma)_{\gamma\in T}$ is the multi-lamination consisting of
elementary laminations associated with the tagged arcs in~$T$.
\end{proposition}

\begin{example}
Let $T=(\gamma_{ap},\gamma_{bp})$ be the triangulation of a
once-punctured digon shown in Figure~\ref{fig:tagged-triang-digon}.
Then $\LL_T=(L_4,L_3)$, in the notation of
Figure~\ref{fig:lamin-digon}.
\end{example}

We next turn to the discussion of cluster algebras with universal
coefficients, which were constructed in \cite[Section~12]{ca4}
for any finite Cartan-Killing (cluster) type.
More explicit versions of this construction for the types $A$ and~$D$ 
have been given in the unpublished work~\cite{fz-classical};
the type~$A$ case was reproduced in \cite[Proposition~7.2]{yang-zel}. 
Here we provide a restatement of the aforementioned results for the 
types $A$ and~$D$.

\begin{proposition}
Let $\SM$ be a marked surface of finite cluster type, so that the set
of tagged arcs $\ATSM$ is finite.
The corresponding cluster algebra with universal coefficients can be
realized as $\Acal(\Surf,\Mark,\LL_\circ)$, where the
multi-lamination~$\LL_\circ$ consists of all elementary laminations
(see Definition~\ref{def:elem-lamin-arc}) 
associated with the tagged arcs $\gamma\in \ATSM$.
\end{proposition}

The straightforward proof of this proposition is omitted.

We note that among simply-laced finite cluster types, only types $A$ and~$D$, 
and their direct products, have topological realizations.
(see \cite[Examples 6.6--6.7, Remark~13.5]{cats1}).

\begin{example}
Let $\SM$ be the once-punctured digon, so that the cluster type is
$A_1\times A_1$.
Then $\LL_\circ=(L_1,L_3,L_4,L_6)$, in the notation of
Figure~\ref{fig:lamin-digon}.
\end{example}

\begin{example}
Let $\SM$ be an unpunctured $(n+3)$-gon (type~$A_n$).
Then $\LL_\circ$ consists of all elementary (i.e., single-curve)
laminations in~$\SM$.
These are the $\frac{n(n+3)}{2}$ curves
connecting non-adjacent midpoints of the sides of the polygon.
\end{example}

\begin{example}
Let $\SM$ be a once-punctured $n$-gon (type~$D_n$).
Again, $\LL_\circ$ consists of all elementary
laminations in~$\SM$.
They include: $(n-3)n$ curves connecting non-adjacent midpoints as in the
previous example (the
number is doubled since we can go on either side of the puncture);
plus $n$ curves connecting adjacent midpoints, going
around the puncture;
plus $2n$ curves starting at one of the $n$ midpoints and spiraling into the
puncture (either clockwise or counterclockwise);
for the grand total of~$n^2$.
\end{example}

\appendix

\chapter{Tropical degeneration and relative lambda lengths}
\label{app:relative-lambda-lengths}

Here we sketch how our main construction of exchange patterns formed
by generalized lambda lengths in ``laminated Teichm\"uller spaces''
can be obtained by \emph{tropical degeneration} from a somewhat more
traditional construction in which laminations are replaced by
hyperbolic structures of ordinary kind.

First, a simple observation concerning (commutative) semifields;
see Definition~\ref{def:n-exch-pattern}.

\begin{lemma}
\label{lem:semifield-product}
A direct product of semifields is a semifield,
under component-wise operations.
\end{lemma}

The tropical semifield of Definition~\ref{def:tropical} is one such
example:
\begin{equation}
\label{eq:trop=product}
\Trop(q_1,q_2,\dots) \cong \Trop(q_1) \times \Trop(q_2) \times \cdots .
\end{equation}

We next describe, somewhat informally, the general procedure of
tropical degeneration.
Fix a real number~$k$, and define the binary operation~$\oplusk$
by
\begin{equation}
\label{eq:oplusk}
y\oplusk z=(y^k + z^k)^{1/k} ,
\end{equation}
where $y$ and $z$ are positive reals, or
positive-valued functions on some fixed set.
It is easy to see that the binary operation $\oplusk$ is commutative,
associative, and distributive with respect to the ordinary
multiplication.
This makes $(\PP,\oplusk,\cdot)$ a semifield,
where $\PP$ is  any set of numbers or functions that is closed under
$\oplusk$ and under multiplication.
The simplest instance of this construction produces the
semifield
$(\RR_{>0},\oplusk,\cdot)$ of positive real numbers, with ordinary
multiplication and with addition~$\oplusk$ defined
by~\eqref{eq:oplusk}.
In the ``tropical limit,'' as $k\to-\infty$, we get the semifield
$(\RR_{>0},\oplus,\cdot)$ with the addition~$\oplus$ given by
\begin{equation}
y\oplus z = \lim_{k\to-\infty} y\oplusk z = \min(y,z).
\end{equation}
This is a close relative of the tropical
addition~\eqref{eq:tropical addition}, in the one-dimensional version,
with the set~$I$ in~\eqref{eq:PP=Trop} consisting of a single element.
To obtain a multi-dimensional version, we apply the direct product
construction of Lemma~\ref{lem:semifield-product},
either after taking the limit (as in~\eqref{eq:trop=product}), or before.
Let us discuss the latter approach in more concrete detail.
Start with a collection of real parameters $\kk=(k_i)_{i\in I}$, and
consider the semifield of tuples $\yy=(y_i)_{i\in I}$
with the ordinary (pointwise) multiplication and with the addition
$\oplus_\kk$ defined by
\[
(\yy \,\oplus_\kk\, \zz)_i
= y_i\,\oplus_{k_i}\, z_i
= (y_i^{k_i} + z_i^{k_i})^{1/k_i}\,.
\]
Taking the limits $k_i\to-\infty$ for all~$i$, we obtain a close
relative of the tropical semifield~\eqref{eq:PP=Trop},
with the addition defined by
\[
(\yy \oplus \zz)_i
= \min(y_i,z_i)\,.
\]


The geometric counterpart of the tropical degeneration procedure
consists in obtaining a lamination, viewed as a point on the Thurston
boundary~\cite{thurston-travaux, thurston-minimal-stretch} of the
appropriate Teichm\"uller space, as a limit of
ordinary hyperbolic structures.
For simplicity, we discuss the case of surfaces with no marked
points in the interior. The general case can be handled in exactly the
same way as before, by lifting arcs to the opened surface.
We still want to get the most general coefficients,
so we cannot use the construction of
Theorem~\ref{th:cluster-lambda} (coefficients coming from boundary
segments); rather, we shall aim at obtaining the construction of
Theorem~\ref{th:cluster-lamin} (coefficients coming from a
multi-lamination) in the restricted generality of surfaces with no
punctures.

Let $\sigma\in\dTeich\SM$ be a decorated hyperbolic structure.
(In the presence of punctures, we would need to consider
$\sigma\in\TTSM$.)
The lambda lengths $\lambda_\sigma(\gamma)$ satisfy generalized
exchange/Ptolemy relations \eqref{eq:ptolemy-lambda}:
\begin{equation}
\label{eq:ptolemy-lambda-sigma}
\lambda_\sigma(\eta) \lambda_\sigma(\theta)
= \lambda_\sigma(\alpha)\lambda_\sigma(\gamma)
  + \lambda_\sigma(\beta)\lambda_\sigma(\delta).
\end{equation}
Choose $k\in\RR$, and set
\[
\lambda_{\sigma,k}(\gamma)=(\lambda_\sigma(\gamma) )^{1/k}.
\]
These quantities satisfy the $\oplusk$-version
of~\eqref{eq:ptolemy-lambda-sigma}:
\begin{equation}
\label{eq:ptolemy-lambda-k}
\lambda_\sigma(\eta) \lambda_\sigma(\theta)
= \lambda_\sigma(\alpha)\lambda_\sigma(\gamma)
  \oplusk \lambda_\sigma(\beta)\lambda_\sigma(\delta).
\end{equation}
Recalling Theorem~\ref{th:decorated-coord},
pick a triangulation~$T$ and a collection of positive reals
$(c(\gamma))_{\gamma \in T\cup\BSM}$.
Then let $\sigma$ and~$k$ vary so that $k\to-\infty$ while the
coordinates $(\lambda_{\sigma,k}(\gamma))$ remain fixed:
\[
\lambda_{\sigma,k}(\gamma)=c(\gamma), \quad \text{for $\gamma \in
  T\cup\BSM$}.
\]
In other words, we set $\lambda_\sigma(\gamma)=c(\gamma)^k$
and let $k\to-\infty$.
As a result, if all of the $c(\gamma) < 1$, then the
$\lambda$-coordinates will go to $+\infty$ and
$\sigma$ goes to a point on the Thurston boundary
of $\dTeich\SM$ which can be identified with a real measured
lamination~$L$ whose tropical lambda lengths match the values we
picked:
\[
c_L(\gamma)=c(\gamma) , \quad \text{for $\gamma \in
  T\cup\BSM$}.
\]
Meanwhile, the exchange relations \eqref{eq:ptolemy-lambda-k} degenerate
into
\[
c_L(\eta) c_L(\theta)
= c_L(\alpha)c_L(\gamma) \oplus c_L(\beta)c_L(\delta),
\]
where $\oplus$ denotes the minimum (the tropical addition).
This fairly standard argument shows how the tropical exchange
relations for the tropical lambda
lengths associated with a single lamination can be obtained by
tropical degeneration from the ordinary exchange relations for the lambda
lengths with respect to a hyperbolic structure.

(If some of the $c(\gamma) > 1$, then some $\lambda$-coordinates will
go to~$0$.  This can still be interpreted as $\sigma$ going to a point
on a generalized ``Thurston boundary'' of Teichm\"uller space, but is
more confusing geometrically.)

To extend this construction to the case of multiple laminations,
we can use the direct product construction of
Lemma~\ref{lem:semifield-product}, yielding a solution of exchange
relations in the tropical semifield with several parameters, as in
Lemmas~\ref{lem:ptolemy-x-trop} and~\ref{lem:exch-digon-trop}.
This solution can then be used to renormalize the ordinary lambda
lengths as in~\eqref{eq:x-lamin},
delivering the requisite normalized exchange pattern by virtue of
Proposition~\ref{pr:rescale-normalize}, as in the proof of
Theorem~\ref{th:cluster-lamin}.

As we have seen, the process of creating a normalized exchange pattern
formed by renormalized lambda lengths consists of three main stages:
\begin{enumerate}
\item
producing a tropical solution of exchange relations from an ordinary
one by means of ``tropical degeneration'';
\item
applying a direct product construction to get a solution depending on
several reference geometries; and
\item
dividing the original non-normalized solution by the one produced at
stage (2) to obtain a normalized pattern.
\end{enumerate}
Alternatively, stage~(1) (tropical
degeneration) can be executed after stage~(2) or even after stage~(3).
To conclude this appendix, we show how to perform stage~(3)
in the
restricted case of a single reference geometry (where stage~(2) is trivial).
To compensate for this restriction, we will consider a general case of
a bordered surface \emph{with} punctures.
This construction bears some similarity to the work of V.~Fock and A.~Goncharov
on cluster varieties, and in particular to their notion of
\emph{symplectic double} \cite[Section~2.2]{fg-quant-dilog}.

Fix a reference geometry, i.e.,
a decorated hyperbolic structure $\varrho\in\TTSM$.
Now, for any other $\sigma\in\TTSM$ satisfying the boundary conditions
\begin{align}
\label{eq:boundary-lambda-match-1}
&\text{$\lambda_\sigma(\beta)=\lambda_\varrho(\beta)$ for all
  $\beta\in\BSM$}, \\
\label{eq:boundary-lambda-match-2}
&\text{$\lambda_\sigma(p)=\lambda_\varrho(p)$ for all $p \in \barM$,}
\end{align}
 the \emph{relative lambda length}
$\lambda_{\sigma/\varrho}(\gamma)$
of an arc $\gamma\in\APSM$ is defined by
\begin{equation}
  \lambda_{\sigma/\varrho}(\gamma)
=
  \lambda_\sigma(\bargamma)/\lambda_\varrho(\bargamma),
\end{equation}
where $\bargamma$ is an arbitrary lift of~$\gamma$ to the opened
surface $\barSM$.
Note that this ratio does not depend on the choice of a
lift~$\bargamma$ since the numerator and the denominator rescale by
the same factor under the twists~$\psi_p$ (see~\eqref{eq:signed-lambda-twist}).

By Theorem~\ref{th:cluster-lambda-opened},
the lambda lengths $\lambda_\varrho(\gamma)$
(respectively, $\lambda_\sigma(\gamma)$)
form a non-normalized exchange pattern with coefficients in the
multiplicative group generated by
the boundary
parameters~\eqref{eq:boundary-lambda-match-1}--\eqref{eq:boundary-lambda-match-2}.
Consequently, by Proposition~\ref{pr:rescale-normalize},
the relative lambda lengths
$\lambda_{\sigma/\varrho}(\gamma)$ form a normalized exchange
pattern over the semifield generated by the lambda lengths for~$\varrho$,
with ordinary addition and multiplication.
To illustrate what this amounts to,
consider a quadrilateral in $\SM$ with sides
$\alpha,\beta,\gamma,\delta$ and diagonals $\eta$ and~$\theta$,
as in Proposition~\ref{pr:geometric-flip}.
Let $\baralpha,\barbeta,\bargamma,\bardelta$ be lifts of
$\alpha,\beta,\gamma,\delta$ to the opened surface $\barSM$,
coordinated so that $\baralpha,\barbeta,\bargamma,\bardelta$ still
form a quadrilateral.
(Cf.\ the discussion surrounding~\eqref{eq:p-eta}.)
Define the cross-ratio $\tau$ of the edge~$\varrho$ by
\begin{equation}
  \tau = \frac{\lambda_\varrho(\baralpha)\lambda_\varrho(\bargamma)}
{\lambda_\varrho(\barbeta)\lambda_\varrho(\bardelta)};
\end{equation}
again, the value of~$\tau$ does not depend on the choice of a
quadruple of coordinated lifts
$\baralpha,\barbeta,\bargamma,\bardelta$.
Then
\begin{equation}
  \lambda_{\sigma/\varrho}(\eta) \lambda_{\sigma/\varrho}(\theta)
=
\frac{\tau}{1+\tau} \lambda_{\sigma/\varrho}(\alpha)
\lambda_{\sigma/\varrho}(\gamma)
+
\frac{1}{1+\tau}\lambda_{\sigma/\varrho}(\beta)\lambda_{\sigma/\varrho}(\delta).
\end{equation}
This exchange relation can of course be obtained by dividing the
corresponding Ptolemy relations for the two sets of lambda lengths
(\emph{viz.}, $\lambda_\sigma$~and~$\lambda_\varrho$) by one another.
Furthermore, the coefficients $\frac{\tau}{1+\tau}$ and
$\frac{1}{1+\tau}$ (or their ratio~$\tau$, in the ``$Y$-pattern
version''~\cite{ca4}) transform under flips according to the appropriate
mutation rules, as predicted by
Propositions~\ref{pr:cv-rescaling}--\ref{pr:rescale-normalize}.
(In the case of ordinary flips on $\SM$, this observation was already
made in~\cite{fock-goncharov1, fock-goncharov2,gsv2}.)

Degenerating the hyperbolic structure $\varrho$ to a lift~$\barL$ of a
lamination~$L$ as explained in the first part of this chapter,
we obtain a normalized exchange pattern of geometric type.
More concretely, we get
\[
\frac{\tau}{1+\tau}
=\frac{\lambda_\varrho(\baralpha)\lambda_\varrho(\bargamma)}
{\lambda_\varrho(\bareta)\lambda_\varrho(\bartheta)}
\to
\frac{c_\barL(\baralpha)c_\barL(\bargamma)}
{c_\barL(\bareta)c_\barL(\bartheta)}
\quad \text{and}\quad
\frac{1}{1+\tau}
=\frac{\lambda_\varrho(\barbeta)\lambda_\varrho(\bardelta)}
{\lambda_\varrho(\bareta)\lambda_\varrho(\bartheta)}
\to
\frac{c_\barL(\barbeta)c_\barL(\bardelta)}
{c_\barL(\bareta)c_\barL(\bartheta)},
\]
recovering the coefficients~\eqref{eq:p-eta}.
It is also easy to see directly that these coefficients are obtained
by exponentiating the shear coordinates of~$L$, conforming to our
general recipe.

\chapter{Versions of Teichm\"uller spaces and coordinates}
\label{app:vers-teichmuller}

This appendix provides a brief guide 
to help the reader navigate among various spaces related to 
the Teichm\"uller space which appear throughout this paper, 
and among different coordinate systems on these spaces. 
There are several independent choices involved in these constructions. 
All of them refer to various
structures that can be put on a compact surface~$\Surf$, possibly with
boundary, with a suitable finite set of marked points~$\Mark$; 
these points may lie in the interior (\emph{punctures}) or on the boundary.

The first choice is between \emph{geometric} and \emph{tropical}
coordinates.  Geometric coordinates, for instance the exponentiated length 
of a curve, are coordinates on a Teichm\"uller space; they satisfy
algebraic relations like the Ptolemy
relation~\eqref{eq:ptolemy-lambda}.  Tropical coordinates, on the
other hand, are coordinates on a space of \emph{measured laminations},
which for us amount to collections of non-intersecting curves, with some 
restrictions.  Such a measured lamination~$L$ can be thought of as
defining a degenerate metric for which all contributions to distances 
come from crossing the lamination~$L$.  
Thus tropical coordinates are some variations of intersection numbers.
The relations they satisfy are naturally written in terms of operations
in an appropriate tropical semifield (see Definition~\ref{def:tropical}).

The next choice is between \emph{cusped} and \emph{opened} surfaces.
Geometrically, a hyperbolic cusped surface has an infinitely long
``horn'' of finite area at each of the punctures in the interior of~$\Surf$, 
while an opened
surface may have a geodesic boundary opened up at these same
punctures.  We also need to orient the boundary at these openings, as
in Definition~\ref{def:partial-Teich}.
In the tropical world, laminations on a cusped surface 
avoid the punctures, while
on an opened surface, they are allowed to spiral around
the opened boundary.

We then consider some version of \emph{lambda lengths}
or \emph{shear coordinates}.  Lambda lengths are in principle simpler,
being essentially the exponentiated lengths of arcs connecting marked
points. Shear coordinates, on the other hand, 
depend on an arc in the context of a particular triangulation.  
Lambda lengths are
coordinates on \emph{decorated} Teichm\"uller spaces, involving some
choice for each marked point, while shear coordinates do not need this
choice.

Finally, each set of coordinates can be extended to work with
\emph{tagged} arcs or triangulations; this is necessary to complete
the cluster algebra structure.

The choices for geometric coordinates (ignoring the
tagging) can be summarized by the following diagram:
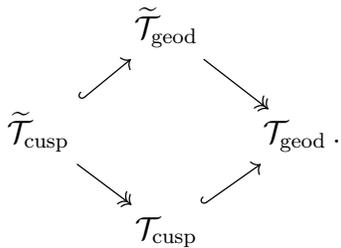
\[
\begin{tikzcd}[>=latex,row sep=scriptsize, column sep=small]
     &\widetilde{\mathcal{T}}_{\operatorname{geod}}
       \arrow[two heads]{dr}\\
  \widetilde{\mathcal{T}}_{\operatorname{cusp}}
       \arrow[hook]{ur} \arrow[two heads]{dr}
  &&   \mathcal{T}_{\operatorname{geod}}\,.\\
  &\mathcal{T}_{\operatorname{cusp}} \arrow[hook]{ur}
\end{tikzcd}
\]
Here $\mathcal{T}_{\operatorname{cusp}}$ is the Teichm\"uller space of
cusped surfaces, while  $\mathcal{T}_{\operatorname{geod}}$ is the
space of opened surfaces with geodesic boundary.  In both cases,
$\widetilde{\mathcal{T}}$ is the decorated Teichm\"uller space, which
allows the definition of lambda lengths, and the coordinates on the
undecorated space $\mathcal{T}$ are shear coordinates.

While each combination of these choices makes sense, we do not explicitly
describe them all in this paper.  Two natural choices are
lambda lengths for
decorated cusped surfaces
($\widetilde{\mathcal{T}}_{\operatorname{cusp}}$ above) and shear
  coordinates for opened surfaces ($\mathcal{T}_{\operatorname{geod}}$ above).
These correspond to the $\mathcal{A}$ and $\mathcal{X}$ spaces of Fock and
Goncharov~\cite{fock-goncharov1}, respectively.  Shear coordinates on
$\mathcal{T}_{\operatorname{cusp}}$ are not independent---they satisfy
a relation for each puncture.  From the
geometric point of view, the principal novelty of this paper is to 
consider lambda lengths for $\widetilde{\mathcal{T}}_{\operatorname{geod}}$
as explained in
Chapter~\ref{sec:lambda-opened}.


\backmatter
\bibliographystyle{amsplain}

\end{document}